\documentclass[11pt]{amsart}
\usepackage[a4paper,margin=1in]{geometry}
\usepackage{amsmath,amssymb,amsthm,mathtools}
\usepackage{enumitem}
\usepackage{xcolor}

\numberwithin{equation}{section}
\newtheorem{theorem}{Theorem}
\newtheorem{lemma}{Lemma}[section]

\newtheorem{remark}{Remark}[section]

\newcommand{\1}{\mathbf 1}

\title[A graph discretization of the Steklov and Neumann problems]{A graph discretization of the Steklov and Neumann problems on domains in complete Riemannian manifolds}
\author{Samuel Audet-Beaumont}
\date{November 11, 2026}

\begin{document}
\begin{abstract}
  We establish quantitative spectral convergence for graph discretizations of the Steklov and Neumann problems for bounded domains in complete Riemannian manifolds. Our construction uses weighted proximity graphs built from maximal $\epsilon$-separated sets, with a correction to the discrete energy accounting for the geometry near the boundary. We consider two discrete Steklov problems: one defined using an averaged boundary trace, and the other using restriction to boundary vertices. Both approximate the continuum Steklov spectrum, but under different conditions on the discretization parameters. The additional scale condition required for the restriction-based problem is also necessary for spectral convergence. As a byproduct, we obtain quantitative convergence of the Neumann spectrum. The estimates are uniform under suitable geometric bounds.
  
\end{abstract}

\maketitle

\section{Introduction}

The goal of this paper is to extend the quantitative graph discretization framework of Burago, Ivanov and Kurylev~\cite{BIK} to the Steklov and Neumann eigenvalue problems on compact domains with smooth boundary in complete Riemannian manifolds. Near the boundary, ambient geodesic balls are truncated by the domain, introducing an anisotropy that must be taken into account in the discrete energy. For the Steklov problem, an additional difficulty arises from the need to approximate the $L^2$-norm on the boundary, rather than in the interior. This leads us to consider two discrete Steklov problems. The first one, using an averaged trace, is more robust. The second one, using restriction to boundary vertices, is more classical but requires additional assumptions on the discretization parameters.

Let $(H,g)$ be a complete smooth Riemannian manifold of dimension
$n\geq2$, and let $M\subset H$ be a compact domain with smooth boundary
$\Sigma=\partial M$. We consider the geometrically controlled class
$\mathcal{M}_n=\mathcal{M}_n(K,D,i_0,r_{tub})$, with bounds on ambient
sectional curvature, on the diameters of $M$ and $\Sigma$, and lower bounds
on the injectivity radii of $H$ and $\Sigma$ and the two-sided normal
injectivity radius of $\Sigma$. The assumptions are specified precisely
in Section~\ref{sec:construction}.

\subsection{Discrete problems and main results}
Choose a maximal $\epsilon$-separated set $B\subset\Sigma$
and extend it to a maximal $\epsilon$-separated set $X\subset M$;
we call $(X,B)$ a maximal $\epsilon$-separated set with boundary.
All distances used for this choice and for the proximity graph are
measured in the ambient manifold $H$. For a second parameter $\rho$,
the graph connects vertices at distance less than $\rho$ and carries
Voronoi volume weights. We associate to it a boundary-corrected energy
$E_\rho^D$ and an averaged trace map
$T_\rho:L^2(X)\to L^2(B)$. Their precise definitions are given in
Section~\ref{sec:construction}. Here $\|\cdot\|_B$ is the norm obtained
from the boundary Voronoi weights.

The generalized discrete Steklov eigenvalues are defined by
the variational principle
$$\sigma_k(X,B):= \min_{\mathcal{L}\subset L^2(X): \text{dim}T_\rho \mathcal{L}=k+1}\max_{u\in \mathcal{L} : T_\rho u \not= 0} R^D(u)$$
with Rayleigh quotient
$$R^D(u):=\frac{E_\rho^D(u)}{\|T_\rho u\|^2_B}.$$

      \begin{theorem}\label{thm1}
            Let $(H,M,g)\in \mathcal{M}_n$ and $(X,B)$ be a maximal $\epsilon$-separated set with boundary. Let $\sigma_k(X,B)$ be defined as before. Then there exist constants $c_k,C_k$ depending only on $\mathcal{M}_n$ and $\sigma_k(M)$ such that if $\rho<c_k$ and $\epsilon< \min\{1/16,1/2n\}\rho,$
            $$|\sigma_k(M)-\sigma_k(X,B)|\leq C_k\left(\sqrt{\rho}+\epsilon+\frac{\epsilon}{\rho}\right).$$
            As a consequence, 
            $$\sigma_k(X,B)\rightarrow \sigma_k(M)$$
            when $\rho,\epsilon,\frac{\epsilon}{\rho}\rightarrow 0.$ If we restrict to the intersection of $\mathcal{M}_n$ with a family of triplets $(H,M,g)$ on which $\sigma_k(M)$ is bounded above, the convergence is uniform.
        \end{theorem}
The averaged trace $T_\rho$ in Theorem~\ref{thm1} is important
for the range of parameters covered there. Indeed, the graph energy is
weighted by $n$-dimensional Voronoi measures, whereas the Steklov
Rayleigh quotient involves $(n-1)$-dimensional boundary measures.
Consequently, functions supported only on boundary vertices may have
small energy compared with the norm of their restriction to $B$.
Averaging over nearby vertices avoids relying solely on those boundary
values and allows convergence under the hypotheses of
Theorem~\ref{thm1}.

On the other hand, the usual Steklov problem on a weighted
graph with boundary is defined using restriction to the boundary
vertices. We therefore also study this problem, replacing the averaged
trace $T_\rho$ with the restriction map
$$R: L^2(X)\rightarrow L^2(B),$$
$$u\mapsto Ru:=u|_B.$$
Its eigenvalues are given by
$$\sigma^N_k(X,B):= \min_{\mathcal{L}\subset L^2(X): \text{dim}R\mathcal{L}=k+1}\max_{u\in \mathcal{L}, Ru \not= 0} R^{D,\text{Nat}}(u)$$
with Rayleigh quotient
$$R^{D,\text{Nat}}(u):=\frac{E_\rho^D(u)}{\| Ru\|^2_B}.$$
For this problem, convergence requires an additional relation
between $\epsilon$ and $\rho$, which severely limits the possible range.
\begin{theorem}\label{thm2}
            Let $(H,M,g)\in \mathcal{M}_n$ and $(X,B)$ be a maximal $\epsilon$-separated set with boundary. Let $\sigma^N_k(X,B)$ be defined as before. Then there exist constants $c_k,C_k$ depending only on $\mathcal{M}_n$ and $\sigma_k(M)$ such that if $\rho^2<c_k\epsilon$ and $\epsilon< \min\{1/16,1/2n\}\rho,$
            $$|\sigma_k(M)-\sigma^N_k(X,B)|\leq C_k\left(\sqrt{\rho}+\epsilon+\frac{\epsilon}{\rho}+\frac{\rho}{\sqrt{\epsilon}}\right).$$
            As a consequence, 
            $$\sigma^N_k(X,B)\rightarrow \sigma_k(M)$$
            when $\rho,\epsilon,\frac{\epsilon}{\rho},\frac{\rho}{\sqrt{\epsilon}}\rightarrow 0.$ The additional assumption $\rho/\sqrt\epsilon\rightarrow 0$ is both sufficient and necessary for spectral convergence of the restriction-based problem. If we restrict our focus to the intersection of $\mathcal{M}_n$ with a family of triplets $(H,M,g)$ on which $\sigma_k(M)$ is bounded above, the convergence is uniform.
        \end{theorem}

Finally, the same discrete energy $E_\rho^D$ also defines a
Neumann graph problem. With the volume-weighted norm $\|\cdot\|_X$
on $L^2(X)$, specified in Section~\ref{sec:construction}, set
$$\lambda_k(X):= \min_{\mathcal{L}\subset L^2(X): \text{dim}\mathcal{L}=k+1}\max_{u\in \mathcal{L} : u \not= 0} R^D_\text{Neu}(u),$$
where
$$R^D_\text{Neu}(u):=\frac{E_\rho^D(u)}{\|u\|^2_X}.$$
The corresponding approximation theorem has the following
form.
\begin{theorem}\label{thm3}
            Let $(H,M,g)\in \mathcal{M}_n$ and $X$ be a maximal $\epsilon$-separated set. Let $\lambda_k(X)$ be defined as before. Then there exist constants $c_k,C_k$ depending only on $\mathcal{M}_n$ and $\lambda_k(M)$ such that if $\rho<c_k$ and $\epsilon< \min\{1/16,1/2n\}\rho,$
            $$|\lambda_k(M)-\lambda_k(X)|\leq C_k\left(\rho+\epsilon+\frac{\epsilon}{\rho}\right).$$
            As a consequence, 
            $$\lambda_k(X)\rightarrow \lambda_k(M)$$
            when $\rho,\epsilon,\frac{\epsilon}{\rho}\rightarrow 0.$ If we restrict our focus to the intersection of $\mathcal{M}_n$ with a family of triplets $(H,M,g)$ on which $\lambda_k(M)$ is bounded above, the convergence is uniform.
        \end{theorem}

        \begin{remark}
          Lu~\cite{JL} extended the quantitative method of Burago, Ivanov
and Kurylev to Neumann eigenvalues on Riemannian manifolds with
boundary, resulting in a result similar to Theorem 3. His construction uses the intrinsic distance and additional
interactions defined through geodesics reflected at the boundary; it
also yields explicit convergence estimates and treats certain spaces
obtained by gluing manifolds. Our construction instead uses distances
in the ambient manifold and corrects the graph weights using distances
between vertices and to the boundary. It avoids the computation of
reflected-geodesic distances and provides a single discrete energy for
both Steklov and Neumann eigenvalues. Lu's Neumann approximation is
already deterministic and quantitative; the distinction here is the
construction and its applicability to both problems, not a general
improvement of the Neumann convergence rate. The geometric assumptions
of the two works are also different.
        \end{remark}
        \begin{remark}
          The construction proposed in the current paper is compatible with the proofs of eigenfunction convergence proved in~\cite{BIK}. We expect such convergence to hold for the three main theorems presented in this section. For the sake of brevity, we do not provide details here.
        \end{remark}

\subsection{Relationship with existing literature}
The comparison of Riemannian eigenvalues with finite-dimensional
models goes back to Dodziuk and Patodi~\cite{Dodziuk}, who used
triangulations. Colin de Verdi\`ere~\cite{ColinDeVerdiere1986}
studied further connections between the spectra of manifolds and graphs.
Mantuano~\cite{Mantuano2005} established uniform two-sided comparisons
between the Laplace eigenvalues of closed Riemannian manifolds and those
of graphs constructed from nets; these estimates do not establish
spectral convergence as the discretization is refined. In contrast,
Aubry~\cite{Aubry2013} obtained quantitative spectral approximation
using geodesic triangulations. Most directly related to the present
construction is the work of Burago, Ivanov and Kurylev~\cite{BIK},
who proved quantitative convergence of Laplace eigenvalues using
weighted proximity graphs on closed manifolds.

The Steklov problem on graphs with boundary has also been studied
independently of approximation questions. Hua, Huang and Wang
\cite{HuaHuangWang2017} obtained Cheeger-type estimates for discrete
Dirichlet-to-Neumann operators. Perrin~\cite{Perrin2019,Perrin2021}
proved geometric lower bounds and isoperimetric upper bounds for
Steklov eigenvalues of graphs. Tschanz~\cite{Tschanz2022,Tschanz2023}
studied Steklov eigenvalue bounds for subgraphs of Cayley graphs and
hyperbolic tilings, making use of comparisons with Riemannian domains.
These works concern geometric inequalities for graph Steklov eigenvalues,
rather than their convergence to the eigenvalues of a prescribed smooth
domain.

For manifolds with boundary, Colbois, Girouard and Raveendran
\cite{CGR} established uniform two-sided comparisons between Steklov
eigenvalues and those of coarse graph discretizations. Their approach
does not involve sending the mesh size to zero, and no spectral
convergence is established. In contrast, the present paper considers
increasingly fine discretizations and proves convergence of individual
Steklov eigenvalues with explicit error estimates. In a different
direction, Vaughn, Berry and Antil~\cite{VaughnBerryAntil2024}
analyzed graph approximations of Dirichlet energy and boundary
integrals, while Peoples and Harlim~\cite{PeoplesHarlim2026} proved
quantitative spectral convergence for Dirichlet and Neumann problems
on sampled embedded manifolds with boundary.

Numerical approximation of Steklov eigenvalues on a given domain
is already a well-developed subject. For example, Andreev and Todorov
\cite{AndreevTodorov2004} analyzed isoparametric finite-element methods,
Armentano and Padra~\cite{ArmentanoPadra2008} established a posteriori
error estimates, Bogosel~\cite{Bogosel2016} used the method of fundamental
solutions, and Wang, Ben-Chen, Polterovich and Solomon
\cite{WangBenChenPolterovichSolomon2018} developed boundary-element
methods for geometry processing. The aim of the present paper is not
to compete with numerical methods tailored to a particular domain,
but to establish approximation estimates for every sufficiently
fine maximal $\epsilon$-separated set, with constants controlled by
the prescribed geometric hypotheses and an upper bound for the
eigenvalue under consideration. This uniformity makes the construction
potentially useful as a theoretical tool in spectral geometry. Aubry's approach
\cite{Aubry2013} shares the emphasis on quantitative control by
geometric hypotheses, but relies on geodesic triangulations of
controlled thinness rather than on separated sets and proximity graphs.

More recently, Lyons~\cite{Lyons} studied discrete Steklov
eigenvalues for independent random samples of Euclidean domains with
$C^2$ boundary and a positive sampling density. The boundary is
represented by a shrinking layer of sampled vertices, and the discrete
Steklov operator is obtained as a Schur complement of the graph
Laplacian. Using transportation and variational methods, Lyons proved
almost-sure convergence of suitably rescaled eigenvalues and results
on eigenvectors and eigenspace projections. In contrast, our vertices
form deterministic $\epsilon$-nets of Riemannian domains and include
vertices on the boundary itself. Theorems~\ref{thm1} and~\ref{thm2}
provide explicit error bounds at finite $\epsilon$ and $\rho$, with
constants controlled by the geometric hypotheses and the corresponding
Steklov eigenvalue. Our approach also approximates Neumann eigenvalues
using the same graph energy. Conversely, Lyons's framework permits
nonconstant sampling densities and establishes eigenvector convergence,
which is not proved here.

Finally, the distinction between our two Steklov problems is
significant: the averaged trace gives convergence for a wider range of
discretization parameters, while restriction to boundary vertices
recovers the usual graph Steklov problem under an additional scale
condition. The use of one boundary-corrected graph energy for these
two Steklov problems and for the Neumann problem is the principal
feature of our approach.

\subsection{Sketch of the proof and plan of the paper}

The proof of Theorem~\ref{thm1} follows the general strategy of Burago, Ivanov and Kurylev~\cite{BIK}, based on comparisons of Rayleigh quotients in both directions. Starting from Steklov eigenfunctions on $M$, we construct discrete functions by averaging over Voronoi cells. We then compare their discrete energy and boundary norms with the corresponding continuous quantities, obtaining an upper bound for the discrete Steklov eigenvalues. For the converse inequality, we construct functions on $M$ from discrete functions and establish analogous comparisons. The two eigenvalue bounds then follow from the min--max principle.

The main new ingredients, compared with the approach of Burago, Ivanov and Kurylev, are the correction of the discrete energy near the boundary and the estimates comparing discrete and continuous boundary norms. Theorem~\ref{thm2} follows from Theorem~\ref{thm1} together with an additional estimate comparing the averaged trace $T_\rho$ with the restriction map $R$. The proof of Theorem~\ref{thm3} follows a parallel argument, using the same discrete energy and replacing the boundary norm comparisons with estimates for the interior $L^2$-norms.

The paper is organized as follows. In Section~\ref{sec:construction}, we give the precise construction of the discrete operators, including the correction of the graph energy near the boundary. Section~3 contains the geometric and analytic preliminary results needed for the approximation estimates. Sections~4 and~5 establish the energy and boundary norm comparisons in both directions, leading to the proof of Theorem~\ref{thm1}. In Section~6, we compare the averaged trace with restriction to boundary vertices and prove Theorem~\ref{thm2}. Finally, Section~7 treats the Neumann problem and establishes Theorem~\ref{thm3}.

\subsection{Acknowledgements}
This project is part of the author's PhD thesis, conducted under the supervision of Alexandre Girouard. The author acknowledges financial support from ISM. The author also thanks Nilima Nigam for interesting discussions on a preliminary version of this paper.

\section{Construction of the discrete operators}\label{sec:construction}

\subsection{Geometric setting, vertex sets, and discrete forms}

We now give the geometric hypotheses and the precise
definitions of the discrete objects appearing in
Theorems~\ref{thm1}--\ref{thm3}.

Let $(H,g)$ be a complete smooth Riemannian manifold of dimension
$n\geq 2$, and let $M\subset H$ be a compact domain with smooth boundary
$\Sigma:=\partial M$. We work within the class
$\mathcal{M}_n=\mathcal{M}_n(K,D,i_0,r_{tub})$ of such triplets
$(H,M,g)$ for which the absolute sectional curvature of $H$ is bounded
by $K$, the diameters of $M$ and $\Sigma$ are at most $D$, the injectivity
radii of $H$ and $\Sigma$ are at least $i_0$, and
$$\text{tub}(\Sigma)\geq r_{tub}.$$
Here $\text{tub}(\Sigma)$ denotes the two-sided normal injectivity
radius of $\Sigma$ in $H$: it is the supremum of the radii $r>0$ for
which the signed normal exponential map
$$\psi : \Sigma \times (-r,r)\rightarrow H,$$
$$(p,t)\mapsto \psi(p,t):=\text{exp}_p(-t\nu_p),$$
is a diffeomorphism onto its image. Throughout the paper,
$\mathcal{M}_n$ is fixed. We use $C$ and $c$ for positive constants
controlled by these geometric data; their values may change from one
occurrence to another.

All distances denoted by $d(\cdot,\cdot)$ are Riemannian
distances in the ambient manifold $H$, unless explicitly indicated
otherwise. For $x\in M$ we write
$$B_r(x):=B^H_r(x)\cap M,$$
where $B_r^H(x)$ is the metric ball in $H$. Thus, in particular,
the balls used to construct the graph are obtained by intersecting
ambient balls with $M$.

Choose parameters $\epsilon<\rho$ satisfying
$$\rho<\min\left\{\frac{i_0}{4}, \frac{r_{tub}}{4}, \frac{1}{2n\sqrt{K}}\right\},\quad \frac{\epsilon}{\rho}<\min\{1/16,1/2n\}.$$
First choose a maximal $\epsilon$-separated set
$B=(x_b)_{b=1}^{N_B}$ on $\Sigma$, using the ambient distance,
and let $(W_b)$ be its boundary Voronoi cells, again defined by the
ambient distance restricted to $\Sigma$. Extend $B$ to a maximal
$\epsilon$-separated set $X=(x_i)_{i=1}^{N}$ in $M$, with interior
Voronoi cells $(V_i)$. We call
$(X,B)$ a maximal $\epsilon$-separated set with boundary.
Note that in general, $W_b\not= V_b\cap \Sigma $ and $W_b$ is not necessarily connected for $x_b\in B$.

We give the cells their natural volume weights:
$$\mu_b:=\text{Vol}_{n-1}(W_b), \quad m_i :=\text{Vol}_n(V_i).$$
As usual, $\sum_b$ and $\sum_i$ denote sums over $B$ and $X$,
respectively. The weights satisfy
$$ \sum_{b} \mu_b =\text{Vol}_{n-1}(\Sigma),\quad  \sum_{i} m_i= \text{Vol}_n(M).$$
We equip the spaces of real-valued functions $L^2(B)$ and
$L^2(X)$ with the weighted inner products
$$\langle u, v\rangle_B :=\sum_b \mu_b u(x_b)v(x_b),\quad \langle u,v\rangle_X :=\sum_i m_i u(x_i)v(x_i),$$
and corresponding norms
$$\| u\|_B^2:= \sum_b \mu_b u(x_b)^2,\quad \|u \|^2_X:= \sum_i m_i u(x_i)^2.$$

The corrected discrete Laplacian on $X$ is defined by
$$(\Delta_D(u))_i:= \frac{4c_\phi }{\rho^{n+2}}\sum_{j:d(x_i,x_j)<\rho} m_j \mathcal{Q}_{ij}^{t(\rho), s}(u_i-u_j).$$
Here, $c_\phi := \frac{n+2}{2\nu_n}$. The weights $\mathcal{Q}_{ij}^{t(\rho),s}>0$ will be made explicit later.
The normalization uses $\nu_n$, the volume of the Euclidean unit
ball. The operator is positive and self-adjoint in $L^2(X)$; its associated
discrete Dirichlet energy is
$$E_\rho^D(u):=\frac{2c_\phi}{\rho^{n+2}} \sum_{ij: d(x_i,x_j)<\rho}m_im_j \mathcal{Q}_{ij}^{t(\rho),s}(u_i-u_j)^2.$$

The first Steklov discretization uses an averaged boundary trace.
For $s>0$, define
$$T_s: L^2(X)\rightarrow L^2(B),$$
$$u\mapsto T_s u,$$
by
$$(T_s u)_b :=\frac{\sum_i k_s(x_b,x_i)m_i u_i
}{\sum_i k_s(x_b,x_i)m_i}  := \hat\Theta_s(b)^{-1} \sum_i k_s(x_b,x_i) m_i u_i,$$
where
$$k_s(x,y):= s^{-n} \phi\left(\frac{d(x,y)}{s}\right),\quad \phi(t):=c_\phi \max\{1-t^2,0\}.$$
Note that $k_s(x,y)\leq Cs^{-n}\mathbf{1}_{B_s(x)}(y)$ .

\subsection{Defining the boundary correcting term $\mathcal{Q}_{ij}^{t(\rho), s}$}
The goal of this subsection is to define $\mathcal{Q}_{ij}^{t(\rho), s}$ while showing some minor results that will be useful for the upcoming proofs. We will also motivate the construction and illustrate why these definitions work particularly well for our purpose.
\\\\
One of the main tools we will use is the quantifiable local flat space comparison afforded by the curvature bounds on the manifold. Thus we must ensure that the construction works as intended on the flat space and use it as a design guide.
Let 
$$\mathbb{H}:= \{\vec{x}\in \mathbb{R}^n : x_n\geq 0\},$$
and 
$$\mathbb{D}_t:=\{x\in B_1(0): x_n\geq -t\}\quad 0\leq t \leq 1.$$
Here, $B_1(0)$ stands for the unit ball in $\mathbb{R}^n$ centered at the origin. For $t\geq 1,$ set $\mathbb{D}_t:=B_1(0).$
For $\vec{x}\in \mathbb{R}^n$ we will write 
$$\vec{x}=(x',x_n),\quad x'\in \mathbb{R}^{n-1}.$$
From this we define the tangential proportion and the normal proportion by
$$\tau(\vec{x}):=\frac{|x'|^2}{|\vec{x}|^2},\quad \nu(\vec{x}):=\frac{|x_n|^2}{|\vec{x}|^2},$$
respectively. From these definitions, it follows that 
\begin{equation}
    \tau(\vec{x})+\nu(\vec{x})=1.
\end{equation}
We can build the flat correction 
$$Q_t(\vec{x}):=a(t)\tau(\vec{x})+b(t)\nu(\vec{x}),$$
where $a,b$ are moment correcting functions. The following will focus on defining those two functions.
    For $t\geq 0,$ define the four universal integrals
    $$A_T(t):=\int_{\mathbb{D}_t}\tau(\vec{x})x_1^2dx=\frac{1}{n-1}\int_{\mathbb{D}_t}\tau(\vec{x})^2|x|^2dx,$$$$ B_T(t):=\int_{\mathbb{D}_t}\nu(\vec{x})x_1^2dx=\frac{1}{n-1}\int_{\mathbb{D}_t}\tau(\vec{x})(1-\tau(\vec{x}))|x|^2dx,$$$$A_N(t):=\int_{\mathbb{D}_t}\tau(\vec{x})x_n^2dx=\int_{\mathbb{D}_t}\tau(\vec{x})(1-\tau(\vec{x}))|x|^2dx,$$$$B_N(t):=\int_{\mathbb{D}_t}\nu(\vec{x})x_n^2dx=\int_{\mathbb{D}_t}(1-\tau(\vec{x}))^2|x|^2dx.$$
    Let 
    $$a_0(t):=\frac{B_N(t)-B_T(t)}{2c_\phi (A_T(t)B_N(t)-A_N(t)B_T(t))},$$
    $$b_0(t):=\frac{A_T(t)-A_N(t)}{2c_\phi (A_T(t)B_N(t)-A_N(t)B_T(t))},$$
    for $0\leq t<1$ and $a_0(t)=b_0(t)=1$ for $t\geq 1$. It follows directly from these definitions that for all $t\geq 0,$
    \begin{equation}
        a_0(t)A_T(t)+b_0(t)B_T(t)=
        a_0(t)A_N(t)+b_0(t)B_N(t)=(2c_\phi)^{-1}.
    \end{equation}
By strict Cauchy-Schwarz,
$$
A_T(t)B_N(t)-A_N(t)B_T(t)\\
$$
\begin{multline}
=\frac{1}{n-1}
\left(\int_{D_t}(1-\tau(\vec{x}))^2|\vec{x}|^2dx\right)
\left(\int_{D_t}\tau(\vec{x})^2\,|\vec{x}|^2dx\right)\\
-\frac{1}{n-1}\left(\int_{D_t}\tau(\vec{x})(1-\tau(\vec{x}))\,|\vec{x}|^2dx\right)^2
>0,
\end{multline}
and thus
$$A_T(t)B_N(t)-A_N(t)B_T(t)\not=0,$$
for all $0\leq t\leq 1$. It follows that $a_0$ and $b_0$ are well defined for all $t\geq 0$. Note also that passing to spherical coordinates and integrating out the tangential angular variables, we can compute that 
$$B_N(t)-B_T(t)>0, \quad A_T(t)-A_N(t)>0.$$
It follows that $a_0(t),b_0(t)>0$. Since $A_T(t),B_T(t),A_N(t),B_N(t)$ are Lipschitz, bounded and continuous on $[0,1]$, and $A_TB_N-A_NB_T>0$ is continuous, it follows that $a_0(t),b_0(t)$ are Lipschitz on $[0,1 ]$ as well. Note that at $t=1$, rotational symmetry gives $a_0(1)=b_0(1)=1.$ Let
$$Q_t^0(\vec{x}):=a_0(t)\tau(\vec{x})+b_0(t)\nu(\vec{x}).$$
The objective of this construction is for the following lemma to hold.
\begin{lemma}
    For every $p\in \mathbb{R}^n$ and $t\geq 0,$
    $$2c_\phi \int_{\mathbb{D}_t}Q_t^0(\vec{x})(p\cdot \vec{x})^2dx=|p|^2.$$
\end{lemma}
 \begin{proof}
     We will use the notation $p=(p',p_n).$ Note that by rotational symmetry in $x'$, 
     $$\int_{\mathbb{D}_t}Q_t^0(\vec{x})x_\alpha x_\beta dx=0\text{ for }\alpha\not=\beta<n,$$
     and
     $$\int_{\mathbb{D}_t}Q_t^0(\vec{x})x_\alpha x_n dx=0\text{ for }\alpha<n,$$
     because the tangential variables are odd. Thus the cross terms have no contribution. We now compute the diagonal terms. For  $\alpha<n$ we have that 
     $$\int_{\mathbb{D}_t}Q_t^0(\vec{x})x_\alpha ^2 dx=a_0(t)A_T(t)+b_0(t)B_T(t)=(2c_\phi )^{-1}.$$
     For the normal direction, we have that 
     $$\int_{\mathbb{D}_t}Q_t^0(\vec{x})x_n ^2 dx=a_0(t)A_N(t)+b_0(t)B_N(t)=(2c_\phi )^{-1}.$$
     Thus 
     $$2c_\phi \int_{\mathbb{D}_t}Q_t^0(\vec{x})(p\cdot \vec{x})^2dx=2c_\phi (2c_\phi)^{-1}|p|^2=|p|^2.$$
 \end{proof}
 This lemma tells us that $a_0$ and $b_0$ yield the right Taylor normalization. However, we will need more than that. We want an abstract dual inequality that will help us greatly when estimating the differentials of smoothed discrete functions. To continue we have to introduce quantities linked to the chosen smoother we will introduce later on. Let
$$\phi(t):=c_\phi \max\{1-t^2,0\},\quad \Theta(t):=\int_{\mathbb{D}_t}\phi(|x|)dx=c_\phi \int_{\mathbb{D}_t}(1-|x|^2)dx.$$
For $0<t<1,$
$$\Theta'(t)= c_\phi \frac{2\nu_{n-1}}{n+1}(1-t^2)^\frac{n+1}{2}.$$
For $t\geq 1,$ we have that 
$$\Theta(t)=1,  \quad \Theta'(t)=0.$$
For $\alpha<n$, we introduce
$$L_\alpha(t,\vec{x})=\frac{2c_\phi}{\Theta(t)}x_\alpha,\quad L_n(t,\vec{x}):=\frac{2c_\phi}{\Theta(t)}x_n- \frac{\Theta'(t)}{\Theta(t)^2}\phi(|x|). $$
These quantities correspond to the derivatives of the smoother in tangential and normal directions, when the smoother is taken in $\mathbb{H}$. Note that we have that 
$$\frac{1}{2}\leq \Theta(t)\leq 1,\quad 0\leq \Theta'(t)\leq C$$
directly from definitions and direct computations. Finally, let 
$$L(t,\vec{x}):=(L_1(t,\vec{x}),L_2(t,\vec{x})\dots L_n(t,\vec{x}))\in \mathbb{R}^n.$$
In what follows, $L(t,\cdot)$ is extended by zero outside of $\mathbb{D}_t$. Since $Q_t^0$ is bounded below away from $0$ and $\Theta,\Theta'$ are bounded and Lipschitz on $[0,1]$, it follows that $L(t,\vec{x})$ is  uniformly bounded for $x\in B_1$. 
We can now define the universal $n\times n$ matrix 
$$G(t):=\int_{\mathbb{D}_t}\frac{L(t,\vec{x})\otimes L(t,\vec{x})}{Q_t^0(\vec{x})}dx.$$
Let 
$$\Lambda(t):=\max\left\{1,\frac{\lambda_{\text{max}}(G(t))}{2c_\phi}\right\}.$$
It follows directly that both $G(t)$ and $\Lambda(t)$ are uniformly bounded and Lipschitz in $t$. Indeed, on $\mathbb{D}_s\cap \mathbb{D}_t$, the integrand in $$|G(t)-G(s)|$$ is uniformly Lipschitz. Furthermore, that same integrand is uniformly bounded on $\mathbb{D}_s\Delta \mathbb{D}_t$ and 
$$|D_s\Delta D_t|\leq C|s-t|.$$
Hence
$$\|G(s)-G(t)\|\leq C|s-t|,\quad |\Lambda(t)-\Lambda(s)|\leq C|s-t|.$$
Finally, set 
$$a(t):=\Lambda(t)a_0(t),\quad b(t):=\Lambda(t)b_0(t)$$
for $0\leq t \leq 1$ and $a(t)=b(t):=1$ for $t\geq 1$. Then, by construction, we have that 
\begin{equation}\label{dualid}
    \int_{\mathbb{D}_t}\frac{L(t,\vec{x})\otimes L(t,\vec{x})}{Q_t(\vec{x})}dx\leq 2c_\phi I
\end{equation}
which is the desired dual inequality. Note that $\Lambda(1)=1$ and that for all $s,t\geq 0$
$$0<c_n\leq a(t),b(t)\leq C_n,\quad |a(t)-a(s)|+|b(t)-b(s)|\leq C_n |t-s|.$$
However, in order to obtain the dual inequality, we pay the price of introducing the new $\Lambda$ term which gives us the following worst Taylor identity (compared to that of $a_0$ and $b_0$) 
\begin{equation}\label{Tayloridentity}
    2c_\phi \int_{\mathbb{D}_t}Q_t(\vec{x})(p\cdot \vec{x})^2dx =\Lambda(t)|p|^2.
\end{equation}
This equation will however prove to be enough for our purpose. Now we have the target flat $Q$-term. We will now define its discrete counterpart, which is the one appearing in our graph Laplacian.
For $x,y\in M$ with $0<d(x,y)<r$, define 
$$\eta_{xy}:=\frac{d(y,\Sigma)-d(x,\Sigma)}{d(x,y)}$$
and $\eta_{xx}=0.$
Note that this quantity is bounded because 
$$|\eta_{xy}|=\left|\frac{d(y,\Sigma)-d(x,\Sigma)}{d(x,y)}\right|\leq \frac{d(x,y)}{d(x,y)}=1.$$
Let $$t=t(r):=\frac{d(x,\Sigma)}{r},$$ and the continuum $Q$-term 
$$\mathcal{Q}_t(x,y)=a\left(t\right)(1-\eta^2_{xy})+b\left(t\right)\eta^2_{xy}.$$
For the graph, if $d(x_i,x_j)<r,$ define
$$\eta_{ij}:=\eta_{x_ix_j}, \quad t_i(r):=\frac{d(x_i,\Sigma)}{r}$$
and 
$$\mathcal{Q}^t_{ij}:=a\left(t_i(r)\right)(1-\eta^2_{ij})+b\left(t_i(r)\right)\eta^2_{ij}.$$
Finally, we define its symmetric counterpart
$$\mathcal{Q}^{t,s}_{ij}:=\frac{1}{2}(\mathcal{Q}^t_{ij}+\mathcal{Q}^t_{ji}).$$
Note that 
$$ \sum_{ij} m_im_j \mathcal{Q}_{ij}^{t,s}(u_i-u_j)^2= \sum_{ij} m_im_j \mathcal{Q}_{ij}^t (u_i-u_j)^2$$
so the discrete energy can be defined either way.
By the boundedness of the functions defining these $Q$ terms,
$$0<c\leq \mathcal{Q}_{ij}^t,\mathcal{Q}_{ij}^{t,s},\mathcal{Q}_t(x,y)\leq C.$$
Note also that all these $Q$-terms are localized near the boundary and are simply equal to $1$ away from the boundary. Their whole purpose is to correct the anisotropy caused by the boundary when integrating over balls.

\section{Preliminary results}
\subsection{Comparative geometry}
First, we introduce a standard curved space comparison result for the tubular neighborhoods of the boundary $\Sigma$.
\begin{lemma}[Uniform collar geometry]\label{collargeom}
There exist constants
$r_*>0$ and  $C_*\geq 1$ 
depending only on $n,K,i_0,r_{\rm tub}$, such that the following properties hold for every $(H,M,g)\in\mathcal M_n$. Let
$$
\psi:\Sigma\times(-r_{\rm tub},r_{\rm tub})\longrightarrow H,
\qquad
\psi(p,t)=\exp_p(-t\nu_p),
$$
be the signed normal exponential map, and set
$$
\Sigma_t:=\psi(\Sigma,t).
$$
For $|t|<r_*$, write $g_t$ for the metric on $\Sigma$ obtained by pulling back the induced metric of $\Sigma_t$ by $\psi(\cdot,t)$. Then
$$
C_*^{-1}g_0\leq g_t\le C_*g_0,
$$
and
$$
(1-C_*|t|)g_0\le g_t\leq(1+C_*|t|)g_0.
$$
If $J(p,t)$ denotes the Jacobian of the Fermi transformation with respect to
$dA_\Sigma\,dt$, then
$$
1-C_*|t|\leq J(p,t)\le1+C_*|t|.
$$
The boundary projection
$$
\pi:\psi(\Sigma\times(-r_*,r_*))\longrightarrow\Sigma,
\qquad
\pi(\psi(p,t))=p,
$$
is $C_*$-Lipschitz. Consequently, for $p,q\in\Sigma$ with
$d_H(p,q)<r_*$,
$$
d_H(p,q)\le d_\Sigma(p,q)\leq C_*\,d_H(p,q).
$$
Moreover, the clamping map
$$
R(\psi(p,t)):=\psi(p,\max\{t,0\})
$$
is $C_*$-Lipschitz on the collar.
\end{lemma}
This result is classical and the reader can refer to~\cite{Tubes}, chapters 2 and 3, for additional details. This result also extends further : we can directly compare the exponential map with model spaces. We will assume from this point forward that $4\rho<r_*$ where $r_*$ is as defined in the previous lemma. Recall our other standing assumptions : $\rho<i_0/4, r_{tub}/4$ and $4K\rho^2<1/n^2$. Fix $x\in M$ . It follows that the Riemannian exponential map 
$$\text{exp}_x:T_xM\rightarrow H,$$
possesses a well-defined inverse $\text{exp}_x^{-1}$ over $B^H_{4\rho}(x)$ the ball of radius $4\rho$ around $x\in M$. It is in fact a diffeomorphism onto its image. For a point $v\in T_x M$ in the image of $B_{4\rho}(x)$ under $\text{exp}^{-1}_x$, set $J_x(v)$ to be the Jacobian of $\text{exp}_x$ at $v$. By the Rauch Comparison Theorem, the relative distortion of metric by $\text{exp}_x$ at $v$ is bounded above and below by a function of $|v|^2$ as a result of the curvature bounds and the scale at which we operate. Thus 
$$(1+C|v|^2)^{-1}\leq J_x(v)\leq 1+C|v|^2.$$
Note that this result also applies to the intrinsic exponential on $\Sigma$ and thus the normal coordinates on Fermi maps at the boundary. This is a central result that we will use multiple times throughout the article, sometimes without explicit mention. 
A direct application of this principle allows us to bound the volume of small balls in the manifold. Indeed, there exist universal constants $c,C>0$ depending only on controlled geometry such that for all $x\in \Sigma$ and all $y\in M$,
$$c s^{n-1}< \text{Vol}_{n-1}(B^\Sigma_s(x))<Cs^{n-1},\quad cs^n< \text{Vol}(B_s(y))<Cs^n,$$
for all $0<s<4\rho.$ We will also need the following result. 
\begin{lemma}\label{scaledfermi}
There exist controlled-geometry constants $\bar{r}_*>0$ and $C_*>0$ such that the following holds. Take $0<r<\bar{r}_*$ and let $x\in M$ satisfy $d(x,\Sigma)<r.$ Let $p=\pi(x)$ be the closest boundary point. Choose an orthonormal basis $e_1,\dots e_{n-1}$ of $T_p\Sigma$. Transport these vectors via parallel transport along the normal geodesic from $p$ to $x$ and complete the basis with the inward unit normal $e_n$. In this way, identify $\mathbb{R}^n$ with $T_x H$. Define
$$\psi_x(z):=\text{exp}_{\text{exp}_p^\Sigma (rz')}\left(-(d(x,\Sigma)+rz_n)\nu_{\text{exp}_p^\Sigma(rz')}\right).$$
Then, uniformly for $z=(z',z_n)$ such that $|z|\leq 2$ and $z_n\geq -d(x,\Sigma)/r$, $\text{exp}_x^{-1}(\psi_x(z))$ is well defined and
$$\left| \frac{1}{r}\text{exp}_x^{-1}(\psi_x(z))-z\right|\leq C_*r|z|,$$
 hence 
$$(1-C_*r)r|z|\leq d(x,\psi_x(z))\leq (1+C_*r)r|z|.$$
\end{lemma}
Indeed with the above identifications, the differential of the scaled Fermi coordinate map differs from parallel transport by $O(r)$, uniformly for $|z|\leq 2$. Combining this Jacobi-field estimate with Rauch comparison theorem gives us that for 
$$F_{x,r}(z):=r^{-1}\text{exp}^{-1}_x (\psi_x(z)),$$
we have 
$$\|DF_{x,r}(z)-I\|\leq Cr.$$
Since $F_{x,r}(0)=0,$ integration along the segment $t\mapsto tz$ gives 
$$|F_{x,r}(z)-z|\leq \int_0^1\|DF_{x,r}(tz)-I\||z|dt\leq Cr|z|.$$
These are standard computations. The reader is invited to consult~\cite{Tubes}, chapters 2–3, and~\cite{Chavel}, chapter IX for further details.\\\\
For the remainder of this article, we will assume that $4\rho<r_*,\bar{r}_*$. Since both quantities are controlled-geometry constants, we are still within the hypothesis of the theorems we wish to prove.
\subsection{Common metric bounds}
In this subsection, we will provide a list of simple metric bounds that will be used throughout the paper, sometimes without stating them explicitly. First, note that by the construction of $X$ and $B$ we directly have that for all $x_b\in B, x_i\in X$, there exist universal constants $c$ and $C$ such that
$$B^\Sigma_{c\epsilon}(x_b)\subset W_b \subset B^\Sigma_{C\epsilon}(x_b), \quad  B_{c\epsilon}(x_i)\subset V_i \subset B_{C\epsilon}(x_i),$$
and thus 
$$c\epsilon^{n-1}\leq \mu_b \leq C\epsilon^{n-1},\quad c\epsilon^n\leq m_i\leq C\epsilon^n.$$
For $V_i$ and $m_i$ this follows directly from the properties of Voronoi sets. For $W_b$ and $\mu_b$ it is a bit more subtle.  Because $B$ is ambient-$\epsilon$-separated, if 
$$d_\Sigma(z,x_b)<\epsilon/3$$
for some $b$, then for every $b'\not=b,$ we have that 
$$d(z,x_{b'})\geq d(x_b,x_{b'})-d(z,x_b)>2\epsilon/3$$ and therefore 
$$B^\Sigma_{\epsilon/3}(x_b)\subset W_b.$$
Furthermore, by the maximality of $B$, for every $z\in \Sigma$ there is a $b'$ such that 
$$d(z,x_{b'})<\epsilon.$$
If $z\in W_b,$ then the Voronoi property gives us that 
$$d(z,x_b)\leq d(z,x_b')<\epsilon.$$ By lemma \ref{collargeom},
$$d_\Sigma(z,x_b)\leq Cd(z,x_{b})\leq C\epsilon,$$
and therefore 
$$W_b\subset B^\Sigma_{C\epsilon}(x_b).$$
As a consequence of this, after possibly decreasing $r_*$ but keeping it universal, there exist controlled-geometry constants $C_0,c,C>0$ such that, for every $z\in \Sigma, x\in M,$ 
\begin{equation}\label{Voronoivolumemu}
c r^{n-1}\leq \sum_{x_b\in B_r^\Sigma(z)}\mu_b \leq C r^{n-1},\end{equation}
for $C_0\epsilon\leq r\leq r_*$ and
\begin{equation}\label{Voronoivolumem}
cr^n\leq \sum_{x_i\in B_r(x)}m_i\leq  Cr^n,\end{equation}
for $2\epsilon\leq r \leq r_*.$
For example, let
$$\hat{\Theta}_\rho(b):=\sum_{i}k_\rho(x_b,x_i)m_i.$$
Then since  $c(1-d(x_b,x_i)^2/\rho^2)\rho^{-n}\mathbf{1}_{B_\rho(x_b)}\leq k_\rho(x_b,x_i)\leq C\rho^{-n}\mathbf{1}_{B_\rho(x_b)},$ 
$$ 0<c\leq \sum_{x_i\in B_\rho(x_b)}c(1-d(x_b,x_i)^2/\rho^2)\rho^{-n}m_i \leq\hat{\Theta}_\rho(b) \leq \sum_{x_i\in B_\rho(x_b)}c\rho^{-n}m_i\leq C.$$
These bounds will be useful later on.
\subsection{Regularity tools}\label{regularitytoolsubs}
The discrete problem is defined so as to obtain the dual bound (\ref{dualid}) at the cost of slightly weakening the Taylor bound, resulting in (\ref{Tayloridentity}). This inexactitude translates further down the road to some distortion in the boundary collar
$$N_{\rho+2\epsilon}=\{x\in M:d(x,\Sigma)<\rho+2\epsilon\}.$$
With enough hypotheses on the functions, this distortion becomes manageable. We will however need the following regularity tools. Let $k\geq 1$ fixed. Consider $f$ such that 
$$f:=\sum_{i=0}^k c_if_i,$$
where the $c_i$ are constants and $f_i$ is an $i$-th Steklov eigenfunction of $M$ orthonormal on $\Sigma$. We have the following uniform classical (see~\cite{Griswald} chapters 1-2) inequalities 
\begin{equation}\label{eq1reg}
\|f\|^2_{L^2(M)}\leq C(\|\nabla f\|^2_{L^2(M)}+\|f|_\Sigma\|^2_{L^2(\Sigma)}),
\end{equation}
\begin{equation}\label{eq2reg}
\|f|_\Sigma\|_{H^{1/2}(\Sigma)}\leq C\|f\|_{H^1(M)},    
\end{equation}
\begin{equation}\label{eq4reg}
\|f\|_{H^{2}(M)}\leq C(\| \Delta f\|_{L^2(M)}+\|\partial_\nu f\|_{H^{1/2}(\Sigma)}+\|f\|_{L^2(M)}).    
\end{equation}
Indeed, by the bounded-geometry assumptions and the uniform collar geometry of Lemma \ref{collargeom}, every $M$ such that $(H,M,g)\in \mathcal{M}_n$ admits boundary charts of radius $r_*$ in which $\Sigma$ is represented by a $C^{1,1}$ graph with uniformly bounded norm. The standard Friedrichs and trace inequalities and the $H^2$ regularity estimate for the Neumann problem therefore hold with constants depending only on the controlled geometry. Since
$$\|\nabla f\|^2_{L^2(M)}=\sum_{i=0}^k \sigma_ic_i^2\leq \sigma_k \|f|_\Sigma\|^2_{L^2(\Sigma)},$$
it follows from (\ref{eq1reg}) that
$$\|f\|^2_{H^1(M)}\leq C(1+\sigma_k)\|f|_\Sigma \|^2_{L^2(\Sigma)}.$$
Applying (\ref{eq2reg}) yields 
\begin{equation}\label{eq3reg}
    \|f|_\Sigma\|^2_{H^{1/2}(\Sigma)}\leq C(1+\sigma_k)\|f|_\Sigma\|^2_{L^2(\Sigma)}.
\end{equation}
Let 
$$F:= \sum_{i=0}^k \sigma_i c_i f_i.$$
Then $F$ is harmonic and 
$$F|_\Sigma = \partial_\nu f.$$
By (\ref{eq3reg}), 
    $$\|\partial_\nu f\|^2_{H^{1/2}(\Sigma)}\leq C(1+\sigma_k)\|\partial_\nu f\|^2_{L^2(\Sigma )}$$
$$=\sum_i \sigma_i^2c_i^2\leq \sigma_k^2\|f|_\Sigma\|^2_{L^2(\Sigma)}.$$
Applying (\ref{eq4reg}) yields 
    \begin{equation}\label{hessianbound}
        \|f\|^2_{H^2(M)}\leq C(1+\sigma_k)^3\|f|_\Sigma \|^2_{L^2(\Sigma)}.
    \end{equation}
We can achieve a similar bound for Neumann. Indeed, assuming $u$ is a linear combination of the $k+1$ first Neumann eigenfunctions, applying the same uniform estimate gives 
\begin{equation}\label{hessianboundneum}
\|u\|^2_{H^2(M)}\leq C(1+\lambda_k)^2\|u\|^2_{L^2(M)}.
\end{equation}
Recall that the moment correcting functions are inexact for the Taylor expansion identity. It will result in an inexact bound for the Dirichlet energy near the boundary. Thus we need to show that it is harmless as long as the inexactitude lives on a tubular neighborhood that is small enough. This is the goal of the following Lemma.
\begin{lemma}\label{bndc1}
   Fix $r_0:= r_*/2$ where $r_*$ is as in Lemma \ref{collargeom}.
   Let $f\in H^2(M)$. Then 
   $$\int_{N_\delta} |\nabla f|^2\leq C\delta \left( \| \nabla f\|^2_{L^2(M)}+\|\nabla^2 f\|^2_{L^2(M)}\right),$$
for $0<\delta < r_0.$
    
\end{lemma}
\begin{proof}
In the Fermi coordinates of Lemma \ref{collargeom}, the corresponding Jacobian satisfies
$$0<C^{-1}\leq J(p,t)\leq C$$
for $0\leq t\leq r_0$.
For fixed $p$, let $\gamma_p(t)=\psi(p,t)$ and let 
$$V_p(t):=P_{t\rightarrow 0} \nabla f(\gamma_p(t)),$$
where $P_{t\rightarrow 0}$ is the parallel transport back to $T_p H$. Since parallel transport is an isometry, we have that 
$$|V_p(t)|=|\nabla f(\gamma_p(t))|.$$
Moreover,
$$V_p'(t)=P_{t\rightarrow 0}(\nabla_{\frac{d\gamma}{dt}} \nabla f),$$
hence 
$$|V'_p(t)|\leq |\nabla^2 f(\gamma_p(t))|.$$
Since $V_p\in H^1(0,r_0),$ for all $t\in (0,r_0)$,
$$V_p(t)=V_p(\tau)+\int_\tau^tV'_p(s)ds,$$
and 
$$|V_p(t)|^2\leq 2|V_p(\tau)|^2+2\left|\int_\tau^t V'_p(s)ds\right|^2$$
$$\leq 2|V_p(\tau)|^2+2|t-\tau|\int_\tau^t |V'_p(s)|^2ds$$
$$\leq 2|V_p(\tau)|^2+2r_0\int_0^{r_0} |V'_p(s)|^2ds.$$
Therefore averaging $\tau$ over $(0,r_0)$ yields 
$$r_0|V_p(t)|^2\leq 2\int_0^{r_0}|V_p(s)|^2ds+2r_0^2\int_0^{r_0}|V'_p(s)|^2ds$$
and it follows that

$$\sup_{0\leq t \leq  r_0} |V_p(t)|^2\leq \frac{2}{r_0}\int_0^{r_0}|V_p(s)|^2ds+2r_0 \int_0^{r_0}|V_p'(s)|^2ds.$$
Therefore 
$$\int_{N_\delta} |\nabla f|^2\leq C \int_\Sigma \int_0^\delta |V_p(t)|^2dt dA(p)$$
$$\leq C\delta \int_\Sigma \left[ \frac{2}{r_0}\int_0^{r_0}|V_p(s)|^2ds+2r_0 \int_0^{r_0}|V_p'(s)|^2ds\right]dA(p)$$
$$\leq C\delta \left( \| \nabla f\|^2_{L^2(M)}+\|\nabla^2 f\|^2_{L^2(M)}\right),$$
for $0<\delta < r_0.$
\end{proof}
\subsection{Discrete energy rescaling bounds}
Our goal here is twofold. First we investigate how the discrete energy behaves under parameter variations. Second, we prove intermediary bounds that will be used in further proofs.
In the following, fix $u\in L^2(X)$ and let 
$$E_\delta (u):=\sum_{i,j:d(x_i,x_j)<\delta}m_im_j(u_i-u_j)^2.$$
We have the following lemmas.
\begin{lemma}\label{sommepascomp}
    For $4\epsilon<\delta\leq \rho$,
    $$E_\delta(u)\leq C\left(\frac{\delta}{\rho}\right)^n E_\rho(u).$$
    Since 
    $$0<c\leq Q_{ij}^{t(\rho)},$$
    it implies that 
    $$E_\delta(u)\leq C\left(\frac{\delta}{\rho}\right)^n\sum_{i,j:d(x_i,x_j)<\rho}m_im_jQ_{ij}^{t(\rho)}(u_i-u_j)^2.$$
\end{lemma}
\begin{proof}
    We do the proof in two cases. First, assume that $\delta \geq \rho/4$. This is the trivial case. Since 
    $$\{i,j: d(x_i,x_j)<\delta\}\subset \{i,j:d(x_i,x_j)<\rho\},$$
    it follows directly that
    $$E_\delta(u)\leq E_\rho(u).$$
    Since $\delta\geq \rho/4,$ then 
    $$\left(\frac{\delta}{\rho}\right)^n\geq 4^{-n}.$$
    Therefore 
    $$E_\delta(u)\leq 4^n\left(\frac{\delta}{\rho}\right)^nE_\rho(u).$$
    This concludes the proof of the first case.
    Now assume that $\delta <\rho/4$. Fix $(i,j)$ such that 
    \begin{equation}\label{ijclose}
        d(x_i,x_j)<\delta.
    \end{equation}
    Define 
    $$A_i:=\{k:d(x_i,x_k)<\rho/2\}. $$
    Since we are at a scale where the Voronoi cells give us volume approximations, by (\ref{Voronoivolumem}) we have that 
    \begin{equation}\label{Aivol}
        \sum_{k\in A_i} m_k\geq C\rho^n.
    \end{equation}
    Note that since
    \begin{equation}\label{jkisanedge}
     d(x_j,x_k)\leq d(x_j,x_i)+d(x_i,x_k)<\rho/4+\rho/2<\rho,  
    \end{equation}
    we have $\rho$-edges $(i,k)$ and $(j,k)$ for $k\in A_i$. Since 
    $$(u_i-u_j)^2\leq 2(u_i-u_k)^2+2(u_k-u_j)^2,$$
    it follows that 
    $$\left(\sum_{k_\in A_i}m_k\right) (u_i-u_j)^2\leq C\sum_{k\in A_i}m_k\left( (u_i-u_k)^2+(u_j-u_k)^2\right).$$
    Then (\ref{Aivol}) implies
    $$(u_i-u_j)^2\leq \frac{C}{\rho^n}\sum_{k\in A_i}m_k\left( (u_i-u_k)^2+(u_j-u_k)^2\right).$$
    Multiplying both sides by $m_im_j$ and summing over all $(i,j)$ such that $d(x_i,x_j)<\delta$ yields 
    \begin{equation}\label{EdeltaleqT1T2}
    E_\delta(u)\leq T_1+T_2
    \end{equation}
    where 
    $$T_1:= \frac{C}{\rho^n}\sum_{i,j:d(x_i,x_j)<\delta}m_im_j \sum_{k\in A_i} m_k(u_i-u_k)^2,$$
    $$T_2:= \frac{C}{\rho^n}\sum_{i,j:d(x_i,x_j)<\delta}m_im_j \sum_{k\in A_i} m_k(u_j-u_k)^2.$$
    Both terms will be estimated separately. We begin with $T_1$. Changing the summation order, by the Voronoi volume comparison (\ref{Voronoivolumem}), we obtain 
    
    $$T_1:= \frac{C}{\rho^n}\sum_i \sum_{k\in A_i} m_im_k(u_i-u_k)^2\left( \sum_{j:d(x_i,x_j)<\delta}m_j\right)$$
    $$\leq C\frac{\delta^n}{\rho^n}\sum_i \sum_{k\in A_i} m_im_k(u_i-u_k)^2$$
    $$= C\frac{\delta^n}{\rho^n}\sum_i \sum_{k:d(x_i,x_k)<\rho/2} m_im_k(u_i-u_k)^2$$
    \begin{equation}\label{T_1bnd}
        \leq C\left(\frac{\delta}{\rho}\right)^n E_\rho(u).
    \end{equation}
    Note that we use the lower bound hypothesis on $\delta$ here in order to be able to use the Voronoi volume comparison.  We will now bound $T_2$. Recall that 
    $$T_2:= \frac{C}{\rho^n}\sum_{i,j:d(x_i,x_j)<\delta}m_im_j \sum_{k:d(x_i,x_k)<\rho/2} m_k(u_j-u_k)^2.$$
    Thus if the triple $(i,j,k)$ contributes to the sum, by (\ref{jkisanedge}), $(j,k)$ is a $\rho$-edge. Changing the order of summation, we obtain
    
    $$T_2:= \frac{C}{\rho^n}\sum_{j,k}m_j m_k(u_j-u_k)^2 \left(\sum_{i:d(x_i,x_k)<\rho/2,\;
                              d(x_i,x_j)<\delta} m_i\right)$$
                              $$=  \frac{C}{\rho^n}\sum_{j,k:d(x_j,x_k)<\rho}m_j m_k(u_j-u_k)^2 \left(\sum_{i:d(x_i,x_k)<\rho/2,\;
                              d(x_i,x_j)<\delta} m_i\right)$$
    $$\leq   \frac{C}{\rho^n}\sum_{j,k:d(x_j,x_k)<\rho}m_j m_k(u_j-u_k)^2 \left(\sum_{i:
     d(x_i,x_j)<\delta} m_i\right)$$
     $$\leq C\frac{\delta^n}{\rho^n}\sum_{d(x_j,x_k)<\rho}m_j m_k(u_j-u_k)^2$$
     \begin{equation}\label{T2bnd}
        =C\left(\frac{\delta}{\rho}\right)^nE_\rho(u).\end{equation}
 Combining (\ref{EdeltaleqT1T2}) with (\ref{T_1bnd}) and (\ref{T2bnd}) yields the desired result.
\end{proof}

\begin{lemma}\label{complik}
    Let 
    $$0<r\leq \rho,\quad \epsilon\leq c\rho.$$
    Then 
    $$\sum_{i,j:d(x_i,x_j)<\rho}\int_{V_i}\int_{V_j}\mathbf{1}_{\{4\epsilon<d(x,y)<r\}}\frac{\epsilon}{d(x,y)}(u_i-u_j)^2dydx \leq C\frac{\epsilon}{\rho}E_\rho(u).$$
    It follows that
    \begin{multline}\sum_{i,j:d(x_i,x_j)<\rho}\int_{V_i}\int_{V_j}\mathbf{1}_{\{4\epsilon<d(x,y)<r\}}\frac{\epsilon}{d(x,y)}(u_i-u_j)^2dydx \\ \leq C\left(\frac{\epsilon}{\rho}\right)\sum_{i,j:d(x_i,x_j)<\rho}m_im_jQ_{ij}^{t(\rho)}(u_i-u_j)^2.        
    \end{multline}
\end{lemma}
    \begin{proof}
        If $r\leq 4\epsilon$, the claim is trivial. Assume $r>4\epsilon$. We decompose the region 
        $$4\epsilon<d(x,y)<r,$$
        in the following way. Let $\delta_\ell:=2^\ell \epsilon$. Choose $L$ such that 
        $$\delta_L \leq r<2\delta_L. $$
        Then, by definition, for fixed $x$
        $$\{y:4\epsilon<d(x,y)<r\}\subset \bigcup_{\ell=2}^L\{y: \delta_\ell <d(x,y)\leq 2 \delta_\ell\}.$$
        Therefore 
        \begin{equation}\label{susmsl}
        \sum_{i,j:d(x_i,x_j)<\rho}\int_{V_i}\int_{V_j}\mathbf{1}_{\{4\epsilon<d(x,y)<r\}}\frac{\epsilon}{d(x,y)}(u_i-u_j)^2dydx \leq \sum_{\ell=2}^L S_\ell,            
        \end{equation}
        where $$S_\ell:=\sum_{i,j:d(x_i,x_j)<\rho}\int_{V_i}\int_{V_j}\mathbf{1}_{\{\delta_\ell<d(x,y)<2\delta_\ell\}}\frac{\epsilon}{d(x,y)}(u_i-u_j)^2dydx.$$
        Since on the annulus 
        $$\{\delta_\ell<d(x,y)\leq 2\delta_{\ell}\},$$
        we have 
        $$\frac{\epsilon}{d(x,y)}\leq \frac{\epsilon}{\delta_\ell},$$
        we obtain 
        \begin{equation}\label{slbnd1}
        S_\ell \leq \frac{\epsilon}{\delta_\ell}\sum_{i,j:d(x_i,x_j)<\rho}\int_{V_i}\int_{V_j}\mathbf{1}_{\{\delta_\ell<d(x,y)<2\delta_\ell\}}(u_i-u_j)^2dydx.            
        \end{equation}
        If $x\in V_i, y\in V_j$ and $d(x,y)\leq 2\delta_\ell,$ we have 
        \begin{equation}\label{premdistxixj}
            d(x_i,x_j)\leq d(x_i,x)+d(x,y)+d(y,x_j)\leq \epsilon+2\delta_\ell+\epsilon.
        \end{equation}
        For $\ell\geq 2,$ 
        $$\delta_\ell=2^\ell\epsilon\geq 4\epsilon.$$
        Therefore
        $$2\epsilon\leq \frac{1}{2}\delta_\ell,$$
        and (\ref{premdistxixj}) yields 
        \begin{equation}\label{xixjnound}
            d(x_i,x_j)\leq 2\delta_\ell+2\epsilon\leq 3\delta_\ell.
        \end{equation}
        Hence the pairs $(i,j)$ that contribute to $S_\ell$ satisfy (\ref{xixjnound}). This bound combined with
        $$\int_{V_i}\int_{V_j}\mathbf{1}_{\{\delta_\ell<d(x,y)<2\delta_\ell\}}(u_i-u_j)^2dydx\leq Cm_im_j (u_i-u_j)^2$$
        applied to (\ref{slbnd1}) yields the following
        \begin{equation}\label{slbnd2}
        S_\ell \leq C\frac{\epsilon}{\delta_\ell}\sum_{i,j:d(x_i,x_j)<\min\{3\delta_\ell,\rho\}}m_im_j(u_i-u_j)^2=C\frac{\epsilon}{\delta_\ell}E_{\min\{3\delta_\ell,\rho\}}(u).\end{equation}
        We now split the proof in two cases. First assume $3\delta_\ell\leq \rho.$ Then lemma \ref{sommepascomp} gives
        $$S_\ell \leq C\frac{\epsilon}{\delta_\ell}E_{3\delta_\ell}(u)\leq C\left(\frac{\epsilon \delta_\ell^{n-1}}{\rho^n}\right)E_\rho(u).$$
        Therefore summing over all $\ell$ such that $3\delta_\ell\leq \rho,$ we obtain 
        $$\sum_{\ell : 3\delta_\ell <\rho}S_\ell \leq C\frac{\epsilon}{\rho^n}\left(\sum_{\ell : 3\delta_\ell <\rho}\delta^{n-1}_\ell\right)E_\rho(u).$$
        Recall that $\delta_\ell=2\ell\epsilon$ and thus the sum is geometric with largest term comparable to $\rho^{n-1}$. It follows that 
        $$\sum_{\ell : 3\delta_\ell <\rho}\delta_\ell^{n-1}\leq C\rho^{n-1}.$$
        Therefore
        \begin{equation}\label{sumlleqrho}
        \sum_{\ell : 3\delta_\ell <\rho}S_\ell \leq C\frac{\epsilon}{\rho^n}\left(\rho^{n-1}\right)E_\rho(u)=C\frac{\epsilon}{\rho}E_\rho(u).
        \end{equation}
        It remains to show what happens with the terms such that $3\delta_\ell>\rho.$ For such $\ell$ we have 
        $$S_\ell\leq C\frac{\epsilon}{\delta_\ell}E_\rho(u)\leq C\frac{\epsilon}{\rho}E_\rho(u),$$
        because
        $$\frac{\epsilon}{\delta_\ell}\leq 3 \frac{\epsilon}{\rho}.$$
        Furthermore, since
        $$\rho/3<\delta_\ell\leq r\leq \rho,$$
        and $\delta_{\ell+1}=2\delta_\ell$, there is at most $3$ such $\ell$. Thus 
        \begin{equation}\label{sumgeqrho}
        \sum_{\ell :3\delta_\ell>\rho}S_\ell\leq C\frac{\epsilon}{\rho}E_\rho(u).
        \end{equation}
        Combining (\ref{sumlleqrho}) and (\ref{sumgeqrho}) yields 
        \begin{equation}\label{fullsumsl}
            \sum_{\ell=2}^LS_\ell\leq C\frac{\epsilon}{\rho}E_\rho(u).
        \end{equation}
        Applying (\ref{fullsumsl}) to (\ref{susmsl}) completes the proof.
    \end{proof}

    \begin{lemma}\label{rhohop}
        Let $A\geq 1$.  Then there exists bounded-geometry constant $C_A>0$ such that if $A\rho< r_*, $
        $$E_{A\rho}(u)\leq C_A E_\rho(u).$$
        Here, $r_*$ is the same quantity present in Lemma \ref{collargeom}.
    \end{lemma}

    \begin{proof}
        Fix a pair $(i,j)$ such that $d(x_i,x_j)<A\rho.$ 
        Let $\gamma_H$ be the ambient minimizing geodesic joining $x_i$ to $x_j$. If $\gamma_H$ leaves $M$, we project its exterior portion onto $\Sigma$ along the normal Fermi coordinates. Since $A\rho<r_*$, by Lemma \ref{collargeom} every point of $\gamma_H$ lies in the tubular neighborhood where this projection is well defined. More precisely, for every s, 
        $$d_H(\gamma_H(s),M)\leq d_H(\gamma_H(s),x_i)\leq d_H(x_i,x_j)<A\rho<r_*.$$
        Let $R$ be the clamping map of Lemma \ref{collargeom} and set 
        $$\gamma:= R\circ \gamma_H.$$
        By Lemma \ref{collargeom} there exist a controlled geometry constant $c^*>0$ such that 
        $$\text{Length}(\gamma)\leq c^*\cdot \text{Length}(\gamma_H)<c^*A\rho.$$
        Let 
        $$N:= \lceil 8c^*A\rceil +1.$$
        Then 
        $$\frac{\text{Length}(\gamma)}{N}<\frac{c^*A\rho}{8c^*A}<\frac{\rho}{8}.$$
        Suppose $\gamma$ is parametrized with constant speed. Let
        $$p_\ell:=\gamma(\ell/N),\quad \ell=0,\dots ,N.$$
        It follows that 
        $$p_0=x_i,\quad p_N=x_j, \quad d(p_{\ell-1},p_\ell)<\frac{\rho}{8}.$$
        For each $1\leq \ell\leq N-1$, define the index set
        $$\mathcal{A}_\ell(i,j):=\{ a: d(x_a,p_\ell)<\rho/8\}.$$
        Let
        $$M_\ell(i,j):= \sum_{a\in A_\ell(i,j)}m_a.$$
        Every $z\in B_{\frac{\rho}{32}}(p_\ell)$ lies in some Voronoi cell $V_a$ whose center satisfies 
        $$d(x_a,p_\ell)\leq d(x_a,z)+d(z,p_\ell)\leq \epsilon+\frac{\rho}{32}<\frac{\rho}{8},$$
        since $\epsilon<\frac{\rho}{16}$. Hence $x_a$ belongs to $\mathcal{A}_\ell$, and therefore 
        $$B_{\rho/32}(p_\ell) \subset \bigcup_{a\in \mathcal{A}_\ell}V_a.$$
        It follows from (\ref{Voronoivolumem}) that
        $$M_\ell(i,j)=\sum_{a\in \mathcal{A}_\ell} m_a\geq \text{Vol}(B_{\rho/32}(p_\ell))\geq c \rho^n$$
        and $M_\ell(i,j)>0$.
        For endpoints, set $a_0=i, a_N=j.$
        For each choice 
        $$a_1\in \mathcal{A}_1(i,j),\dots a_{N-1}\in \mathcal{A}_{n-1}(i,j),$$
        we have that $(a_{\ell-1},a_\ell)$ is a $\rho$-edge on the graph. Indeed, for $1<\ell<N,$
        $$d(x_{a_{\ell-1}},x_{a_\ell})\leq d(x_{a_{\ell-1}},p_{\ell-1})+d(p_{\ell-1},p_\ell)+d(p_\ell,x_{a_\ell})$$
        \begin{equation}\label{disttoAl}
        \leq \rho/8+\rho/8+\rho/8<\rho,\end{equation}
        while for $\ell=1$,
        $$d(x_i,x_{a_1})\leq d(p_0,p_1)+d(p_1,x_{a_1})<\rho/8+\rho/8<\rho,$$
        and for $\ell=N$
        
        $$d(x_j,x_{a_{N-1}})\leq d(p_{N-1},p_N)+d(p_{N-1},x_{a_{N-1}})<\rho/8+\rho/8<\rho.$$
        For any such admissible chain, we have 
        $$(u_i-u_j)^2=\left(\sum_{\ell=1}^Nu_{a_{\ell-1}}-u_{a_\ell}\right)^2\leq N\sum_{\ell=1}^N(u_{a_{\ell-1}}-u_{a_\ell})^2.$$
        Averaging over all such chains, we obtain 
        $$(u_i-u_j)^2\leq N \sum_{\ell=1}^N\sum_{a_1\in \mathcal{A}_1(i,j)}\dots \sum_{a_{N-1}\in \mathcal{A}_{N-1}(i,j)}\left( \prod_{q=1}^{N-1} \frac{m_{a_q}}{M_q(i,j)}\right)(u_{a_{\ell-1}}-u_{a_\ell})^2.$$
        Multiplying by $m_im_j$ and summing over all ordered pairs with $d(x_i,x_j)<A\rho$ yields 
        \begin{equation}\label{earholeqnsumsl}
         E_{A\rho}(u)\leq N\sum_{\ell=1}^N S_\ell,
        \end{equation}
        where
        $$S_\ell:= \sum_{i,j:d(x_i,x_j)<A\rho}m_im_j\sum_{a_1\in \mathcal{A}_1(i,j)}\dots \sum_{a_{N-1}\in \mathcal{A}_{N-1}(i,j)}\left( \prod_{q=1}^{N-1} \frac{m_{a_q}}{M_q(i,j)}\right)(u_{a_{\ell-1}}-u_{a_\ell})^2.$$
        It remains to show that $$S_\ell\leq C_A E_\rho(u),$$
        for every $\ell$. Fix $2\leq \ell\leq N-1.$ Then summing over all vertices except for the $a_\ell$ and $a_{\ell-1}$ terms give $1$ due to the normalization. Therefore
        $$S_\ell =\sum_{i,j:d(x_i,x_j)<A\rho}m_im_j\sum_{p\in \mathcal{A}_{\ell-1}(i,j)}\sum_{q\in \mathcal{A}_{\ell}(i,j)}\frac{m_{p}}{M_{\ell-1}(i,j)}\frac{m_{q}}{M_\ell(i,j)}(u_{q}-u_{p})^2.$$
         Since
         $$M_{\ell-1}(i,j)\geq c\rho^n,\quad M_\ell(i,j)\geq c \rho^n,$$
        we have 
        $$S_\ell \leq C \rho^{-2n}\sum_{i,j:d(x_i,x_j)<A\rho}m_im_j\sum_{p\in \mathcal{A}_{\ell-1}(i,j)}\sum_{q\in \mathcal{A}_{\ell}(i,j)}m_pm_q(u_{q}-u_{p})^2.$$
        $$=C \rho^{-2n}\sum_{pq}m_pm_q (u_{q}-u_{p})^2\sum_{i,j: d(x_i,x_j)<A_\rho, p\in \mathcal{A}_{\ell-1}(i,j),q\in \mathcal{A}_\ell(i,j)}m_im_j$$
        \begin{equation}\label{slbnd3}
        =C \rho^{-2n}\sum_{p,q:d(x_p,x_q)<\rho}m_pm_q (u_{q}-u_{p})^2\sum_{i,j: d(x_i,x_j)<A_\rho, p\in \mathcal{A}_{\ell-1}(i,j),q\in \mathcal{A}_\ell(i,j)}m_im_j.
        \end{equation}
        The last equality comes from the fact that if $p\in \mathcal{A}_{\ell-1}$ and $q\in \mathcal{A}_\ell(i,j)$ then $d(x_p,x_q)<\rho,$ as shown before in equation (\ref{disttoAl}).
        We will bound the inner sum. First, if $p\in \mathcal{A}_{\ell-1}(i,j),$ then 
        $$d(x_p,p_{\ell-1})<\rho/8.$$
        We also have
        $$d(x_i,p_{\ell-1})\leq \text{Length}(\gamma)<CA\rho.$$
        Therefore 
        $$d(x_i,x_p)\leq d(x_i,p_{\ell-1})+d(p_{\ell-1},x_p)<C_A\rho.$$
        Hence for fixed index $p$, if $i$ is such that $p\in A_{\ell-1}(i,j)$ and $d(x_i,x_j)<A\rho$ , then 
        $$x_i\in B_{C_A\rho}(x_p).$$
        The exact same argument holds for $q$ and $j$. If $j$ is such that $q\in A_{\ell}(i,j)$ and $d(x_i,x_j)<A\rho$, then 
        $$x_j\in B_{C_A\rho}(x_q).$$
        Thus 
        \begin{multline}\label{innearsumsl}
         \sum_{i,j: d(x_i,x_j)<A_\rho, p\in \mathcal{A}_{\ell-1}(i,j),q\in \mathcal{A}_\ell(i,j)}m_im_j\\ \leq \sum_{i: d(x_i,x_p)<C_A\rho} m_i \sum_{j: d(x_j,x_q)<C_A\rho}m_j\leq C_A\rho^{2n}.   
        \end{multline}
        Thus applying (\ref{innearsumsl}) to (\ref{slbnd3}) for $2\leq \ell\leq N-1$,
        \begin{equation}\label{midsl}
        S_\ell \leq C_A \sum_{p,q:d(x_p,x_q)<\rho}m_pm_q (u_{q}-u_{p})^2\leq C_A E_\rho(u).    
        \end{equation}
        We now have to show that this remains true for both extremities. Set $\ell=1.$ Then
        $$S_1=\sum_{i,j:d(x_i,x_j)<A\rho}m_im_j \sum_{q\in A_1(i,j)}\frac{m_q}{M_1(i,j)}(u_i-u_q)^2$$
        $$\leq C\rho^{-n}\sum_{i,j:d(x_i,x_j)<A\rho}m_im_j \sum_{q\in A_1(i,j)}m_q(u_i-u_q)^2$$
        $$= C\rho^{-n} \sum_{i,q}m_im_q(u_i-u_q)^2\sum_{j:d(x_i,x_j)<A\rho, q\in \mathcal{A}_1(i,j)}m_j$$
        
        $$= C\rho^{-n} \sum_{i,q: d(x_i,x_q)<\rho}m_im_q(u_i-u_q)^2\sum_{j:d(x_i,x_j)<A\rho, q\in \mathcal{A}_1(i,j)}m_j.$$
        Here, the inequality follows from the lower bound on $M_1(i,j)$. The last equality follows from the fact that if $q\in \mathcal{A}_1(i,j)$, then $d(x_i,x_q)<\rho$ as shown before in equation (\ref{disttoAl}). Since 
        $$\sum_{j:d(x_i,x_j)<A\rho, q\in \mathcal{A}_1(i,j)}m_j\leq \sum_{j:d(x_i,x_j)<A\rho}m_j\leq C_A\rho^n,$$
        we obtain 
        \begin{equation}\label{S_1bnd111}
        S_1\leq C_A\sum_{i,q: d(x_i,x_q)<\rho}m_im_q(u_i-u_q)^2\leq C_AE_\rho(u).    
        \end{equation}
        Finally, we need the bound for $S_N$. This case is very similar. We have 
        $$S_N =\sum_{i,j:d(x_i,x_j)<A\rho}m_im_j \sum_{p\in \mathcal{A}_{N-1}(i,j)}\frac{m_p}{M_{N-1}(i,j)}(u_p-u_j)^2$$
        $$\leq C\rho^{-n}\sum_{i,j:d(x_i,x_j)<A\rho}m_im_j \sum_{p\in \mathcal{A}_{N-1}(i,j)}m_p(u_p-u_j)^2$$
        $$=C\rho^{-n} \sum_{p,j} m_pm_j(u_p-u_j)^2\sum_{i:d(x_i,x_j)<A\rho, p\in \mathcal{A}_{N-1}(i,j)}m_i$$
        $$=C\rho^{-n} \sum_{p,j:d(x_p,x_j)<\rho} m_pm_j(u_p-u_j)^2\sum_{i:d(x_i,x_j)<A\rho, p\in \mathcal{A}_{N-1}(i,j)}m_i.$$
        Here, the inequality follows from the lower bound on $M_{N-1}(i,j).$ The last equality follows from the fact that if $p\in \mathcal{A}_{N-1}(i,j)$, then $d(x_j,x_p)<\rho$ as shown before in equation (\ref{disttoAl}). Since 
        $$\sum_{i:d(x_i,x_j)<A\rho, p\in \mathcal{A}_{N-1}(i,j)}m_i\leq \sum_{i:d(x_i,x_j)<A\rho}m_i\leq C_A\rho^n,$$
        we obtain 
        \begin{equation}\label{S_Nbnd}
        S_{N}\leq C_A \sum_{p,j:d(x_p,x_j)<\rho} m_pm_j(u_p-u_j)^2\leq C_A E_\rho(u).   
        \end{equation}
        It follows from (\ref{midsl}), (\ref{S_1bnd111}) and (\ref{S_Nbnd}) that for $1\leq \ell\leq N$,
        \begin{equation}\label{allsl}
        S_\ell \leq C_A E_\rho(u).    
        \end{equation}
        Combining (\ref{earholeqnsumsl}) with (\ref{allsl}) completes the proof.
    \end{proof}

\section{Upper bound on $\sigma_k(X,B)$ in terms of $\sigma_k(M)$ : From continuous to discrete}

The goal of this section is to obtain an upper bound on $\sigma_k(X,B)$ in terms of $\sigma_k(M)$. The way we obtain such a bound is by using the variational characterization of $\sigma_k(X,B)$. Therefore, we need to have discrete test functions for $\sigma_k(X,B)$ that somehow related to $\sigma_k(M)$. A natural choice is therefore to build these test functions from smooth Steklov eigenfunctions. The first step is thus to describe how we will build the discrete analogues of smooth functions.
\\\\
Fix $f\in C^2(\overline{M})$. Let
$$P : C^2(\overline{M})\rightarrow L^2(X),$$
$$f\mapsto Pf,$$
where for $i\in I$,
$$(Pf)_i:= \frac{1}{m_i}\int_{V_i}fdV_n.$$
We split this section in three distinct parts. First will be the energy estimates, where we will compare the discrete energy of $Pf$ with the usual Dirichlet energy of $f$. Second, we will compare the $B$-norm of the discrete trace of $Pf$ with the $L^2(\Sigma)$-norm of the classical trace of $f$. In the third subsection, we will specialize those results to functions which are a finite sum of Steklov eigenfunctions of $M$. We will then bound the Rayleigh quotient of such functions to obtain one side of the final eigenvalue bound. 
\subsection{Discrete energy estimates of the discretized smooth functions}
In the discrete Rayleigh quotient, the discrete energy appears in the numerator. Thus if we take $Pf$ as a test function, we need to bound its discrete energy in terms of its usual Dirichlet energy. This subsection is concerned with exactly those estimates. For general $f\in C^2(\overline{M})$, such uniform estimates can be difficult to obtain. Therefore, we do the process in two ways. We show weaker bounds, where we bound $E_\rho^D(Pf)$ in terms of the Dirichlet energy of $Pf$ and $L^2$-hessian norm of $f$ and allow for some near-boundary distortion in the Dirichlet energy . Since these bounds will be applied to a linear combination of smooth eigenfunctions, the regularity tools introduced in the previous section will allow us to bridge this gap for this very specific case down the line and obtain pure Dirichlet energy bounds.
\\\\
By definition, 
$$(Pf)_i-(Pf)_j=\frac{1}{m_im_j}\int_{V_i}\int_{V_j}(f(x)-f(y))dxdy.$$
Therefore, by Cauchy-Schwarz
$$((Pf)_i-(Pf)_j)^2\leq\frac{1}{m_im_j}\int_{V_i}\int_{V_j}(f(x)-f(y))^2dxdy.$$
It follows that 
$$E^D_\rho(Pf)=\frac{2c_\phi}{\rho^{n+2}}\sum_{i,j:d(x_i,x_j)<\rho}m_im_j\mathcal{Q}^{t(\rho),s}_{ij}((Pf)_i-(Pf)_j)^2$$
\begin{equation}\label{eqerhotointerpol}
\leq\frac{2c_\phi}{\rho^{n+2}}\sum_{i,j:d(x_i,x_j)<\rho}\mathcal{Q}^{t(\rho)}_{ij}\int_{V_i}\int_{V_j}(f(x)-f(y))^2dxdy. 
\end{equation}
To get closer to the continuous setting, we want to bound $\mathcal{Q}^{t(\rho)}_{ij}$ by its continuous counterpart. In order to do so we will need the following pointwise bound. 
\begin{lemma}\label{swapQdiscrtocont}
Let $R>0$ such that $|R-\rho|\leq 4\epsilon.$ Suppose $x\in V_i$ and $y\in V_j$ with $d(x_i,x_j)<\rho$ and  $4\epsilon<d(x,y)<R$.  Then
    $$|\mathcal{Q}^{t(\rho)}_{ij}-\mathcal{Q}_{t(R)}(x,y)|\leq C\frac{\epsilon}{\rho}+C\frac{\epsilon}{d(x,y)}.$$
\end{lemma}
\begin{proof}
This is a straightforward Lipschitz bound. First recall that
$$\eta_{xy}:=\frac{d(y,\Sigma)-d(x,\Sigma)}{d(x,y)}, \eta_{xx}:=0,$$ 
$$\mathcal{Q}_{t(r)}(x,y)=a\left(\frac{d(x,\Sigma)}{r}\right)(1-\eta^2_{xy})+b\left(\frac{d(x,\Sigma)}{r}\right)\eta^2_{xy},$$
and
$$\mathcal{Q}^{t(r)}_{ij}:=a\left(\frac{d(x_i,\Sigma)}{r}\right)(1-\eta^2_{x_ix_j})+b\left(\frac{d(x_i,\Sigma)}{r}\right)\eta^2_{x_ix_j}.$$
Note that 
$$|d(x_i,x_j)-d(x,y)| \leq 2 \epsilon,$$
because the voronoi cells are contained in balls of radius $\epsilon$. It follows that
$$d(x_i,x_j)\geq d(x,y)-2\epsilon\geq \frac{d(x,y)}{2}.$$
Since $d(\cdot ,\Sigma)$ is $1$-Lipschitz, we have that
$$|d(x_i,\Sigma)-d(x,\Sigma)|\leq \epsilon,\quad |d(x_j,\Sigma)-d(y,\Sigma)|\leq \epsilon.$$
We start by bounding 
$$\left|a\left(\frac{d(x_i,\Sigma)}{\rho}\right)-a\left(\frac{d(x,\Sigma)}{R}\right) \right|,\quad  \left|b\left(\frac{d(x_i,\Sigma)}{\rho}\right)-b\left(\frac{d(x,\Sigma)}{R}\right) \right|.$$
First note that since $a(t)=b(t)=1$ for $t\geq 1$, it follows that 
$$a\left(\frac{d(z,\Sigma)}{s}\right)=a\left(\min\left\{\frac{d(z,\Sigma)}{s},1\right\}\right),\quad b\left(\frac{d(z,\Sigma)}{s}\right)=b\left(\min\left\{\frac{d(z,\Sigma)}{s},1\right\}\right).$$
Considering all different cases, we have that
$$\left| \min\left\{\frac{d(x_i,\Sigma)}{\rho},1\right\}-\min\left\{\frac{d(x,\Sigma)}{\rho},1\right\}\right|\leq \frac{\epsilon}{\rho},$$
$$\left| \min\left\{\frac{d(x,\Sigma)}{\rho},1\right\}-\min\left\{\frac{d(x,\Sigma)}{R},1\right\}\right|\leq C\frac{|R-\rho|}{\rho}\leq C\frac{\epsilon}{\rho},$$
and therefore
$$\left| \min\left\{\frac{d(x_i,\Sigma)}{\rho},1\right\}-\min\left\{\frac{d(x,\Sigma)}{R},1\right\}\right|\leq C\frac{\epsilon}{\rho}.$$
Since $a$ and $b$ are Lipschitz, it follows directly that 
$$\left|a\left(\frac{d(x_i,\Sigma)}{\rho}\right)-a\left(\frac{d(x,\Sigma)}{R}\right) \right|\leq C\frac{\epsilon}{\rho},\quad \left|b\left(\frac{d(x_i,\Sigma)}{\rho}\right)-b\left(\frac{d(x,\Sigma)}{R}\right) \right|\leq C\frac{\epsilon}{\rho}.$$
We now bound the $\eta$ terms
$$|\eta_{x_ix_j}-\eta_{xy}|=\left|\frac{d(x_j,\Sigma)-d(x_i,\Sigma)}{d(x_i,x_j)}-\frac{d(y,\Sigma)-d(x,\Sigma)}{d(x,y)}\right|$$
$$\leq \left|\frac{(d(x_j,\Sigma)-d(x_i,\Sigma))-(d(y,\Sigma)-d(x,\Sigma))}{d(x_i,x_j)}\right|+|d(y,\Sigma)-d(x,\Sigma)|\left| \frac{1}{d(x_i,x_j)}-\frac{1}{d(x,y)}\right|.$$
For the first term, we can bound the numerator by 
$$|(d(x_j,\Sigma)-d(x_i,\Sigma))-(d(y,\Sigma)-d(x,\Sigma))|$$$$\leq |d(x_j,\Sigma)-d(y,\Sigma)|+|d(x_i,\Sigma)-d(x,\Sigma)|\leq C\epsilon,$$
and the denominator by 
$$d(x_i,x_j)\geq \frac{d(x,y)}{2}.$$
The second term can be bounded similarly by
$$|d(y,\Sigma)-d(x,\Sigma)|\left| \frac{|d(x,y)-d(x_i,x_j)|}{d(x_i,x_j)d(x,y)}\right|$$$$\leq
|d(y,\Sigma)-d(x,\Sigma)|\frac{2\epsilon}{d(x_i,x_j)d(x,y)}\leq \frac{2\epsilon}{d(x_i,x_j)}\leq C\frac{\epsilon}{d(x,y)}.$$
It follows that 
$$|\eta_{x_ix_j}-\eta_{xy}|\leq C\frac{\epsilon}{d(x,y)}.$$
Since both $\eta$ terms are in $[-1,1]$, we have 
$$|\eta^2_{x_ix_j}-\eta^2_{xy}|\leq C\frac{\epsilon}{d(x,y)}.$$
Since $a$ and $\eta$ are bounded,
$$\left |a\left(\frac{d(x,\Sigma)}{R}\right)\eta^2_{xy}-a\left(\frac{d(x_i,\Sigma)}{\rho}\right)\eta^2_{x_ix_j}\right|$$
$$\leq \left|a\left(\frac{d(x,\Sigma)}{R}\right)-a\left(\frac{d(x_i,\Sigma)}{\rho}\right)\right|\eta_{xy}^2+a\left(\frac{d(x_i,\Sigma)}{\rho}\right)|\eta^2_{xy}-\eta^2_{x_ix_j}|$$
$$\leq C\frac{\epsilon}{\rho}+C\frac{\epsilon}{d(x,y)}.$$
The same holds for $b$. Consequently,
$$|\mathcal{Q}^{t(\rho)}_{ij}-\mathcal{Q}_{t(R)}(x,y)|\leq \left|a\left(\frac{d(x_i,\Sigma)}{\rho}\right)-a\left(\frac{d(x,\Sigma)}{R}\right) \right|$$$$+\left|a\left(\frac{d(x,\Sigma)}{R}\right)\eta^2_{xy}-a\left(\frac{d(x_i,\Sigma)}{\rho}\right)\eta^2_{x_ix_j}\right|$$$$+\left|b\left(\frac{d(x,\Sigma)}{R}\right)\eta^2_{xy}-b\left(\frac{d(x_i,\Sigma)}{\rho}\right)\eta^2_{x_ix_j}\right|$$
$$\leq C\frac{\epsilon}{\rho}+C\frac{\epsilon}{d(x,y)}.$$
\end{proof}
We must now link the continuous interpolation energy with the usual Dirichlet energy. Furthermore, we need to bound effectively the error terms introduced by the pointwise bounds. We proceed by introducing the following technical lemma.

\begin{lemma}\label{directcomparison}
Let $f\in C^2 (\overline M), 0<r\leq 2\rho.$ Then 
\begin{equation}\label{asser1}
\frac{1}{r^2}\int_M \int_{B_r(x)} (f(y)-f(x))^2dydx
\leq Cr^n\|\nabla f\|^2_{L^2(M)}+Cr^{n+2} \|f\|^2_{H^2(M)},
\end{equation}

\begin{equation}\label{asser2}
    \frac{1}{r^2}\int_M \int_{B_r(x)} (d(x,y))^{-1}(f(y)-f(x))^2dydx \leq Cr^{n-1}\|\nabla f\|^2_{L^2(M)}+Cr^{n+1} \|f\|^2_{H^2(M)},
\end{equation}
    \begin{multline}\label{asser3}
        \frac{1}{r^2}\int_M \int_{B_r(x)}\mathcal{Q}_{t(r)}(x,y)(f(y)-f(x))^2dydx\\ \leq \frac{r^{n}}{2 c_\phi}(1+Cr) \| \nabla f \|^2_{L^2(\Lambda(r)}+Cr^{n+1}\| f \|^2_{H^2(M)}.
    \end{multline}
  \\\\  where 
    $$\|\nabla f\|^2_{L^2(\Lambda(r))}:=\int_M \Lambda\left(\frac{d(x,\Sigma)}{r}\right)|\nabla f(x)|^2dx.$$
\end{lemma}

\begin{proof}
The direct comparison is based on the Taylor expansion of $f$. In order to do so, we will need to extend $f$ slightly beyond the boundary of $M$. Let 
$$U:=\{ x\in H \; :\; d_H(x,M)<r\}.$$
By the curvature and injectivity-radius assumptions, the ambient manifolds admit harmonic coordinate charts of a fixed controlled radius in which the metric is uniformly $C^{1,\alpha}$-controlled (see Theorem 1.1 of Anderson's work ~\cite{Anderson1990}). Moreover, the tubular-radius lower bound, together with the sectional-curvature bound, gives a uniform bound on the second fundamental form of $\Sigma$ (see Theorem A of Guijarro and Wilhelm's work~\cite{Guijarro}); hence, after decreasing the chart radius if necessary, $M$ admits boundary charts of radius $r_*$ in which $\Sigma$ is represented by a $C^{1,1}$ graph with uniformly controlled norm. The Sobolev extension theorem for Lipschitz domains (see~\cite{Stein}, Chapter 6 \S3) therefore gives, after localization and a uniformly locally finite partition of unity, a bounded extension operator $$E:H^2(M)\rightarrow H^2(U),$$
such that $Ef(x)=f(x)$ for all $x\in M$ and satisfying
\begin{equation}\label{ExtHess}
  \| Ef\|^2_{H^2(U)}\leq C\| f\|^2_{H^2(M)}.
\end{equation}

Let $x\in M, y\in B_r(x)$. Let $\gamma(s):= \text{exp}_x(s\xi),s\in[0,1],$ where $\xi:=\text{exp}^{-1}_x(y)$. Note that $\gamma$ is in $U$ because 
$$d(\gamma(s),M)\leq d(\gamma(s),x)=sd(x,y)<r.$$ In particular, $\gamma$ is a geodesic from $x$ to $y$. By density, it suffices to prove the following estimates for smooth functions. We will therefore assume $Ef$ to be smooth and obtain estimates depending only on $H^2(U)$ norm, so that the result holds in $H^2(U)$ after taking the limits of approximating smooth functions. Let
$$F(s):=Ef(\gamma(s)).$$
By the chain rule and the Levi-Civita connection's definition, we have that
$$F'(s):=\langle \nabla Ef(\gamma(s)),\gamma'(s)\rangle.$$
Differentiating once again yields 
$$F''(s)=\frac{d}{ds}\langle \nabla Ef(\gamma(s)),\gamma'(s)\rangle$$
$$=\langle \nabla_{\gamma'(s)}\nabla Ef(\gamma(s)),\gamma'(s)\rangle+\langle \nabla Ef(\gamma(s)),\nabla_{\gamma'(s)}\gamma'(s)\rangle$$
$$=\langle \nabla_{\gamma'(s)}\nabla Ef(\gamma(s)),\gamma'(s)\rangle$$
$$=\nabla^2 Ef_{\gamma(s)}(\gamma'(s),\gamma'(s)).$$
The cancellation in the third equality follows from the fact that $\gamma$ is a geodesic. The last equality follows from the definition of the Hessian. To keep the dependencies clear, we will write 
$$\gamma'(s)=P_s(\xi)$$
where $P_s:T_xM\rightarrow T_{\gamma(s)}H$ is the parallel transport along $\gamma$. In particular, note that $|P_s \xi|=|\xi|=d(x,y).$
\\\\
In order to exhibit the linear nature of the Hessian operator, we denote by
$$\xi^t\nabla^2Ef(\gamma(s))\xi:=(\nabla^2Ef)_{\gamma(s)}(P_s\xi,P_s\xi).$$
For a twice differentiable function $g$, we have that 
\begin{equation}\label{identity1}
g(1)-g(0)=g'(0)+\int_0^1(1-s)g''(s)ds.
\end{equation}
Applying (\ref{identity1}) to $F$ and dividing by $r$ yields 
$$\frac{Ef(\text{exp}_x(\xi))-Ef(x)}{r}=\frac{\langle \nabla Ef(x),\xi\rangle}{r}+\frac{1}{r}\int_0^1(1-s)\xi^t\nabla^2Ef(\gamma(s))\xi ds$$
$$=: L(x,\text{exp}_x(\xi))+R(x,\text{exp}_x(\xi)).$$
and thus for all $\delta>0,$ 
\begin{equation}\label{deltabnd}
\left(\frac{Ef(\text{exp}_x(\xi))-Ef(x)}{r}\right)^2\leq (1+\delta)L(x,\text{exp}_x(\xi))^2+\left(1+\frac{1}{\delta }\right)R(x,\text{exp}_x(\xi))^2.
\end{equation}
We start by introducing a first change of variables in $y$ by exponential coordinates around $x$. We will note $\xi_n$ the component of $\xi$ in the direction of the boundary when $d(x,\Sigma)\leq r$. In other contexts, $\xi_n$ can be any direction. We have that 
$$R(x,\text{exp}_x(\xi))^2\leq \frac{1}{r^2}\left(\int_0^1(1-s)\xi^t\nabla^2Ef(\gamma(s))\xi ds\right)^2$$
$$\leq \frac{1}{r^2}\left(\int_0^1 (1-s)ds\right)\left(\int_0^1(1-s)(\xi^t\nabla^2Ef(\gamma(s))\xi )^2ds\right)$$
$$=\frac{1}{r^2}\int_0^1(1-s)(\xi^t\nabla^2Ef(\gamma(s))\xi )^2ds$$
$$\leq \frac{1}{r^2}\int_0^1|\xi|^4\|\nabla^2Ef(\gamma(s))\|^2_Fds$$
\begin{equation}\label{identity2}
\leq \frac{1}{r^2}\int_0^1|\xi|^4\|\nabla^2Ef(\gamma(s))\|^2_Fds.
\end{equation}
The first inequality follows from Cauchy-Schwarz. Here, $\|\cdot\|_F$ denotes the Frobenius operator norm. We now introduce another change of variable. Let $z:=\text{exp}_x(s\xi)$. It follows that $|\xi|=\frac{d(x,z)}{s}$. Thus integrating (\ref{identity2}) over an $r$ ball and using the known Jacobian bounds, we have that
$$\int_{B_r(x)}R(x,y)^2dy\leq (1+Cr)\frac{1}{r^2}\int_0^1\int_{d(x,z)\leq sr}\frac{d(x,z)^4}{s^{n+4}}\|\nabla^2Ef(z)\|^2_F dz ds$$
$$\leq (1+Cr)\frac{1}{r^2}\int_0^1\int_{U}\mathbf{1}_{\{d(x,z)\leq sr\}}\frac{d(x,z)^4}{s^{n+4}}\|\nabla^2Ef(z)\|^2_F dz ds.$$
We change into geodesic coordinates around $z$. Let $v:=\text{exp}_z^{-1}(x).$ We use once again Jacobian bounds to obtain 
$$\int_M\int_{B_r(x)}R(x,y)^2dydx\leq (1+Cr)\frac{1}{r^2}\int_0^1\int_U\int_{|v|\leq sr}\frac{|v|^4}{s^{n+4}}\|\nabla^2Ef(z)\|^2_F dv dzds.$$

$$\leq C(1+Cr^2)r^{n+2}\int_0^1\int_U\|\nabla^2Ef(z)\|^2_F dxds$$
$$\leq  Cr^{n+2} \|\nabla^2 Ef\|^2_{L^2(U)}$$
\begin{equation}\label{intMintbrRcarré}
\leq Cr^{n+2}\|f\|^2_{H^2(M)}.
\end{equation}
We proceed to bound $L$. Note that since $x\in M$, and $Ef=f$ on $M$, we have 
$$\nabla_H Ef(x)=\nabla f (x).$$
In geodesic normal coordinates at $x$, with $\xi:=\text{exp}_x^{-1}(y)$, we have that 
$$\int_M\int_{B_r(x)}L(x,y)^2dydx\leq (1+Cr)\frac{1}{r^{2}}\int_M\int_{|\xi|\leq r}\langle \xi,\nabla f(x)\rangle^2 d\xi dx$$
$$=(1+Cr)\left(\frac{1}{r^{2}}\int_M \sum_{i,j}\partial_if(x)\partial_jf(x)\int_{|\xi|\leq r}\xi_i\xi_jd\xi dx\right)$$
\begin{equation}\label{intMintbrLcarré}
\leq Cr^{n}\|\nabla f\|^2_{L^2(M)}.
\end{equation}
From (\ref{deltabnd}) with $\delta=2$ combined with (\ref{intMintbrRcarré}) and (\ref{intMintbrLcarré}),  $$\frac{1}{r^2}\int_M \int_{B_r(x)} (f(y)-f(x))^2dydx\leq C\int_M\int_{B_r(x)} L(x,y)^2dxdy+C\int_M\int_{B_r(x)} R(x,y)^2dxdy  $$
$$\leq Cr^n\|\nabla f\|^2_{L^2(M)}+Cr^{n+2} \|f\|^2_{H^2(M)}.$$
This completes the proof of (\ref{asser1}). We now prove (\ref{asser2}). We have  
$$\frac{1}{r^2}\int_M \int_{B_r(x)} (d(x,y))^{-1}(f(y)-f(x))^2dydx\leq 2\int_M\int_{B_r(x)} (d(x,y))^{-1}L(x,y)^2dxdy$$$$+2\int_M\int_{B_r(x)}(d(x,y))^{-1} R(x,y)^2dxdy.$$
So it remains to bound the left term precisely. 
By the same proof as before, 
$$ \int_M\int_{B_r(x)}(d(x,y))^{-1}R(x,y)^2dydx\leq (1+Cr)\frac{1}{r^2}\int_0^1\int_U\int_{|v|\leq sr}\frac{|v|^3}{s^{n+3}}\|\nabla^2Ef(z)\|^2_F dv dzds.$$

$$\leq C(1+Cr^2)r^{n+1}\int_0^1\int_U\|\nabla^2Ef(z)\|^2_F dxds$$
$$\leq  Cr^{n+1} \|\nabla^2Ef\|^2_{L^2(U)}$$
$$\leq Cr^{n+1} \|f\|^2_{H^2(M)}$$
and 
$$\int_M\int_{B_r(x)}(d(x,y))^{-1}L(x,y)^2dydx\leq (1+Cr)\frac{1}{r^{2}}\int_M\int_{|\xi|\leq r}\frac{\langle  \xi,\nabla f(x)\rangle^2}{|\xi|} d\xi dx$$
$$=(1+Cr)\left(\frac{1}{r^{2}}\int_M \sum_{i,j}\partial_if(x)\partial_jf(x)\int_{|\xi|\leq r}\frac{\xi_i\xi_j}{|\xi|}d\xi dx\right)$$
$$\leq Cr^{n-1}\|\nabla f\|^2_{L^2(M)}.$$
Applying (\ref{deltabnd}) with $\delta=2$ once again, it follows that 

$$\frac{1}{r^2}\int_M \int_{B_r(x)} (d(x,y))^{-1}(f(y)-f(x))^2dydx\leq Cr^{n-1}\|\nabla f\|^2_{L^2(M)}+Cr^{n+1} \|f\|^2_{H^2(M)},$$
completing the proof of (\ref{asser2}). Now for (\ref{asser3}), applying (\ref{deltabnd}) once again but with $\delta=r$ this time around, 
$$\frac{1}{r^2}\int_M \int_{B_r(x)}\mathcal{Q}_{t(r)}(x,y)(f(y)-f(x))^2dydx\leq (1+r)\int_M\int_{B_r(x)} \mathcal{Q}_{t(r)}(x,y)L(x,y)^2dxdy$$$$+(1+r^{-1})\int_M\int_{B_r(x)}\mathcal{Q}_{t(r)}(x,y) R(x,y)^2dxdy$$
$$\leq (1+r)\int_M\int_{B_r(x)} \mathcal{Q}_{t(r)}(x,y)L(x,y)^2dxdy+C(1+r^{-1})\int_M\int_{B_r(x)}R(x,y)^2dxdy$$
$$\leq (1+r)\int_M\int_{B_r(x)} \mathcal{Q}_{t(r)}(x,y)L(x,y)^2dxdy+ Cr^{n+1}  \|f\|^2_{H^2(M)}.$$
Away from the boundary, we have that $\mathcal{Q}\equiv1$ and thus
$$\int_{M-N_r}\int_{B_r(x)} \mathcal{Q}_{t(r)}(x,y)L(x,y)^2dxdy$$
$$\leq (1+Cr^2)\frac{1}{r^2}\int_{M-N_r}\int_{|\xi|\leq r} \langle \nabla f (x),\xi\rangle^2 d\xi dx$$
$$=(1+Cr^2)\int_{M-N_r}\frac{r^n}{2c_\phi}|\nabla f(x)|^2$$
\begin{equation}\label{awayfrombndry}
\leq (1+Cr^2)\frac{ r^n}{2c_\phi}\|\nabla f \|^2_{L^2(M-N_r)}. 
\end{equation}
When $x\in N_r$ is near the boundary, we have that 
\begin{multline}
 \int_{B_r(x)} \mathcal{Q}_{t(r)}(x,y)L(x,y)^2dxdy\\ =r^{n-2}\int_{\{z:\psi_x(z)\in B_r(x)\}}\mathcal{Q}_t(x,\psi_x(z))\langle \nabla f(x), \text{exp}^{-1}_x(\psi_x(z))\rangle^2J_x(z)dz.
\end{multline}
 By Lemma \ref{scaledfermi},
 $$\text {exp}_x^{-1}(\psi_x(z))=r(z+E_x(z)),$$
 where $$|E_x(z)|\leq Cr.$$ 
 By Cauchy Schwarz, we have that 
 \begin{equation}\label{innprodestimate}
\langle \nabla f(x),z+E_x(z)\rangle ^2\leq (1+Cr)\langle \nabla f(x),z\rangle^2+Cr|\nabla f(x)|^2.
 \end{equation}
 It follows from (\ref{awayfrombndry}) and (\ref{innprodestimate}) that 
  $$\int_{B_r(x)} \mathcal{Q}_{t(r)}(x,y)L(x,y)^2dxdy\leq (1+Cr)r^{n}\int_{\{z:\psi_x(z)\in B_r(x)\}}\mathcal{Q}_t(x,\psi_x(z))\langle\nabla f(x),z\rangle^2J_x(z)dz$$
  $$+Cr^{n+1}|\nabla f (x)|^2\int_{\{z:\psi_x(z)\in B_r(x)\}}\mathcal{Q}_t(x,\psi_x(z))J_x(z)dz$$
  $$\leq (1+Cr)r^{n}\int_{\{z:\psi_x(z)\in B_r(x)\}}\mathcal{Q}_t(x,\psi_x(z))\langle\nabla f(x),z\rangle^2dz+Cr^{n+1}|\nabla f (x)|^2$$
  \begin{equation}\label{eq31}
\leq (1+Cr)r^{n}\int_{\{z:\psi_x(z)\in B_r(x)\}}\mathcal{Q}_t(z)\langle\nabla f(x),z\rangle^2dz+Cr^{n+1}|\nabla f (x)|^2.
  \end{equation}
Here, the last inequality follows from lemma \ref{klem}. Since 
$$|\{z:\psi_x(z)\in B_r(x)\}\Delta \mathbb{D}_t|\leq Cr,$$
and on the symmetric difference 
$$\mathcal{Q}_t(z)\leq C,\quad \langle \nabla f (x),z\rangle ^2\leq C |\nabla f(x)|^2 .$$
It follows that 
$$\int_{\{z:\psi_x(z)\in B_r(x)\}}\mathcal{Q}_t(z)\langle \nabla f (x),z\rangle^2 dz \leq \int_{\mathbb{D}_t}\mathcal{Q}_t(z)\langle \nabla f (x), z \rangle ^2dz + Cr|\nabla f(x)|^2.$$
By (\ref{Tayloridentity}), 
\begin{equation}\label{eq32}
\int_{\mathbb{D}_t}\mathcal{Q}_t(z)\langle \nabla f (x), z \rangle ^2dz\leq \frac{1}{2c_\phi}\Lambda(t)|\nabla f(x)|^2.
\end{equation}
Therefore, combining (\ref{eq31}) and (\ref{eq32}), we obtain
$$\int_{M}\int_{B_r(x)} \mathcal{Q}_{t(r)}(x,y)L(x,y)^2dxdy\leq \frac{r^n}{2c_\phi}(1+Cr)\int_M \Lambda(t)|\nabla f (x)|^2dx+Cr^{n+1} \| \nabla f \|^2_{L^2(M)}.$$
Hence
$$\frac{1}{r^2}\int_M \int_{B_r(x)}\mathcal{Q}_{t(r)}(x,y)(f(y)-f(x))^2dydx$$$$\leq \frac{r^{n}}{2 c_\phi}(1+Cr) \| \nabla f \|^2_{L^2(\Lambda)}+Cr^{n+1} \|f\|^2_{H^2(M)}.$$
\end{proof}
Note that the presence of the distorted Dirichlet energy near the collar is a direct consequence of the inexactitude of the Taylor identity (\ref{Tayloridentity}). We can now combine the two preceding lemmas to get a good bound on the discrete energy.
\begin{lemma}\label{energyPu}

For all $f\in C^2(\overline{M}),$ we have that 
$$\frac{2c_\phi}{\rho^{n+2}}\sum_{i,j:d(x_i,x_j)<\rho}\mathcal{Q}^{t(\rho)}_{ij}\int_{V_i}\int_{V_j}(f(x)-f(y))^2dxdy\leq \left(1+C\rho +C \frac{\epsilon}{\rho}\right)\|\nabla f\|^2_{L^2(\Lambda)}+C\left(\rho+\frac{\epsilon}{\rho}\right) \|f\|^2_{H^2(M)}$$
where 
$$\|\nabla f\|^2_{L^2(\Lambda)}:= \int_M \Lambda\left(\frac{d(x,\Sigma)}{\rho+2\epsilon}\right)|\nabla f(x)|^2dx.$$
It follows that 
$$E_\rho^D(Pf)\leq \left(1+C\rho+C \frac{\epsilon}{\rho}\right)\|\nabla f\|^2_{L^2(\Lambda)}+C\left(\rho+\frac{\epsilon}{\rho}\right) \|f\|^2_{H^2(M)}.$$
    
\end{lemma}
\begin{proof}
From inequality (\ref{eqerhotointerpol}),
$$E^D_\rho(Pf)\leq\frac{2c_\phi}{\rho^{n+2}}\sum_{i,j:d(x_i,x_j)<\rho}\mathcal{Q}^{t(\rho)}_{ij}\int_{V_i}\int_{V_j}(f(x)-f(y))^2dxdy$$
$$=\frac{2c_\phi}{\rho^{n+2}}\sum_{i,j:d(x_i,x_j)<\rho}\mathcal{Q}^{t(\rho)}_{ij}\int_{V_i}\int_{V_j}\mathbf{1}_{\{d(x,y)\leq 4\epsilon\}}(f(x)-f(y))^2dxdy$$$$
+\frac{2c_\phi}{\rho^{n+2}}\sum_{i,j:d(x_i,x_j)<\rho}\mathcal{Q}^{t(\rho)}_{ij}\int_{V_i}\int_{V_j}\mathbf{1}_{\{d(x,y)> 4\epsilon\}}(f(x)-f(y))^2dxdy.$$
For the first term, we have 
$$\frac{2c_\phi}{\rho^{n+2}}\sum_{i,j:d(x_i,x_j)<\rho}\mathcal{Q}^{t(\rho)}_{ij}\int_{V_i}\int_{V_j}\mathbf{1}_{\{d(x,y)\leq 4\epsilon\}}(f(x)-f(y))^2dxdy$$
$$\leq \frac{C}{\rho^{n+2}}\sum_{i,j}\int_{V_i}\int_{V_j}\mathbf{1}_{\{d(x,y)\leq 4\epsilon\}}(f(x)-f(y))^2dxdy$$
$$= \frac{C}{\rho^{n+2}}\int_{M}\int_{M\cap B_{4\epsilon(x)}}(f(x)-f(y))^2dxdy$$
$$= \frac{C \epsilon^{n+2}}{\rho^{n+2}\epsilon^{n+2}}\int_{M}\int_{B_{4\epsilon(x)}}(f(x)-f(y))^2dxdy$$

$$\leq C\frac{\epsilon}{\rho}\|\nabla f\|^2_{L^2(M)}+C\frac{\epsilon^2}{\rho^2}  \|f\|^2_{H^2(M)}. $$
Here, the last inequality follows from lemma \ref{directcomparison} with $r=4\epsilon$.
For the second term, we can apply lemma \ref{swapQdiscrtocont} with $R=\rho+2\epsilon$ in order to obtain
$$\frac{2c_\phi}{\rho^{n+2}}\sum_{i,j:d(x_i,x_j)<\rho}\mathcal{Q}^{t(\rho)}_{ij}\int_{V_i}\int_{V_j}\mathbf{1}_{\{d(x,y)> 4\epsilon\}}(f(x)-f(y))^2dxdy$$
$$\leq \frac{2c_\phi}{\rho^{n+2}}\sum_{i,j:d(x_i,x_j)<\rho}\int_{V_i}\int_{V_j}\mathcal{Q}_{t(\rho+2\epsilon)}(x,y)(f(x)-f(y))^2dxdy$$
$$+\left(\frac{\epsilon}{\rho}\right)\frac{C}{\rho^{n+2}}\sum_{i,j:d(x_i,x_j)<\rho}\int_{V_i}\int_{V_j}(f(x)-f(y))^2dxdy$$
$$+ \epsilon\frac{C}{\rho^{n+2}}\sum_{i,j:d(x_i,x_j)<\rho}\int_{V_i}\int_{V_j}\frac{1}{d(x,y)}(f(x)-f(y))^2dxdy$$
$$\leq \left(2c\phi +C\frac{\epsilon}{\rho}\right)\frac{1}{\rho^{n+2}}\sum_{i,j:d(x_i,x_j)<\rho}\int_{V_i}\int_{V_j}\mathcal{Q}_{t(\rho+2\epsilon)}(x,y)(f(x)-f(y))^2dxdy$$
$$+ \epsilon\frac{C}{\rho^{n+2}}\sum_{i,j:d(x_i,x_j)<\rho}\int_{V_i}\int_{V_j}\frac{1}{d(x,y)}(f(x)-f(y))^2dxdy$$
$$\leq \left(2c\phi +C\frac{\epsilon}{\rho}\right)\frac{1}{\rho^{n+2}}\int_{M}\int_{B_{\rho+2\epsilon}(x)}\mathcal{Q}_{t(\rho+2\epsilon)}(x,y)(f(x)-f(y))^2dxdy$$
$$+ \epsilon\frac{C}{\rho^{n+2}}\int_{M}\int_{B_{\rho+2\epsilon}(x)}\frac{1}{d(x,y)}(f(x)-f(y))^2dxdy$$
$$\leq \left(1+C\rho+C \frac{\epsilon}{\rho}\right)\|\nabla f\|^2_{L^2(\Lambda)}+C\left(\frac{\epsilon}{\rho}+\rho\right) \|f\|^2_{H^2(M)}.$$
Once again, the last inequality follows from lemma \ref{directcomparison}. 
\end{proof}
It follows immediately from lemma \ref{energyPu} and equation (\ref{eqerhotointerpol}) that \begin{equation}\label{energyctod}
E^D_\rho(Pf)\leq \left(1+C\rho+C \frac{\epsilon}{\rho} \right)\|\nabla f\|^2_{L^2(\Lambda)}+C\left( \rho+\frac{\epsilon}{\rho}\right) \|f\|^2_{H^2(M)}.
\end{equation}
\subsection{Quadrature bound on the discretized smooth functions}
Since the boundary norm appears in the denominator of the Rayleigh quotient, we must bound the $b$ norm of $T_\rho Pf$ below by the $L^2(\Sigma)$ norm of $f$. Similarly as in the energy case, we allow for some other norms to appear here as well, as they will be taken care of by elliptic regularity tools when the result is specialized to a linear combination of eigenfunctions. The operator $P$ defined in the beginning of this section is based on interior weights and interior Voronoi cells. This operator is thus not the most natural option for this quadrature step. We will therefore introduce $P_\Sigma$, a more natural operator for our context, which will work as a stepping stone between $Pf$ and $T_\rho Pf$ in the upcoming proofs. 
\\\\
Let 
$$P_\Sigma : H^1(M)\rightarrow L^2(B),$$
$$f\mapsto P_\Sigma f,$$
where 
$$(P_\Sigma f)_b:=\frac{1}{\mu_b}\int_{W_b}f(z)dz.$$
Then for all $f\in H^1(M),$
$$
    \|P_\Sigma f\|^2_B= \sum_b \mu_b (P_\Sigma f)^2_b=\sum_b\mu_b \left(\frac{1}{\mu_b}\int_{W_b}f(z)dz\right)^2$$\begin{equation}\label{PSigmacontract}
    \leq \sum_b \int_{W_b}f(z)^2dz = \|f|_\Sigma\|^2_{L^2(\Sigma)}.
\end{equation}
Let
$$P^*_\Sigma : L^2(B)\rightarrow L^2(\Sigma),$$
$$u\mapsto P^*_\Sigma u,$$
where
$$P^*_\Sigma u(x):= \sum_b \mathbf{1}_{W_b}(x) u_b.$$
Then, for all $f\in H^1(M)$,
\begin{equation}
    \|P_\Sigma f\|_B= \|P_\Sigma^*P_\Sigma f\|_{L^2(\Sigma)}.
\end{equation}
Furthermore,
\begin{equation}\label{L^2eq}
    \|f|_\Sigma\|^2_{L^2(\Sigma)}=\| P_\Sigma^* P_\Sigma f\|^2_{L^2(\Sigma)}+\|f|_\Sigma-P_\Sigma^*P_\Sigma f\|^2_{L^2(\Sigma)},
\end{equation}
and thus 
\begin{equation}\label{secondhalf}
    \|P_\Sigma f\|^2_B = \| f|_\Sigma \|^2_{L^2(\Sigma)}-\|f|_\Sigma - P^*_\Sigma P_\Sigma f\|^2_{L^2(\Sigma)}.\end{equation}
Let $\mathcal{O}:= I-P_\Sigma^* P_\Sigma$. Then, as an operator from $L^2(\Sigma)\mapsto L^2(\Sigma)$, by (\ref{L^2eq}), 
$$\|\mathcal{O}\|_{L^2(\Sigma)\rightarrow L^2(\Sigma)}\leq C.$$
As an operator from $H^1(\Sigma)$ to $L^2(\Sigma)$, the proof of lemma 4.2 of~\cite{BIK}, which only relies on the fact that
$$B^\Sigma_{c\epsilon}(x_b)\subset W_b \subset B^\Sigma_{C\epsilon}(x_b),$$
applied to the closed manifold $\Sigma$ implies
$$\|\mathcal{O}\|_{H^1(\Sigma)\rightarrow L^2(\Sigma)}\leq C\epsilon.$$
Then by interpolation,
$$\|\mathcal{O}\|_{H^{1/2}(\Sigma)\rightarrow L^2(\Sigma)}\leq \| \mathcal{O} \|^{1/2}_{L^2(\Sigma)\rightarrow L^2(\Sigma)}\|\mathcal{O}\|^{1/2}_{H^1(\Sigma)\rightarrow L^2(\Sigma)}\leq  C\epsilon^{1/2}, $$
and thus 
$$\|\mathcal{O} f\|^2_{L^2(\Sigma)}\leq C\epsilon\|f|_\Sigma \|^2_{H^{1/2}(\Sigma)}.$$
Using this last inequality in (\ref{secondhalf}) yields 
\begin{equation}\label{lastb}
     \|P_\Sigma f\|^2_B \geq  \| f|_\Sigma \|^2_{L^2(\Sigma)}-C\epsilon\|f|_\Sigma \|^2_{H^{1/2}(\Sigma)}.
\end{equation}
This is a promising bound. All we need going forward is to replace the $P_\Sigma f$  by the appropriate $T_\rho Pf$ trace term. The following lemma shows how to proceed. 
\begin{lemma}\label{first b}
    For all $f\in C^2(\overline{M}),$
    $$\|T_\rho(Pf)-P_\Sigma f\|^2_B\leq C\rho\|\nabla f\|^2_{L^2(M)} + C\rho \|f|_\Sigma \|^2_{H^{1/2}(\Sigma)}.$$
\end{lemma}
\begin{proof}
Recall that
$$(T_\rho Pf)_b:=\hat{\Theta}_\rho(b)^{-1}\sum_ik_\rho(b,x_i)\int_{V_i}f(y)dy.$$
Since
$$\hat{\Theta}_\rho(b)^{-1}\sum_i m_i k_\rho(b,x_i)=1,$$
and $$\frac{1}{\mu_b}\int_{W_b}1 dV_{n-1}=1,$$
we have that 
\begin{equation}\label{eq33}
(T_\rho Pf)_b-(P_\Sigma f)_b=\hat{\Theta}_\rho(b)^{-1}\frac{1}{\mu_b}\sum_i k_\rho(b,x_i)\int_{V_i}\int_{W_b}(f(y)-f(z))dzdy.
\end{equation}
Since 
$$\hat{\Theta}_\rho(b)^{-1}\frac{1}{\mu_b}\sum_i k_\rho(b,x_i)\int_{V_i}\int_{W_b} 1 dzdy=1,$$
applying Cauchy-Schwarz to (\ref{eq33}) yields

$$|(T_\rho Pf)_b-(P_\Sigma f)_b|^2\leq \hat{\Theta}_\rho(b)^{-1}\frac{1}{\mu_b}\sum_i k_\rho(b,x_i)\int_{V_i}\int_{W_b}|f(y)-f(z)|^2dzdy.$$
Now from the two kernel estimates 
$$\hat{\Theta}_\rho(b)\geq c>0, \quad k_\rho(b,x_i)\leq C\rho^{-n}\mathbf{1}_{\{d(b,x_i)<\rho\}},$$
it follows that
$$\mu_b|(T_\rho Pf)_b-(P_\Sigma f)_b|^2\leq \hat{\Theta}_\rho(b)^{-1}\sum_i k_\rho(b,x_i)\int_{V_i}\int_{W_b}|f(y)-f(z)|^2dzdy,$$
\begin{equation}\label{eq34}
\leq C\rho^{-n}\sum_{i: d(b,x_i)<\rho}\int_{V_i}\int_{W_b}|f(y)-f(z)|^2dzdy.
\end{equation}
Summing (\ref{eq34}) over $b$ yields 
\begin{equation}\label{eq35}
\|T_\rho(Pf)-P_\Sigma f\|^2_B \leq C\rho^{-n} \sum_b\sum_{i: d(b,x_i)<\rho}\int_{V_i}\int_{W_b}|f(y)-f(z)|^2dzdy. 
\end{equation}
If $d(b,x_i)<\rho, y\in V_i$ and $z\in W_b,$ then 
$$d(y,x_i)\leq \epsilon,\quad  d(z,b)\leq \epsilon.$$
Therefore 
$$d(y,z)\leq d(y,x_i)+d(x_i,b)+d(b,z)\leq \rho+2\epsilon<2\rho.$$
Thus the total integration region of (\ref{eq35}) is contained in 
$$\{(z,y)\in \Sigma\times M : d(y,z)<2\rho\}.$$
It implies that 
$$\|T_\rho(Pf)-P_\Sigma f\|^2_B \leq C\rho^{-n} \int_{\Sigma }\int_{B_{2\rho}(z)\subset M}|f(y)-f(z)|^2dydz. $$
The integral here takes entirely place in the collar 
$$N_{2\rho}:=\{x\in M :d(x,\Sigma)<2\rho \}.$$
Since $2\rho<r_*,$ use Fermi coordinates of Lemma \ref{collargeom}, 
$$y=\psi(p,t),\quad 0\leq t\leq 2\rho.$$
 It follows from Lemma \ref{collargeom} that
$$C^{-1}\leq J(p,t)\leq C.$$
For $y= \psi(p,t)\in B_{2\rho(z)},$ then 
$$t\leq 2\rho,\quad d_\Sigma(p,z)\leq C \rho.$$
The second inequality follows from Lemma \ref{collargeom} as well since 
$$d_\Sigma(p,z)\leq Cd(p,z)\leq C(d(p,y)+d(y,z))=C(t+d(y,z))< C\rho .$$
It follows that 
$$\rho^{-n} \int_{\Sigma }\int_{B_{C\rho}(z)\subset M}|f(y)-f(z)|^2dydz\leq  C\rho^{-n}\int_{\Sigma }\int_{B^\Sigma_{C\rho}(z)}\int_0^{2\rho}|f(\psi(p,t))-f(z)|^2dtdpdz$$
$$\leq C\rho^{-n}\int_{\Sigma }\int_{B^\Sigma_{C\rho}(z)}\int_0^{2\rho}|f(\psi(p,t))-f(p)|^2dtdpdz$$$$+C\rho^{-n}\int_{\Sigma }\int_{B^\Sigma_{C\rho}(z)}\int_0^{2\rho}|f(p)-f(z)|^2dtdpdz$$
\begin{equation}\label{splitnormtang}
=: I_N + I_T.    
\end{equation}
We bound each of these quantities separately. We begin with the normal term $I_N$. Let $f(p,t):= f(\psi(p,t))$. We have that 
$$|f(p,t)-f(p)|=\int_0^t \partial_s f(p,s)ds$$
and thus 
$$|f(p,t)-f(p)|^2\leq t \int_0^t |\partial_s f(p,s)|^2ds \leq C\rho \int_0^{2\rho}|\partial_s f(p,s)|^2ds.$$
It follows that 
$$I_N \leq C\rho^{-n}\int_{\Sigma}\int_{B^\Sigma_{C\rho}(z)}\rho^2\int_0^{2\rho}|\partial_sf(p,s)|^2dsdpdz.$$
Switching the integral in $p$ with the integral in $z$ yields 
$$I_N \leq C\rho\int_{\Sigma}\int_0^{2\rho}|\partial_sf(p,s)|^2dsdp$$
\begin{equation}\label{bndIN}
\leq C\rho\|\nabla f\|^2_{L^2(N_{2\rho})}\leq C\rho\|\nabla f \|^2_{L^2(M)}.
\end{equation}
We now bound the tangential term $I_T$. We have that 
$$I_T \leq  C\rho^{-n+1}\int_{\Sigma }\int_{B^\Sigma_{C\rho}(z)}|f(p)-f(z)|^2dpdz$$
$$=:\rho^{-n+1} C E_{C\rho}^\text{BIK}(f).$$
Here, $E_{C\rho}^\text{BIK}(f)$ is the interpolation energy of $f$ with respect to the closed manifold $\Sigma$ as described in~\cite{BIK}. The Remark 3.2 and the Lemma 3.3 of~\cite{BIK} give the following bounds
\begin{equation}\label{BIKgiven}
I_T\leq C \|f|_\Sigma \|^2_{L^2(\Sigma)}, \quad I_T\leq C\rho^2\|  f|_\Sigma \|^2_{H^1(\Sigma)}.
\end{equation}
For $g\in L^2(\Sigma)$, let 
$$\mathcal{O}_\rho(g)(z,p):= \rho^{-(n-1)/2}\mathbf{1}_{\{d_\Sigma(z,p)<C\rho\}}(g(p)-g(z)).$$
Then, by (\ref{BIKgiven}), 
$$\|\mathcal{O}_\rho g\|^2_{L^2(\Sigma\times \Sigma)}\leq C\|g\|^2_{L^2(\Sigma)}, \quad \|\mathcal{O}_\rho g\|^2_{L^2(\Sigma\times \Sigma)}\leq C\rho^2\|g\|^2_{H^1(\Sigma)}.$$
In other words, 
$$\mathcal{O}_\rho : L^2(\Sigma) \rightarrow L^2(\Sigma)^2,$$
has norm bounded by $C$ while 
$$\mathcal{O}_\rho: H^1(\Sigma)\rightarrow L^2(\Sigma)^2,$$
has norm bounded by $C\rho.$ Thus by the interpolation identity, 
$$\|\mathcal{O}_\rho \|_{H^{1/2}(\Sigma)\rightarrow L^2(\Sigma)^2}\leq \| \mathcal{O}_\rho \|^{1/2}_{L^2(\Sigma)\rightarrow L^2(\Sigma)^2}\|\mathcal{O}_\rho \|^{1/2}_{H^1(\Sigma)\rightarrow L^2(\Sigma)^2}\leq  C\rho^{1/2}. $$
It follows that 
\begin{equation}\label{bndIT}
I_T\leq C\rho \|f|_\Sigma \|^2_{H^{1/2}(\Sigma)}.
\end{equation}
Combining (\ref{splitnormtang}) with (\ref{bndIN}) and (\ref{bndIT}) concludes the proof.
\end{proof}
We can now use this last lemma to obtain the desired bound and achieve the goal of this subsection. We have that 
$$\| T_\rho(Pf)\|^2_B\geq (1-\sqrt{\rho})\|P_\Sigma f\|^2_B -\frac{1}{\sqrt{\rho}}\| T_\rho(Pf)-P_\Sigma f \|^2_B $$
$$\geq  (1-\sqrt{\rho})\|P_\Sigma f\|^2_B-C\sqrt{\rho}(\|\nabla f \|^2_{L^2(M)}+ \| f|_\Sigma \|^2_{H^{1/2}(\Sigma)})$$
\begin{equation}\label{quadctod}
\geq  (1-\sqrt{\rho})\| f\|_{L^2(\Sigma)}^2-C\sqrt{\rho}\|\nabla f \|^2_{L^2(M)}-C(\sqrt{\rho}+\epsilon) \| f|_\Sigma \|^2_{H^{1/2}(\Sigma)}.
\end{equation}
Here, the second inequality follows from lemma \ref{first b} and the third inequality follows from (\ref{lastb}).

\subsection{Proof of the upper bound on $\sigma_k(X,B)$ in terms of $\sigma_k(M)$}
In this subsection, we will use all previous results of this section in order to obtain our first eigenvalue bound. Let 
$$E_k(M):= \text{span}\{f_0,f_1,\dots,f_k\},$$
be the first $k+1$ Steklov eigenfunction span. Fix $f\in E_k(M)$ such that $f|_\Sigma\not=0$. It follows immediately that 
$$\|\nabla f\|^2_{L^2(M)}\leq \sigma_k(M)\|f\|^2_{L^2(\Sigma)}.$$
Furthermore, we can apply all regularity results (see section \ref{regularitytoolsubs}) to $f$. In particular,
$$\|f|_\Sigma\|^2_{H^{1/2}}\leq C(1+\sigma_k(M))\|f\|^2_{L^2(\Sigma)}.$$
Therefore, by (\ref{quadctod}), 
$$\|T_\rho(Pf)\|^2_B \geq (1-\sqrt{\rho)}\|f\|^2_{L^2(\Sigma)}-C\sqrt{\rho}\sigma_k(M)\|f\|^2_{L^2(\Sigma)}-C(\sqrt{\rho}+\epsilon)(1+\sigma_k(M))\|f\|^2_{L^2(\Sigma)}$$
\begin{equation}\label{drayleighquad}
    =: D_k(\rho,\epsilon)\|f\|^2_{L^2(\Sigma)}.
\end{equation}
Thus, for $\epsilon,\rho$ smaller than a constant depending only on the family of manifolds and the eigenvalue $\sigma_k(M),$ $T_\rho\circ P$ is injective and $T_\rho P(E_k(M))$ has dimension $k+1$ in $L^2(B)$. Thus $Pf$ is an admissible test function. We now need to estimate its discrete energy. Recall from section \ref{regularitytoolsubs} the following regularity bounds
$$ \|f\|^2_{H^2(M)}\leq C(1+\sigma_k(M))^3\|f\|^2_{L^2(\Sigma)}, \quad \|\nabla f\|^2_{L^2(N_{\rho+2\epsilon})}\leq C(\rho+2\epsilon)(1+\sigma_k)^{3}\|f|_\Sigma \|^2_{L^2(\Sigma)}.$$
    Since $\Lambda\left(\frac{d(x,\Sigma)}{\rho+2\epsilon}\right)\equiv 1$ on $M-N_{\rho+2\epsilon}$ and is bounded on $M$, we have that 
    $$\|\nabla f\|^2_{L^2(\Lambda)}:=\int_M \Lambda\left(\frac{d(x,\Sigma)}{\rho+2\epsilon}\right)|\nabla f(x)|^2dx\leq \|\nabla f\|^2_{L^2(M)}+C\|\nabla f\|^2_{L^2(N_{\rho+2\epsilon})}$$
    and thus  
    \begin{equation}\label{LambdaL2bound}
    \|\nabla f\|^2_{L^2(\Lambda)}\leq  \|\nabla f\|^2_{L^2(M)}+C(\rho+2\epsilon)(1+\sigma_k(M))^{3}\|f|_\Sigma \|^2_{L^2(\Sigma)}.
    \end{equation}
Substituting these bounds in (\ref{energyctod}) yields 
$$E_\rho^D(Pf)\leq \left(1+C\rho+C\frac{\epsilon}{\rho}\right)(\sigma_k(M)+C(\rho+\epsilon)(1+\sigma_k(M))^{3})\|f\|^2_{L^2(\Sigma)}$$
$$+ C\left(\rho+\frac{\epsilon}{\rho}\right) (1+\sigma_k(M))^3\|f\|^2_{L^2(\Sigma)}$$
\begin{equation}\label{drayleighenergy}
=: N_k(\rho,\epsilon)\|f\|^2_{L^2(\Sigma)}.  \end{equation}
It follows from the Rayleigh quotient for $\sigma_k(X,B)$ combined with equations (\ref{drayleighquad}) and (\ref{drayleighenergy}) that 
$$\sigma_k(X,B)\leq \max_{0\not=f\in E_k(M)}R^D(Pf)\leq \frac{N_k(\rho,\epsilon)}{D_k(\rho,\epsilon)}$$
$$= \frac{\left( 1+C\rho+C\frac{\epsilon}{\rho}\right) (\sigma_k(M)+C(\rho+\epsilon)(1+\sigma_k(M))^{3})+C\left(\rho+\frac{\epsilon}{\rho}\right)(1+\sigma_k(M))^3}{1-C(1+\sigma_k(M))(\sqrt{\rho+\epsilon)}}.$$
Note that this bound tends to $\sigma_k(M)$ as $\epsilon, \rho$ and $\frac{\epsilon}{\rho}$ go to $0$.

\begin{section}{ Upper bound on $\sigma_k(M)$ in terms of $\sigma_k(X,B)$ : From discrete to continuous}
The goal of this section is to bound $\sigma_k(M)$ in terms of $\sigma_k(X,B)$.  Similarly as in the previous section, we will proceed by using the variational characterization of the eigenvalue we want to bound. In this case, we will use the continuum Steklov variational characterization. In order to do so, we will once again need to construct test functions. This time around, the test functions will have to be in $H^1(M)$ and carry information about the graph's discrete Steklov eigenvalue. A natural choice is therefore to build these test functions from the graph's Steklov eigenfunctions. We will need to introduce a way to take a function $u$ from $L^2(X)$ and build a function on $H^1(M)$ in a way that uniformly conserves the appropriate amount of information. This part is generally trickier than the previous section, which relies mostly on direct Lipschitz quadrature bounds and elliptic regularity. This section will instead rely heavily on the finer design details of the discretization. The smoothing of discrete function will proceed in two separate steps. First, let 
$$P^* : L^2(X)\rightarrow L^2(M),$$
$$u\mapsto P^* u,$$
where 
$$P^* u := \sum_i \mathbf{1}_{V_i}u_i.$$
In general, $P^*u$ is not weakly differentiable and therefore we need to introduce a smoothing step. There are multiple options here. However, in order to obtain the bounds we need, we must find a way to smooth $P^*u$ in such a way that the Dirichlet energy  of the smoothed function is very close to the discrete energy $E_\rho^D(u)$ of the initial function. Define
$$k_r(x,y)=r^{-n}\phi\left(\frac{d(x,y)}{r}\right),\quad \Theta_r(x):=\int_M k_r(x,y)dy.$$
Then the smoother is given by 
$$\Lambda_r : L^2(M)\rightarrow H^1(M),$$
where
$$\Lambda_r(f):=\Theta_r(x)^{-1}\int_M k_r(x,y)f(y)dy.$$
Then, if $u\in L^2(X)$ is a graph function, its smoothed counterpart is given by 
$$I_\rho u := \Lambda_{\rho-4\epsilon}(P^*u).$$

\begin{subsection}{Preliminary smoother estimates}
    The focus of this subsection is to set the stage for the key lemma proved in the next subsection by proving basic estimates surrounding kernel specific functions. These estimates will be split between the two following lemmas.
        \begin{lemma}\label{interim}
            There exist constants $r_0>0$ and $C\geq 1$ depending only on the controlled geometry such that for all $r<r_0,$ all collar points $x$ such that $d(x,\Sigma)<r$, and all $z$ such that $|z|\leq 1+Cr, z_n\geq -t$ we have that
            \begin{equation}\label{ass1}
                |\Theta_r(x)-\Theta(t)|\leq Cr, \quad |\Theta_r(x)^{-1}-\Theta(t)^{-1}|\leq Cr, \quad |\Theta_r(x)^{-2}-\Theta(t)^{-2}|\leq Cr,
            \end{equation}
            \begin{equation}\label{ass2}
                |rd\Theta_r(x)-\Theta'(t)e_n|\leq Cr.
            \end{equation}
        \end{lemma}
        \begin{proof}
            We start by showing (\ref{ass1}). We use the scaled Fermi map $\psi_x$ of Lemma \ref{scaledfermi}. Let
            $$\Omega_x := \{z\in \mathbb{R}^n :\psi_x(z)\in B_r(x)\}.$$
            Because 
        $$\psi_x(z)\in M \iff z_n\geq -t,$$
        and by Lemma \ref{scaledfermi},
        $$(1-Cr)r|z|\leq d(x,\psi_x(z))\leq (1+Cr)r|z|,$$
        we have the inclusions
        $$\{z_n \geq -t, |z|\leq 1-Cr\}\subset \Omega_x\subset \{z_n \geq -t, |z|\leq 1+Cr\}.$$
        This implies that 
        $$\text{Vol}(\Omega_x\Delta \mathbb{D}_t)\leq Cr.$$
            Consider 
        $$\Omega :=\Omega_x\cap \mathbb{D}_t.$$
            Recall that 
            $$\Theta_r(x)=\int_{\Omega_x}\phi\left(\frac{d(x,\psi_x(z))}{r}\right)J_x(z)dz,\quad \Theta(t)= \int_{\mathbb{D}_t}\phi(|z|)dz.$$
            Therefore, 
            $$|\Theta_r(x)-\Theta(t)|\leq \int_\Omega \left|\phi\left(\frac{d(x,\psi_x(z))}{r}\right)J_x(z)-\phi(|z|)\right|dz$$
            $$+\int_{\Omega_x-\mathbb{D}_t}|\phi\left(\frac{d(x,\psi_x(z))}{r}\right)J_x(z)| + \int_{\mathbb{D}_t-\Omega_x}|\phi(|z|)|dz$$
            $$ \leq \int_\Omega \left|\phi\left(\frac{d(x,\psi_x(z))}{r}\right)J_x(z)-\phi(|z|)\right|dz+C\text{Vol}(\Omega_x\Delta \mathbb{D}_t)$$
            \begin{equation}\label{eq41}
            \leq \int_\Omega \left|\phi\left(\frac{d(x,\psi_x(z))}{r}\right)J_x(z)-\phi(|z|)\right|+Cr.\end{equation}
Here, the second inequality follows from the boundedness of both integrands as well as the inclusion of each domain in the symmetric difference. Furthermore, since $\phi$ is bounded and Lipschitz, 
$$\left|\phi\left(\frac{d(x,\psi_x(z))}{r}\right)J_x(z)-\phi(|z|)\right|$$$$\leq \left|\phi\left(\frac{d(x,\psi_x(z))}{r}\right)\right||J_x(z)-1|+\left|\phi\left(\frac{d(x,\psi_x(z))}{r}\right)-\phi(|z|)\right|$$
\begin{equation}\label{eq42}
\leq C|J_x(z)-1|+C\left|\frac{d(x,\psi_x(z))}{r}-|z| \right|\leq Cr.
\end{equation}
It follows from (\ref{eq41}) and (\ref{eq42}) that 
$$|\Theta_r(x)-\Theta(t)|\leq \int_\Omega Crdz+Cr\leq Cr.$$
        Since both $\Theta_r(x)$ and $\Theta(t)$ are bounded away from $0$, we also have that 
        $$|\Theta_r(x)^{-1}-\Theta(t)^{-1}|\leq Cr,\quad |\Theta_r(x)^{-2}-\Theta(t)^{-2}|\leq Cr.$$
        This completes the proof of (\ref{ass1}). We will now prove (\ref{ass2}).
        Direct computation yields 
        $$rd\Theta_r(x)=2c_\phi \int_{\Omega_x}\frac{\text{exp}_x^{-1}(\psi_x(z))}{r}J_x(z)dz,$$
        while 
        $$\Theta'(t)e_n=2c_\phi\int_{\mathbb{D}_t}zdz.$$
        We now want to compare the domain of these integrals. 
        By Lemma \ref{scaledfermi} and the Jacobian estimate
        $$\int_\Omega  \left|\frac{\text{exp}_x^{-1}(\psi_x(z))}{r}J_x(z)-z\right|dz$$$$\leq \int_\Omega  \left|\frac{\text{exp}_x^{-1}(\psi_x(z))}{r}-z\right||J_x(z)|dz+\int_\Omega |z||J_x(z)-1| 
        \leq Cr.$$
         Since both integrands are bounded on the symmetric difference of the domains $ \Omega_x$ and $\mathbb{D}_t$ we have that
         $$|rd\Theta_r(x)-\Theta'(t)e_n|$$$$\leq C\int_\Omega  \left|\frac{\text{exp}_x^{-1}(\psi_x(z))}{r}J_x(z)-z\right|dz+C\text{ Vol}(\{z:\psi_x(z)\in B_r(x)\}\Delta \mathbb{D}_t)
         \leq Cr.$$
        \end{proof}
        We will need to eventually bound the Dirichlet energy of smoothed discrete functions. In this process, we will have to deal with the derivatives of the kernel. The following lemma will help streamline this step. First, we need to introduce some notation. Let  $$\mathcal{B}_r(x,y):=\Theta_r(x)^{-1}d_xk_r(x,y)+d(\Theta_r(x)^{-1})k_r(x,y),$$ be the curved kernel derivative covector. Its natural rescaling in the Fermi coordinates is $r^{n+1}\mathcal{B}_r(x,\psi_x(z))$. Note that
        \begin{equation}\label{intbto0}
        \int_M \mathcal{B}_r(x,y)dy=d_x\left[\Theta_r(x)^{-1}\int_M k_r(x,y)dy\right]=d_x1=0.  
        \end{equation}
        Also, since $k_r(x,y)$ is supported where $d(x,y)<r,$
        \begin{equation}
            \mathcal{B}_r(x,y)=0\quad \text{when }d(x,y)\geq r
        \end{equation}
        and at fixed $r$,
        \begin{equation}\label{Bsup}
            |\mathcal{B}_r(x,y)|\leq Cr^{-n-1}\mathbf{1}_{\{d(x,y)<r\}}.
        \end{equation}
        Since 
        \begin{equation}\label{dlambd}
            d\Lambda_rf(x)=\int_Mf(y)B_r(x,y)dy,
        \end{equation}
        holds for smooth $f$, it follows from (\ref{Bsup}) and density that (\ref{dlambd}) holds weakly for $f\in L^2(M)$ as well. The next lemma aims to compare it with its flat counterpart, which is embedded nicely in the design of our discrete energy.
    \begin{lemma}\label{klem}
        There exist constants $r_0>0$ and $C\geq 1$ depending only on the controlled geometry such that for all $r<\min\{2\rho,r_0\} ,$ all collar points $x$ such that $d(x,\Sigma)<r$, and all $z$ such that $|z|\leq 1+Cr, z_n\geq -t$. Let 
        $$\mathcal{Q}_t(z):=a(t)\left(1-\frac{z_n^2}{|z|^2}\right)+b(t)\frac{z_n^2}{|z|^2}.$$
        Then
        \begin{equation}\label{firstasser}
            (1-Cr)\mathcal{Q}_t(z)\leq\mathcal{Q}_t(x,\psi_x(z))\leq(1+Cr)\mathcal{Q}_t(z),
        \end{equation}
        and for every unit covector $w\in T_x^*M$,
        \begin{equation}\label{secondasser}
       \int_{z_n\geq-t} \frac{(r^{n+1}\langle \mathcal{B}_r(x,\psi_x(z)),w\rangle)^2}{\mathcal{Q}_t(z)}dz\leq (1+r)\int_{z_n\geq-t} \frac{\langle L(t,z),w\rangle^2}{\mathcal{Q}_t(z)}dz+Cr.
        \end{equation}
        Here, $t$ and $z$ are given by the aforementioned change of coordinates.
    \end{lemma}

    \begin{proof}
        First, note that 
        $$d(\psi_x(z),\Sigma)-d(x,\Sigma)=d(x,\Sigma)+rz_n-d(x,\Sigma)= rz_n.$$
        Therefore 
        $$\eta_{xy}=\frac{d(\psi_x(z),\Sigma)-d(x,\Sigma)}{d(x,\psi_x(z))}=\frac{rz_n }{d(x,\psi_x(z))}.$$
        By the distance comparison of Lemma \ref{scaledfermi},
        $$(1-Cr)r|z|\leq d(x,\psi_x(z))\leq (1+Cr)r|z|.$$
        Therefore we have 
        $$|d(x,\psi_x(z))-r|z||\leq Cr^2|z|,$$
        and thus 
        $$\left|\eta_{x\psi_x(z)}-\frac{z_n}{|z|} \right|=\left|\frac{d(\psi_x(z),\Sigma)-d(x,\Sigma)}{d(x,\psi_x(z))}-\frac{z_n}{|z|} \right|=\left|\frac{rz_n}{d(x,\psi_x(z))}-\frac{z_n}{|z|} \right|$$
        $$=|z_n|\left|\frac{r|z|-d(x,\psi_x(z))}{d(x,\psi_x(z))|z|}\right|\leq |z_n|\frac{Cr^2}{d(x,\psi_x(z))}$$$$\leq Cr\frac{|z_n|}{|z|}\leq Cr.$$
        Then 
       $$|\mathcal{Q}_t(x,\psi_x(z))-\mathcal{Q}_t(z)|=|b(t)-a(t)|\left|\eta^2_{x\psi_x(z)}-\frac{z_n^2}{|z|^2}\right|$$ 
       \begin{equation}\label{Qdiffbnd}
        \leq Cr|b(t)-a(t)|\leq Cr.
       \end{equation}
       The last inequality follows from the fact that $a(t)$ and $b(t)$ are bounded functions. Since 
       $$\mathcal{Q}_t(z)\geq C_n,$$
       for a universal constant $C_n$, we have that
       $$Cr\leq Cr\mathcal{Q}_t(z).$$
       This last bound combined with (\ref{Qdiffbnd}) completes the proof of (\ref{firstasser}). We will now prove (\ref{secondasser}). By definition, we have that 
       $$r^{n+1}\mathcal{B}_r(x,\psi_x(z))=r^{n+1}\Theta_r(x)^{-1}d_xk_r(x,\psi_x(z))+r^{n+1}d(\Theta_r(x)^{-1})k_r(x,\psi_x(z))$$
       $$=2c_\phi\Theta_r(x)^{-1}\frac{\text{exp}_x^{-1}(\psi_x(z))}{r}\mathbf{1}_{\{\psi_x(z)\in B_r(x)\}}+rd(\Theta_r(x)^{-1})\phi\left(\frac{d(x,\psi_x(z))}{r}\right)$$
       $$=2c_\phi\Theta_r(x)^{-1}\frac{\text{exp}_x^{-1}(\psi_x(z))}{r}\mathbf{1}_{\{\psi_x(z)\in B_r(x)\}}-r\Theta_r(x)^{-2}d\Theta_r(x)\phi\left(\frac{d(x,\psi_x(z))}{r}\right).$$
    Recall that 
    $$L(t,z)=\mathbf{1}_{\mathbb{D}_t}(z)\left(\frac{2c_\phi }{\Theta(t)}z-\frac{\Theta'(t)}{\Theta(t)^2}\phi(|z|)e_n\right).$$
    We proceed by showing that
    \begin{equation}\label{alldir}
    \frac{\left| \Theta_r(x)^{-1}\frac{\text{exp}_x^{-1}(\psi_x(z))}{r}\mathbf{1}_{\{\psi_x(z)\in B_r(x)\}}-\Theta(t)^{-1}z\mathbf{1}_{\mathbb{D}_t}\right|^2}{\mathcal{Q}_t(z)}\leq Cr^2,   
    \end{equation}
    in $ \{z:\psi_x(z)\in B_r(x)\}\cap \mathbb{D}_t=:\Omega$
    and 
    \begin{equation}\label{normdir}
        \frac{\left|r\Theta_r(x)^{-2}d\Theta_r(x)\phi\left(\frac{d(x,\psi_x(z))}{r}\right)-\frac{\Theta'(t)}{\Theta(t)^2}\phi(|z|)e_n\right|^2}{\mathcal{Q}_t(z)}\leq Cr^2 \end{equation}
    We begin by showing (\ref{alldir}). Suppose $z\in \{z:\psi_x(z)\in B_r(x)\}\cap \mathbb{D}_t=:\Omega.$
The triangle inequality and Lemma \ref{scaledfermi} yield
$$\left|\Theta_r(x)^{-1}\frac{\text{exp}_x^{-1}(\psi_x(z))}{r}-\Theta(t)^{-1}z\right|\leq |\Theta_r(x)^{-1}|\left|\frac{\text{exp}_x^{-1}(\psi_x(z))}{r}-z\right|+|z||\Theta_r(x)^{-1}-\Theta(t)^{-1}|$$
$$\leq Cr|\Theta_r(x)^{-1}|+|z||\Theta_r(x)^{-1}-\Theta(t)^{-1}|$$
$$\leq Cr+|\Theta_r(x)^{-1}-\Theta(t)^{-1}|$$
$$\leq Cr.$$
Here, the last inequality follows from lemma \ref{interim}. Inequality (\ref{alldir}) follows directly from the uniform positive lower bound on $\mathcal{Q}_t(z).$ 
We will now show (\ref{normdir}).
        We start by using the triangle inequality followed by the boundedness of $\phi$,$\Theta'(t)$ and the boundedness of $\Theta$ away from $0$ to get 
        $$\left|r\Theta_r(x)^{-2}d\Theta_r(x)\phi\left(\frac{d(x,\psi_x(z))}{r}\right)-\frac{\Theta'(t)}{\Theta(t)^2}\phi(|z|)e_n\right|$$
        $$\leq \left|r\Theta_r(x)^{-2}d\Theta_r(x)-\frac{\Theta'(t)}{\Theta(t)^2}e_n\right|\phi\left(\frac{d(x,\psi_x(z))}{r}\right)+ \frac{|\Theta'(t)|}{\Theta(t)^2}\left|\phi\left(\frac{d(x,\psi_x(z))}{r}\right)-\phi(|z|)\right|$$
        $$\leq C\left|r\Theta_r(x)^{-2}d\Theta_r(x)-\frac{\Theta'(t)}{\Theta(t)^2}e_n\right|+ Cr\frac{|\Theta'(t)|}{\Theta(t)^2}$$
        $$\leq C|\Theta_r(x)^{-2}-\Theta(t)^{-2}||rd\Theta_r(x)|+\Theta(t)^{-2}|rd\Theta_r(x)-\Theta'(t)e_n|+ Cr\frac{|\Theta'(t)|}{\Theta(t)^2}$$
        $$\leq Cr+C|rd\Theta_r(x)-\Theta'(t)e_n|+ Cr$$
        $$=Cr+C|rd\Theta_r(x)-\Theta'(t)e_n|$$
        $$\leq Cr.$$
        Here, we used lemma \ref{interim} multiple times to bound the following quantities
        $$|rd\Theta_r(x)|\leq |\Theta'(t)|+Cr\leq C, $$
        $$|\Theta_r(x)^{-2}-\Theta(t)^{-2}|\leq Cr,$$
        $$|rd\Theta_r(x)-\Theta'(t)e_n|\leq Cr,$$
        in the fourth and fifth inequalities. Once again, (\ref{normdir}) follows directly from the lower bound on $\mathcal{Q}_t(z)$. Combining (\ref{alldir}) and (\ref{normdir}) we obtain the following
\begin{equation}\label{intest}
    \int_{\{z_n\geq -t\}\cap \Omega} \frac{(r^{n+1}\langle \mathcal{B}_r(x,\psi_x(z)),w\rangle-\langle L(t,z),w\rangle)^2}{\mathcal{Q}_t(z)}dz\leq Cr^2.
\end{equation}
 On the support of either of the two terms, if $z\not \in \Omega,$ then $z\in \{z:\psi_x(z)\in B_r(x)\}\Delta \mathbb{D}_t=:\Omega^c$. Note also that 
$$\left|\Theta_r(x)^{-1}\frac{\text{exp}_x^{-1}(\psi_x(z))}{r}\right|\leq C,$$
on $\{z:\psi_x(z)\in B_r(x)\}$ and that 
$$|\Theta(t)^{-1}z|\leq C,$$
on $\mathbb{D}_t$. Thus since $\mathcal{Q}_t(z)$ is positive and bounded away from zero, we have that
$$\int_{\{z_n\geq -t\}\cap \Omega^c}\frac{\left| \Theta_r(x)^{-1}\frac{\text{exp}_x^{-1}(\psi_x(z))}{r}\mathbf{1}_{\{\psi_x(z)\in B_r(x)\}}\right|^2}{\mathcal{Q}_t(z)}dz\leq C\text{ Vol}(\Omega^c)   \leq Cr,$$
and therefore 
$$\int_{\{z_n\geq -t\}\cap \Omega^c} \frac{(r^{n+1}\langle \mathcal{B}_r(x,\psi_x(z)),w\rangle)^2}{\mathcal{Q}_t(z)}dz\leq Cr.$$
It follows that
$$\int_{z_n\geq-t} \frac{(r^{n+1}\langle \mathcal{B}_r(x,\psi_x(z)),w\rangle)^2}{\mathcal{Q}_t(z)}dz\leq \int_{\{z_n\geq-t\}\cap\Omega} \frac{(r^{n+1}\langle \mathcal{B}_r(x,\psi_x(z)),w\rangle)^2}{\mathcal{Q}_t(z)}dz+Cr$$
$$=\int_{\{z_n\geq -t\}\cap \Omega} \frac{(r^{n+1}\langle \mathcal{B}_r(x,\psi_x(z)),w\rangle-\langle L(t,z),w\rangle+\langle L(t,z),w\rangle )^2}{\mathcal{Q}_t(z)}dz+Cr$$
$$\leq (1+r)\int_{\{z_n\geq -t\}\cap \Omega} \frac{\langle L(t,z),w\rangle^2}{\mathcal{Q}_t(z)}dz$$$$+(1+r^{-1})\int_{\{z_n\geq -t\}\cap \Omega} \frac{(r^{n+1}\langle \mathcal{B}_r(x,\psi_x(z)),w\rangle-\langle L(t,z),w\rangle)^2}{\mathcal{Q}_t(z)}dz+Cr$$
$$\leq (1+r)\int_{\{z_n\geq -t\}\cap \Omega} \frac{\langle L(t,z),w\rangle^2}{\mathcal{Q}_t(z)}dz+(1+r^{-1})Cr^2+Cr$$
$$= (1+r)\int_{\{z_n\geq -t\}\cap \Omega} \frac{\langle L(t,z),w\rangle^2}{\mathcal{Q}_t(z)}dz+Cr$$
$$\leq (1+r)\int_{z_n\geq-t} \frac{\langle L(t,z),w\rangle^2}{\mathcal{Q}_t(z)}dz+Cr.$$
Here, the first inequality follows from the $\Omega^c$ estimation above. The second to last inequality follows from (\ref{intest}).
    \end{proof}
\end{subsection}

\begin{subsection}{The key lemma}
The goal of this subsection is to prove the following lemma. 
\begin{lemma}\label{smoothedenergy}
    For every $f\in L^2(M)$, we have that
    $$\|d\Lambda_r f\|^2_{L^2(M)}\leq (1+Cr)\frac{2c_\phi}{r^{n+2}}\int_M\int_{B_r (x)}\mathcal{Q}_{t(r)}(x,y)(f(y)-f(x))^2dydx.$$
\end{lemma}
    \begin{proof}
        Let 
        $$F:=\Lambda_rf.$$
    By (\ref{dlambd}) and identity (\ref{intbto0}),
    $$dF(x)=\int_{B_r(x)}f(y)\mathcal{B}_r(x,y)dy=\int_{B_r(x)}(f(y)-f(x))\mathcal{B}_r(x,y)dy.$$
    If $d(x,\Sigma)\geq r,$ then the ball integral is not truncated by the boundary of $M$. Thus the Burago lemma 5.5 full proof applies. We automatically get that 
    $$|dF(x)|^2\leq (1+Cr)\frac{2c_\phi}{r^{n+2}}\int_{B_r(x)}(f(y)-f(x))^2dy$$
    $$=(1+Cr)\frac{2c_\phi}{r^{n+2}}\int_{B_r(x)}\mathcal{Q}_{t(r)}(x,y)(f(y)-f(x))^2dy.$$
    Thus we can focus on what happens near the boundary. Suppose $d(x,\Sigma)<r$. We have that 
    $$dF=\Theta_r(x)^{-1}\int_{B_r(x)}(f(y)-f(x))d_xk_r(x,y)dy+d(\Theta_r(x)^{-1})\int_{B_r(x)}k_r(x,y)(f(y)-f(x))$$
    $$=\int_{B_r(x)}(f(y)-f(x))\mathcal{B}_r(x,y)dy,$$
    where $\mathcal{B}_r(x,y)$ is the kernel derivative covector defined in the last subsection. Then for a unit covector $w\in T_x^*M$, we have that 
    $$\langle dF(x),w\rangle =\int_{B_r(x)}(f(y)-f(x))\langle \mathcal{B}_r(x,y),w\rangle dy.$$
    By Cauchy-Schwarz, we have that 
    $$\langle dF(x),w\rangle ^2\leq \left(\int_{B_r(x)}\mathcal{Q}_{t(r)}(x,y)(f(y)-f(x))^2 dy\right)\left(\int_{B_r(x)} \frac{\langle \mathcal{B}_r(x,y),w\rangle^2}{\mathcal{Q}_{t(r)}(x,y)}dy\right).$$
    Therefore, our goal is to prove 
    \begin{equation}\label{keygoal}
     \int_{B_r(x)} \frac{\langle \mathcal{B}_r(x,y),w\rangle^2}{\mathcal{Q}_{t(r)}(x,y)}dy\leq (1+Cr)\frac{2c_\phi}{r^{n+2}}.   
    \end{equation}
    By the change of coordinates presented in the last section followed by the bounds on its jacobian, we have that
    $$\int_{B_r(x)} \frac{\langle \mathcal{B}_r(x,y),w\rangle^2}{\mathcal{Q}_{t(r)}(x,y)}dy\leq (1+Cr)r^n\int_{z_n\geq -t(r)}\frac{\langle \mathcal{B}_r(x,\psi_x(z)),w\rangle^2}{\mathcal{Q}_{t(r)}(x,\psi_x(z))}dz.$$
    Using lemma \ref{klem} two times, applying the first assertion to the denominator then the second to the numerator, we obtain that
        $$\int_{B_r(x)} \frac{\langle \mathcal{B}_r(x,y),w\rangle^2}{\mathcal{Q}_{t(r)}(x,y)}dy\leq (1+Cr)r^n\int_{z_n\geq -t}\frac{\langle \mathcal{B}_r(x,\psi_x(z)),w\rangle^2}{\mathcal{Q}_t(z)}dz$$
        $$\leq (1+Cr)r^ {-n-2}\int_{z_n\geq -t}\frac{\langle r^{n+1} \mathcal{B}_r(x,\psi_x(z)),w\rangle^2}{\mathcal{Q}_t(z)}dz$$
        $$\leq (1+Cr)r^{-n-2}\int_{z_n\geq-t} \frac{\langle L(t,z),w\rangle^2}{\mathcal{Q}_t(z)}dz+Cr^{-n-1}$$
        $$\leq (1+Cr)r^{-n-2}2c_\phi+Cr^{-n-1}$$
     Here, the second to last inequality follows from lemma \ref{klem} while the last inequality follows from (\ref{dualid}). This completes the proof of (\ref{keygoal}) and therefore 
            $$|dF(x)|^2\leq (1+Cr)\frac{2c_\phi}{r^{n+2}}\int_{B_r(x)}\mathcal{Q}_{t(r)}(x,y)(f(y)-f(x))^2dy,$$
            for all $x\in M$. Integrating this last inequality over $M$ yields the desired result.
    \end{proof}
\end{subsection}

\begin{subsection}{Smoothed functions energy estimate}
The key lemma of the preceding subsection did most of the work in getting the energy estimates. To conclude the energy bound we only need to prove the following lemma.
\begin{lemma}
Let $r:= \rho-4\epsilon$. Then for all $u\in L^2(X),$ 
    $$ \frac{2c_\phi}{r^{n+2}}\int_M\int_{B_r(x)}\mathcal{Q}_{t(r)}(x,y)(P^*u(y)-P^*u(x))^2dy\leq \left(1+C\frac{\epsilon}{\rho}\right) E_\rho^D(u).$$
\end{lemma}
\begin{proof}
    If $x\in V_i, y\in V_j,$ then $P^*u(x)=u_i$ and $P^*u(y)=u_j.$ Thus 
$$\int_M\int_{B_r(x)}\mathcal{Q}_{t(r)}(x,y)(P^*u(y)-P^*u(x))^2dy=\sum_{i,j}\int_{V_i}\int_{V_j}\mathbf{1}_{\{d(x,y)<r\}}\mathcal{Q}_{t(r)}(x,y)(u_j-u_i)^2dydx.$$
Since $d(x,y)<r=\rho-4\epsilon$, then 
$$d(x_i,x_j)\leq d(x_i,x)+d(x,y)+d(y,x_j)\leq \epsilon+(\rho-4\epsilon)+\epsilon=\rho-2\epsilon<\rho.$$
Hence
$$\int_M\int_{B_r(x)}\mathcal{Q}_{t(r)}(x,y)(P^*u(y)-P^*u(x))^2dy=\sum_{d(x_i,x_j)<\rho}\int_{V_i}\int_{V_j}\mathbf{1}_{\{d(x,y)<r\}}\mathcal{Q}_{t(r)}(x,y)(u_j-u_i)^2dydx$$
\begin{equation}\label{Sdecomp}=S_S+S_L,
\end{equation}
where 
$$S_S:=\sum_{(i,j):d(x_i,x_j)<\rho}\int_{V_i}\int_{V_j}\mathbf{1}_{\{d(x,y)\leq 4\epsilon\}}\mathcal{Q}_{t(r)}(x,y)(u_j-u_i)^2dydx,$$
$$S_L:=\sum_{(i,j):d(x_i,x_j)<\rho}\int_{V_i}\int_{V_j}\mathbf{1}_{\{4\epsilon<d(x,y)< r\}}\mathcal{Q}_{t(r)}(x,y)(u_j-u_i)^2dydx.$$
We will bound each term separately, beginning with $S_S$.  Since $\mathcal{Q}_{t(r)}$ is bounded, we have 
$$S_S\leq C\sum_{(i,j):d(x_i,x_j)<\rho}\int_{V_i}\int_{V_j}\mathbf{1}_{\{d(x,y)\leq 4\epsilon\}}(u_j-u_i)^2dydx.$$
If $d(x,y)\leq 4\epsilon$, then
$$d(x_i,x_j)\leq d(x_i,x)+d(x,y)+d(y,x_j)\leq 6\epsilon,$$
which implies
$$ S_S\leq C \sum_{(i,j):d(x_i,x_j)<6\epsilon}\int_{V_i}\int_{V_j}(u_j-u_i)^2dydx$$
\begin{equation}\label{S_Sbnd}
=C\sum_{(i,j):d(x_i,x_j)<6\epsilon}m_im_j(u_j-u_i)^2=CE_{6\epsilon}(u).
\end{equation}
Thus applying  lemma \ref{sommepascomp} to (\ref{S_Sbnd}) yields
\begin{equation}\label{S_Sbnd1}
S_S\leq C\left(\frac{\epsilon}{\rho}\right)\sum_{(i,j):d(x_i,x_j)<\rho}m_im_j\mathcal{Q}_{ij}^{t(\rho)}(u_j-u_i)^2.
\end{equation}
We will now bound $S_L$. We begin by applying lemma \ref{swapQdiscrtocont} with $R=\rho-4\epsilon$ to $S_L$ and thus obtaining
\begin{equation}\label{S_Ldecomp}
S_L\leq S_0+S_1+S_2,
\end{equation}
where 
$$S_0:= \sum_{(i,j):d(x_i,x_j)<\rho}\int_{V_i}\int_{V_j}\mathbf{1}_{\{4\epsilon<d(x,y)< r\}}\mathcal{Q}_{ij}^{t(\rho)}(u_j-u_i)^2dydx,$$
$$S_1:=C \frac{\epsilon}{\rho} \sum_{(i,j):d(x_i,x_j)<\rho}\int_{V_i}\int_{V_j}(u_j-u_i)^2\mathbf{1}_{\{4\epsilon<d(x,y)< r\}}dydx,$$
$$S_2:= C\sum_{(i,j):d(x_i,x_j)<\rho}\int_{V_i}\int_{V_j}(u_j-u_i)^2\mathbf{1}_{\{4\epsilon<d(x,y)< r\}}\frac{\epsilon}{d(x,y)}dydx.$$
We will bound each term one by one. First, for $S_0$, simply bounding the indicator function by $1$ yields 
\begin{equation}\label{S_0bnd}
S_0\leq \sum_{(i,j):d(x_i,x_j)<\rho}\int_{V_i}\int_{V_j}\mathcal{Q}_{ij}^{t(\rho)}(u_j-u_i)^2dydx=\sum_{(i,j):d(x_i,x_j)<\rho}m_im_j\mathcal{Q}_{ij}^{t(\rho)}(u_j-u_i)^2.
\end{equation}
For $S_1$, we proceed in a similar manner, bounding the characteristic function by $1$ and introducing $\mathcal{Q}_{ij}^{t(\rho)}$. Note that this can be done because $\mathcal{Q}_{ij}^{t(\rho)}$ is bounded below away from $0$. We obtain 
\begin{equation}\label{S_1bnd}
S_1\leq C\frac{\epsilon}{\rho}\sum_{(i,j):d(x_i,x_j)<\rho}m_im_j\mathcal{Q}_{ij}^{t(\rho)}(u_j-u_i)^2.
\end{equation}
Finally, for $S_2$ we can apply the same bound on the indicator function followed by lemma \ref{complik} to obtain
\begin{equation}\label{S_2bnd}
S_2\leq C\frac{\epsilon}{\rho}\sum_{(i,j):d(x_i,x_j)<\rho}m_im_j\mathcal{Q}_{ij}^{t(\rho)}(u_j-u_i)^2.
\end{equation}
Therefore combining (\ref{S_0bnd}), (\ref{S_1bnd}) and (\ref{S_2bnd}) with (\ref{S_Ldecomp}) yields
\begin{equation}\label{S_Lbnd1}
<S_L\leq \left(1+C\frac{\epsilon}{\rho}\right)\sum_{(i,j):d(x_i,x_j)<\rho}m_im_j\mathcal{Q}_{ij}^{t(\rho)}(u_j-u_i)^2.
\end{equation}
Finally, combining (\ref{S_Sbnd1}) and (\ref{S_Lbnd1}) with (\ref{Sdecomp}) yields
$$\int_M\int_{B_r(x)}\mathcal{Q}_{t(r)}(x,y)(P^*u(y)-P^*u(x))^2dy\leq \left(1+C\frac{\epsilon}{\rho}\right)\sum_{(i,j):d(x_i,x_j)<\rho}m_im_j\mathcal{Q}_{ij}^{t(\rho)}(u_j-u_i)^2.$$
Since $r=\rho-4\epsilon$, we have 
$$r^{-(n+2)}\leq \left(1+C\frac{\epsilon}{\rho}\right) \rho^{-(n+2)}$$
and
$$\frac{2c_\phi}{r^{n+2}}\int_M\int_{B_r(x)}\mathcal{Q}_{t(r)}(x,y)(P^*u(y)-P^*u(x))^2dy\leq \left(1+C\frac{\epsilon}{\rho}\right)\frac{2c_\phi}{\rho^{n+2}}\sum_{(i,j):d(x_i,x_j)<\rho}m_im_j\mathcal{Q}_{ij}^{t(\rho)}(u_j-u_i)^2$$
$$=\left(1+C\frac{\epsilon}{\rho}\right)E_\rho^D(u).$$
\end{proof}
\end{subsection}
Combining this last lemma with lemma \ref{smoothedenergy}, we obtain the following result.
\begin{lemma}\label{energydtoc}
For all $u\in L^2(X),$
   $$\|\nabla  I_\rho u\|^2_{L^2(M)}\leq \left(1+C\left(\rho+\frac{\epsilon}{\rho}\right)\right)E_\rho^D(u).$$ 
\end{lemma}
  
    \begin{subsection}{Quadrature bound on the smoothed discrete functions}
        Fix $u\in L^2(X).$ In the Rayleigh quotient for the continuum Steklov problem, the $L^2(\Sigma)$ norm of the test function appears at the denominator. Thus if we want to use the variational characterization, we need to prove a lower bound on the $L^2(\Sigma)$ norm of the trace of the smoothed function $I_\rho( u)$ on $M$ in terms of the $L^2(B)$ norm of its discrete counterpart $T_\rho u$ and the discrete energy of $u$. Let 
        $$U:=P^*u,\quad F:= I_\rho u.$$
        Set $r:=\rho-4\epsilon$. Then by definition, 
        $$F(x)=\Theta_r(x)^{-1}\int_Mk_r(x,y)U(y)dy=\sum_i A_i(x) u_i$$
        where 
        $$A_i(x):= \Theta_r(x)^{-1}\int_{V_i}k_r(x,y)dy.$$
        By the definition of $\Theta_r(x),$ we immediately have that for all $x\in M,$
        \begin{equation}\label{sumAi=1}
            \sum_i A_i(x)=1.
        \end{equation}
        Similarly, we can write $T_\rho u$ as
        $$(T_\rho u)_b=\hat{\Theta}_\rho (b)^{-1}\sum_i m_ik_\rho (b,x_i)u_i=\sum_i \hat{A}_i(b)u_i$$
        where 
        $$\hat{A}_i(b):= \hat{\Theta}_\rho(b)^{-1}m_i k_\rho(b,x_i).$$
        Note that by the definition of $\hat{\Theta}_\rho (b)$, for all $x_b\in B$, we have that
        \begin{equation}\label{sumAichapo=1}
            \sum_i \hat{A}_i(b)=1.
        \end{equation}
        Let 
$$P_\Sigma : H^1(M)\rightarrow L^2(B),$$
$$f\mapsto P_\Sigma f,$$
where 
$$(P_\Sigma f)_b:=\frac{1}{\mu_b}\int_{W_b}f(z)dz.$$
Then for all $f\in H^1(M),$
$$
    \|P_\Sigma f\|^2_B= \sum_b \mu_b (P_\Sigma f)^2_b=\sum_b\mu_b \left(\frac{1}{\mu_b}\int_{W_b}f(z)dz\right)^2$$\begin{equation}\label{PSigmacontract}
    \leq \sum_b \int_{W_b}f(z)^2dz = \|f|_\Sigma \|^2_{L^2(\Sigma)}.
\end{equation}
Applying $P_\Sigma$ to $F$ yields
$$(P_\Sigma F)_b=\sum_i \bar{A}_i(b) u_i,$$
where 
$$\bar{A}_i(b):=\frac{1}{\mu_b}\int_{W_b}A_i(z)dz.$$
The boundary norm of this function will act as an intermediary between the boundary norm of $T_\rho u$ and the $L^2(\Sigma)$ norm of $I_\rho u$. Equation (\ref{sumAi=1}) implies
\begin{equation}
    \sum_i \bar{A}_i(b)=1.
\end{equation}
In order to streamline the upcoming proofs, we begin by proving some intermediary technical results. Let
    $$e_i(b):=\bar{A}_i(b)-\hat{A}_i(b).$$
    It follows directly that 
    $$(P_\Sigma F)_b-(T_\rho u)_b=\sum_i e_i(b) u_i.$$
    By (\ref{sumAi=1}) and (\ref{sumAichapo=1}), we obtain
    \begin{equation}\label{sumei=0}
        \sum_i e_i(b)=0.
    \end{equation}
    The next lemma introduces two important bounds related to this quantity. These identities will then be used in the following lemma to show the closeness between $P_\Sigma F$ and $T_\rho u$ in $L^2(B)$.
    \begin{lemma}\label{lemmbndonei}
        For all $x_b\in B,$
        \begin{equation}\label{sumabsoluteeibound}
            \sum_i |e_i(b)|\leq C\frac{\epsilon}{\rho},
        \end{equation}
        and
        \begin{equation}\label{pointwiseeibound}
            |e_i(b)|\leq C\frac{\epsilon}{\rho} m_i k_{2\rho}(x_b,x_i).
        \end{equation}
    \end{lemma}
    \begin{proof}
        For all $x\in M$ and all $b\in B,$
         \begin{equation}\label{Thetabbound}
        c\leq \Theta_r(x)\leq C,\quad c\leq \hat{\Theta}_r(b)\leq C,\quad  c\leq \hat{\Theta}_\rho(b)\leq C. \end{equation}
        Indeed, the bounds on $\hat{\Theta}_r(b)$ follow from 
        $$\hat{\Theta}_r(b)=\sum_i m_ik_r(x_b,x_i)\leq Cr^{-n} \sum_{i:d(x_b,x_i)<r}m_i\leq Cr^{-n}\text{Vol}(B_{r+2\epsilon}(x_b))\leq C,$$
        and 
        $$\hat{\Theta}_r(b)\geq cr^{-n}\sum_{i:d(b,x_i)<r/2}m_i\geq c r^{-n}\text{Vol}(B_{r/2}(x_b))\geq c.$$
        The bound on $\hat{\Theta}_\rho(b)$ is similar.
        We begin by estimating 
        $$|\frac{1}{\mu_b}\int_{W_b }\int_{V_i}k_r(z,y)dydz-m_i k_\rho (x_b,x_i)|.$$
        For $y\in V_i, z\in W_b$, we have $d(y,x_i)\leq \epsilon$ and $d(z,b)\leq \epsilon$. Thus 
        \begin{equation}\label{diffdist}
        |d(z,y)-d(x_b,x_i)|\leq C\epsilon,\end{equation}
        since 
        $$d(z,y)\leq d(z,x_b)+d(x_b,x_i)+d(x_i,y)= d(x_b,x_i)+2\epsilon,$$
        and 
        $$d(x_b,x_i)\leq d(x_b,z)+d(z,y)+d(y,x_i)=d(z,y)+2\epsilon.$$
        Also, 
        $$|k_r(z,y)-k_\rho(x_b,x_i)|=|r^{-n}\phi(d(z,y)/r)-\rho^{-n}\phi(d(x_b,x_i)/\rho)|$$
        $$\leq |r^{-n}\phi(d(z,y)/r)-r^{-n}\phi(d(x_b,x_i)/r)|+|r^{-n}\phi(d(x_b,x_i)/r)-\rho^{-n}\phi(d(x_b,x_i)/\rho)|.$$
        Since $\phi$ is $C$-Lipschitz, applying the Lipschitz inequality followed by (\ref{diffdist}) yields
        $$|r^{-n}\phi(d(z,y)/r)-r^{-n}\phi(d(x_b,x_i)/r)|\leq Cr^{-n-1}|d(z,y)-d(x_b,x_i)|\leq C\epsilon r^{-n-1}\leq  C\frac{\epsilon}{\rho}\rho^{-n}.$$
        Let
        $$G(s):= s^{-n}\phi(d(x_b,x_i)/s),$$
        for $s\in [r,\rho]$. Then either $d(x_b,x_i)<s$ or $G(s)=0$. If the former is true, then
        $$G'(s)=-ns^{-n-1}\phi(d(x_b,x_i)/s)-d(x_b,x_i)s^{-n-2}\phi'(d(x_b,x_i)/s)$$
        $$\leq Cs^{-n-1}+C d(x_b,x_i) s^{-n-2}$$
        $$\leq Cs^{-n-1}$$
        $$\leq C \rho^{-n-1}.$$
        Then, by the mean value theorem, 
        $$|r^{-n}\phi(d(x_b,x_i)/r)-\rho^{-n}\phi(d(x_b,x_i)/\rho)|=|G(r)-G(\rho)|$$
        $$\leq |r-\rho|\sup_{s\in [r,\rho]}|G'(s)|\leq C\frac{\epsilon}{\rho}\rho^{-n}.$$
        Therefore
        \begin{equation}\label{krbdirectbnd}
            |k_r(z,y)-k_\rho(x_b,x_i)|\leq  C\frac{\epsilon}{\rho}r^{-n},
        \end{equation}
        whenever at least one of the kernels is nonzero. 
        If either $k_r(z,y)\not=0$ for some $y\in V_i$ or $k_\rho(x_b,x_i)\not=0$ then by (\ref{diffdist}),
        $$d(x_b,x_i)\leq \rho.$$
        Hence 
        $$k_{2\rho}(x_b,x_i)=(2\rho )^{-n}\phi\left(\frac{d(x_b,x_i)}{2\rho}\right)\geq cr^{-n},$$
        on $B_\rho(x_b)$, which contains the support of $|k_r(z,y)-k_\rho(x_b,x_i)|$. It follows from (\ref{krbdirectbnd}) that 
        \begin{equation}\label{krbk2directbnd}
            |k_r(z,y)-k_\rho(x_b,x_i)|\leq  C\frac{\epsilon}{\rho}k_{2\rho}(x_b,x_i).
        \end{equation}
        Therefore
        $$\left|\frac{1}{\mu_b}\int_{W_b }\int_{V_i}k_r (z,y)dydz  - m_ik_\rho(x_b,x_i)\right|=\frac{1}{\mu_b}\left|\int_{W_b }\int_{V_i}(k_r(z,y) - k_\rho(x_b,x_i))dydz\right|.$$
        \begin{equation}\label{locestim}
            \leq C\frac{\epsilon}{\rho}m_ik_{2\rho}(x_b,x_i).
        \end{equation}
        We also have 
        \begin{equation}
           \sum_i\left|\frac{1}{\mu_b}\int_{W_b }\int_{V_i}k_r(z,y)dydz  - m_ik_\rho (x_b,x_i)  \right|\leq C\frac{\epsilon}{\rho}.
        \end{equation}
        Let 
        $$\bar{\Theta}_r(b):=\frac{1}{\mu_b}\int_{W_b }\Theta_r(z)dz.$$
        Then 
        $$|\bar{\Theta}_r(b)-\hat{\Theta}_\rho(b)|=\left|\sum_i \frac{1}{\mu_b}\int_{W_b }\int_{V_i}k_r(z,y)dydz-m_i k_\rho(x_b,x_i)\right|$$
        $$\leq \sum_i \left|\frac{1}{\mu_b}\int_{W_b }\int_{V_i}k_r(z,y)dydz-m_i k_\rho (x_b,x_i)\right|$$
        $$\leq C\frac{\epsilon}{\rho}.$$
        Since both $\bar{\Theta}_r$ and $\hat{\Theta}_\rho$ are bounded below away from $0$, 
        $$\left| \bar{\Theta}_r(b)^{-1}-\hat{\Theta}_\rho(b)^{-1}\right|\leq C\frac{\epsilon}{\rho}.$$
        If $z,z'\in W_b $, then $d(z,z')\leq 2\epsilon$. Since $\Theta_r$ is Lipschitz at scale $r$, 
        $$|\Theta_r(z)-\Theta_r(z')|\leq C\frac{\epsilon}{\rho},$$
        and thus 
        $$|\Theta_r(z)^{-1}-\bar{\Theta}_r(b)^{-1}|\leq C\frac{\epsilon}{\rho}.$$
        It follows that 
        $$\left|\bar{A}_i(b)-\bar{\Theta}_r(b)^{-1}\frac{1}{\mu_b}\int_{W_b }\int_{V_i}k_r(z,y)dydz\right|\leq \frac{1}{\mu_b}\int_{W_b} |\Theta_r(z)^{-1}-\bar{\Theta}_r(b)^{-1}|\int_{V_i}k_r(z,y)dydz$$
        $$\leq C\frac{\epsilon}{\rho}\frac{1}{\mu_b}\int_{W_b}\int_{V_i} k_r(z,y)dydz$$
        $$\leq C\frac{\epsilon}{\rho}\frac{1}{\mu_b}\int_{W_b}\int_{V_i} k_{2\rho}(x_b,x_i)dydz $$
        \begin{equation}\label{Aivsthetabar}
        =C\frac{\epsilon}{\rho}m_i k_{2\rho}(x_b,x_i).
        \end{equation}
        We also have that 
$$\left|\bar{\Theta}_r(b)^{-1}\frac{1}{\mu_b}\int_{W_b }\int_{V_i}k_r(z,y)dydz-\hat{A}_i(b) \right|  $$
        $$=\left|\bar{\Theta}_r(b)^{-1}\frac{1}{\mu_b}\int_{W_b }\int_{V_i}k_r(z,y)dydz-\hat{\Theta}_\rho^{-1}(b)m_ik_\rho(x_b,x_i) \right|  $$
        $$\leq \bar{\Theta}_r^{-1}(b)\left|\frac{1}{\mu_b}\int_{W_b }\int_{V_i}k_r(z,y)dydz-m_ik_\rho (x_b,x_i) \right| $$
        $$+|\bar{\Theta}_r^{-1}(b)-\hat{\Theta}_\rho^{-1}(b)|m_ik_\rho (x_b,x_i) $$
        \begin{equation}\label{Aivsthetachapo}
        \leq C\frac{\epsilon}{\rho}m_i k_{2\rho}(x_b,x_i).\end{equation}
 Hence, applying inequalities (\ref{Aivsthetabar}) and (\ref{Aivsthetachapo}) yields 
 $$|e_i(b)|=| \bar{A}_i(b)-\hat{A}_i(b)|$$
$$\leq |\bar{A}_i(b)-\bar{\Theta}_r(b)^{-1}\frac{1}{\mu_b}\int_{W_b }\int_{V_i}k_r(z,y)dydz|$$
$$+ \left|\bar{\Theta}_r(b)^{-1}\frac{1}{\mu_b}\int_{W_b }\int_{V_i}k_r(z,y)dydz-\hat{A}_i(b) \right|$$
$$\leq C\frac{\epsilon}{\rho}m_i k_{2\rho}(x_b,x_i).$$
        Furthermore, summing over $i$  yields 
        $$\sum_i |e_i(b)| \leq C\frac{\epsilon}{r}\sum_i m_i k_{2\rho}(x_b,x_i)\leq C \frac{\epsilon}{r}\leq C\frac{\epsilon}{\rho}.$$
    \end{proof}

        \begin{lemma}\label{PsigmaFclosetoTrhou}
        For all $u\in L^2(X)$,
$$\| P_\Sigma I_\rho u- T_\rho u\|^2_B=\sum_b \mu_b |(P_\Sigma I_\rho u)_b-(T_\rho u)|_b|^2 \leq C\left(\frac{\epsilon}{\rho}\right)^2\rho E_\rho^D(u).$$        \end{lemma}
        \begin{proof}
           Set $F:= I_\rho u$. First, recall that
            $$(P_\Sigma F)_b-(T_\rho u)_b =\sum_i e_i(b)u_i.$$
            By (\ref{sumei=0}), \begin{equation}\label{Fb-Trub}
            (P_\Sigma F)_b-(T_\rho u)_b=\sum_i e_i(b)u_i=\sum_i e_i(b) (u_i -(T_{2\rho} u)_b).            
            \end{equation}
            Then, by Cauchy-Schwarz followed by lemma \ref{lemmbndonei},
            $$|(P_\Sigma F)_b-(T_\rho u)_b|^2=\left|\sum_i e_i(b)(u_i-(T_{2\rho} u)_b)\right|^2$$
            $$\leq \left(\sum_i |e_i(b)|\right)\left(\sum_i |e_i(b)| |u_i-(T_{2\rho} u)_b|^2\right)$$
            $$\leq C \left(\frac{\epsilon}{\rho}\right)^2\sum_i m_i k_{2\rho}(x_b,x_i)|u_i-(T_{2\rho} u)_b|^2.$$
            Thus
            \begin{equation}\label{finito}
            \sum_b \mu_b |(P_\Sigma F)_b-(T_\rho u)|_b|^2 \leq C \left(\frac{\epsilon}{\rho}\right)^2\sum_b \mu_b\sum_i m_i k_{2\rho}(x_b,x_i)|u_i-(T_{2\rho} u)_b|^2.    
            \end{equation}
                Let's now focus on the inner sum. By definition 
                $$\sum_i m_i k_{2\rho}(x_b,x_i)|u_i-(T_{2\rho} u)_b|^2=\sum_i m_i k_{2\rho}(x_b,x_i)\left|u_i-\left(\hat{\Theta}_{2\rho}(b)^{-1}\sum_j m_jk_{2\rho}(x_b,x_j)u_j\right)\right|^2$$
                $$=\hat{\Theta}_{2\rho}(b)^{-1}\sum_i m_i k_{2\rho}(x_b,x_i)\left|\sum_j m_jk_{2\rho}(x_b,x_j)u_i-\left(\sum_j m_jk_{2\rho}(x_b,x_j)u_j\right)\right|^2$$
                $$\leq C\sum_{ij} m_i m_j k_{2\rho}(x_b,x_i)k_{2\rho}(x_b,x_j)\left|u_i-u_j\right|^2.$$
                It follows that
                $$\sum_b \mu_b\sum_i m_i k_{2\rho}(x_b,x_i)|u_i-(T_{2\rho} u)_b|^2$$$$\leq  C\sum_{ij} m_i m_j \left|u_i-u_j\right|^2\sum_b \mu_b k_{2\rho}(x_b,x_i)k_{2\rho}(x_b,x_j).$$
                If 
                $$k_{2\rho}(x_b,x_i)k_{2\rho}(x_b,x_j)\not=0,$$
                then $d(x_i,x_j)\leq 4\rho.$ Since 
                $$k_{2\rho}(x_b,x_i),k_{2\rho}(x_b,x_j)\leq C\rho^{-n},$$
                and the boundary points $x_b$ contributing to the inner sum lie in a boundary ball of radius $C\rho$ of volume at most $C\rho^{n-1}$, we have 
                $$\sum_b \mu_b k_{2\rho}(x_b,x_i)k_{2\rho}(x_b,x_j)\leq C \rho^{n-1}\rho^{-n}\rho^{-n} \mathbf{1}_{\{d(x_i,x_j)\leq C\rho\}},$$
                and 
                $$\sum_b \mu_b\sum_i m_i k_{2\rho}(x_b,x_i)|u_i-(T_{2\rho} u)_b|^2\leq  C\rho^{-n-1}\sum_{ij: d(x_i,x_j)<C\rho} m_i m_j \left|u_i-u_j\right|^2$$
                $$\leq C\rho^{-n-1}\sum_{ij: d(x_i,x_j)<\rho} m_i m_j \left|u_i-u_j\right|^2.$$
                Here, the last inequality follows from lemma \ref{rhohop}. Since $\mathcal{Q}_{ij}^{t(\rho),s}$ is bounded below, we have 
                $$\sum_b \mu_b\sum_i m_i k_{2\rho}(x_b,x_i)|u_i-(T_{2\rho} u)_b|^2\leq C\rho^{-n-1}\sum_{ij: d(x_i,x_j)<\rho} m_i m_j\mathcal{Q}_{ij}^{t(\rho),s} \left|u_i-u_j\right|^2$$
                $$\leq C\rho^{n+2}\rho^{-n-1} E_\rho^D(u)$$
                $$\leq C\rho E_\rho^D(u).$$
                Combining this last inequality with (\ref{finito}) yields 
                $$\sum_b \mu_b |(P_\Sigma F)_b-(T_\rho u)|_b|^2 \leq C\left(\frac{\epsilon}{\rho}\right)^2\rho E_\rho^D(u).$$
                
        \end{proof}
        We now prove the final bound of this subsection.
        \begin{lemma}\label{quadraturedtos}
        Let $u\in L^2(X).$ Then
        $$\|I_\rho u\|^2_{L^2(\Sigma)}\geq \left(1-\rho-\frac{\epsilon}{\rho}\right) \| T_\rho u\|^2_B - C \epsilon E_\rho^D(u).$$
        \end{lemma}
        \begin{proof}
        Set $F:=I_\rho u$. 
        We have that 
        $$\|P_\Sigma F\|^2_B\geq \left(1-\rho-\frac{\epsilon}{\rho}\right) \| T_\rho u\|^2_B - \left(\rho+\frac{\epsilon}{\rho}\right)^{-1}\|P_\Sigma F-T_\rho u\|^2_B.$$
        Applying lemma \ref{PsigmaFclosetoTrhou}, 
        $$\|P_\Sigma F\|^2_B\geq \left(1-\rho-\frac{\epsilon}{\rho}\right) \| T_\rho u\|^2_B -C \left(\rho+\frac{\epsilon}{\rho}\right)^{-1}\left(\frac{\epsilon}{\rho}\right)^2 \rho   E_\rho^D(u)$$
        $$\geq \left(1-\rho-\frac{\epsilon}{\rho}\right) \| T_\rho u\|^2_B - C \epsilon  E_\rho^D(u).$$
        Finally, by inequality (\ref{PSigmacontract}), 
        $$\|F\|^2_{L^2(\Sigma)}\geq \left(1-\rho-\frac{\epsilon}{\rho}\right) \| T_\rho u\|^2_B - C \epsilon E_\rho^D(u).$$
        \end{proof}
    \end{subsection}

\subsection{Upper bound on $\sigma_k(M)$ in terms of $\sigma_k(X,B)$}
In this subsection we will use all previous results of this section in order to obtain our second eigenvalue bound. Let 
$$E_k(X):= \text{span}\{u_0,u_1,\dots,u_k\},$$
be the first $k+1$ Steklov eigenfunction span. Fix $u\in E_k(X)$ such that $T_\rho u\not=0$. It follows that 
$$E_\rho^D(u)\leq \sigma_k(X,B) \|T_\rho u\|^2_B.$$
Then by lemma \ref{quadraturedtos},
\begin{equation}\label{crayleighquad}
    \|I_\rho u\|^2_{L^2(\Sigma)}\geq \left(1-\rho-\frac{\epsilon}{\rho}\right) \| T_\rho u\|^2_B - C \epsilon \sigma_k(X,B)  \|T_\rho u\|^2_B.\end{equation}
Thus, by the bound on $\sigma_k(X,B)$ in terms of $\sigma_k(M)$, for $\epsilon,\rho$ smaller than a constant depending only on the family of manifolds and the eigenvalue $\sigma_k(M),$ $I_\rho$ is injective and $I_\rho(E_k(X))$ has dimension $k+1$ in $L^2(\Sigma)$. Thus $I_\rho u$ is an admissible test function. We now need to estimate its Dirichlet energy. From lemma \ref{energydtoc}, 
  \begin{equation}\label{crayleighenergy}
  \|\nabla  I_\rho u\|^2_{L^2(M)}\leq \left(1+C\left(\rho+\frac{\epsilon}{\rho}\right)\right)\sigma_k(X,B) \|T_\rho u\|^2_B.\end{equation}
It follows from the Rayleigh quotient for $\sigma_k(M)$ combined with equations (\ref{crayleighquad}) and (\ref{crayleighenergy}) that 
$$\sigma_k(M)\leq \max_{0\not=u\in E_k(X)}R^D(I_\rho u)\leq \frac{\left(1+C\left(\rho+\frac{\epsilon}{\rho}\right)\right)\sigma_k(X,B)}{1-\rho-\frac{\epsilon}{\rho}- C\epsilon \sigma_k(X,B)}.$$

Note that this bound tends to $\sigma_k(X,B)$ as $\epsilon, \rho$ and $\frac{\epsilon}{\rho}$ go to $0$. Combining the final result of the preceding continuous-to-discrete comparison with this eigenvalue comparison, we immediately get Theorem 1.
    \end{section}

\begin{section}{Proof of Theorem 2}
The result is proved by comparing the natural discrete Steklov spectrum with the discrete generalized spectrum, for which we have Theorem 1 as a convergence result. First, let us recall some definitions. For $h\in L^2(X),$ we have that 
$$(T_\rho h)_b := \hat{\Theta}_\rho(b)^{-1} \sum_i m_i k_\rho(x_b,x_i)h_i,\quad \hat{\Theta}_\rho(b) :=\sum_i m_ik_\rho(x_b,x_i),$$
Since $T_\rho$ preserves constants we have that 
$$(T_\rho h)_b-h_b =\hat{\Theta}_\rho (b)^{-1} \sum_i m_i k_\rho(x_b,x_i)(h_i-h_b).$$
By the definition of $\hat{\Theta}_\rho$, applying Cauchy-Schwarz yields 
$$|(T_\rho h)_b-h_b|^2\leq \hat{\Theta}_\rho (b)^{-1} \sum_i m_i k_\rho(x_b,x_i)(h_i-h_b)^2.$$
Multiplying by $\mu_b$ and summing over $b$ yields
$$\|T_\rho h-h\|^2_B =\sum_b \mu_b |(T_\rho h)_b-h_b|^2 \leq \sum_b \mu_b \hat{\Theta}_\rho(b)^{-1} \sum_i m_i k_\rho(x_b,x_i)(h_i-h_b)^2$$
$$\leq C\rho^{-n} \sum_b \sum_{i:d(x_b,x_i)<\rho} \mu_b m_i (h_i-h_b)^2$$
$$\leq C\epsilon^{-1}\rho^{-n} \sum_b \sum_{i:d(x_b,x_i)<\rho} m_b m_i (h_i-h_b)^2.$$
Here, the second inequality follows from
$$\hat{\Theta}_\rho(b)\geq c >0, \quad k_\rho(x_b,x_i)\leq C\rho^{-n} \1_{\{d(x_b,x_i)<\rho\}},$$
while the third inequality follows from the fact that
$$\mu_b\leq C\frac{m_b}{\epsilon}.$$
This last inequality is exactly the structural obstruction forcing the additional $\rho^2/\ \epsilon\rightarrow 0$  assumption of Theorem 2. Since the discrete $Q$ weight is bounded below away from $0$, we have that 
\begin{equation}\label{nicesteklovkey}
\|T_\rho h-h\|^2_B\leq C\frac{\rho^2}{\epsilon}E^D_\rho(h).\end{equation}
This is the key identity that will allow us to prove Theorem 2. Let 
$$E_k^G(X):=\text{span}\{u_0,u_1,\dots u_k\},$$
be the first $k+1$ generalized Steklov eigenfunction span. Fix $u\in E_k^G(X)$ such that $T_\rho u \not=0$. It follows that 
$$E_\rho^D(u)\leq \sigma_k(X,B)\|T_\rho u\|^2_B.$$
Then, by (\ref{nicesteklovkey}), 
$$\|u\|_B\geq \|T_\rho u\|_B-\| T_\rho u - u\|_B$$
$$\geq \left(1-C\sqrt{\frac{\rho^2}{\epsilon}\sigma_k(X,B)}\right)\|T_\rho u\|_B.$$
By Theorem 1, $\sigma_k(X,B)\leq C_k$ where $C_k$ depends only on the controlled geometry and $\sigma_k(M)$. Hence if $\rho^2/\epsilon<c_k$ with $c_k$ chosen sufficiently small, then 
$$C\sqrt{\frac{\rho^2}{\epsilon}\sigma_k(X,B)}\leq \frac{1}{2},$$
and the restriction $R$ to the boundary is injective on $E_k^G(X)$ with $R(E_k^G(X))$ of dimension $k+1$. Thus $u$ can be taken as a test function for $\sigma^N_k(X,B)$. Doing so yields 
$$\sigma_k^N(X,B)\leq \max_{0\not= u\in E_k^G(X)}R^{D,\text{Nat}}(u)\leq \frac{\sigma_k(X,B)}{\left(1-C\sqrt{\frac{\rho^2}{\epsilon}\sigma_k(X,B)}\right)^2}.$$
In particular, $$\sigma_k^N(X,B)\leq 4\sigma_k(X,B)\leq C_k.$$ Now we do the other bound with the same strategy. 
Let 
$$E_k^N(X):=\text{span}\{v_0,v_1,\dots v_k\},$$
be the first $k+1$ natural Steklov eigenfunction span. Fix $v\in E_k^N(X)$ such that $v|_B \not=0$. It follows that 
$$E_\rho^D(v)\leq \sigma^N_k(X,B)\| v\|^2_B.$$
Then, by (\ref{nicesteklovkey}), 
$$\|T_\rho v\|_B\geq \| v\|_B-\| T_\rho v - v\|_B$$
$$\geq \left(1-C\sqrt{\frac{\rho^2}{\epsilon}\sigma^N_k(X,B)}\right)\| v\|_B.$$
Thus by the first comparison combined with  Theorem 1, for $\rho^2/\epsilon<c_k$ where $c_k$ is a constant depending only on controlled geometry and $\sigma_k(M)$, we have that
$$ C\sqrt{\frac{\rho^2}{\epsilon}\sigma^N_k(X,B)}\leq \frac{1}{2},$$ and thus $T_\rho$ is injective on $E_k^N(X)$ and $T_\rho(E_k^N(X))$ has dimension $k+1$. Thus $v$ can be taken as a test function for $\sigma_k(X,B)$. Doing so yields 
$$\sigma_k(X,B)\leq \max_{0\not=v\in E_k^N(X)}R^{D}(v)\leq \frac{\sigma^N_k(X,B)}{\left(1-C\sqrt{\frac{\rho^2}{\epsilon}\sigma^N_k(X,B)}\right)^2}.$$
Since in both Rayleigh computations the denominator factors are at least $1/2$ and both eigenvalues are uniformly bounded, the universal inequality
$$(1-a)^{-2}\leq 1+Ca$$
yields 
$$|\sigma_k^N(X,B)-\sigma_k(X,B)|\leq C_k \frac{\rho}{\sqrt{\epsilon}}.$$
The convergence result then follows from Theorem 1. 
\\\\
We now show that the additional hypothesis of Theorem 2 is indeed necessary. Suppose that $\rho,\epsilon,\epsilon/\rho\to0$, but that $\rho^2/\epsilon\not\rightarrow 0$. We show that $\sigma^N_1(X,B)$ does not converge to $\sigma_1(M)$ by building a well-suited test function. Note that uniformity of bounds is no longer a cause of concern for this part. Passing to a subsequence if necessary, we assume that there exists a constant $\delta$ such that
$$\frac{\rho^2}{\epsilon}\geq \delta>0.$$
Let $f$ be a first nonconstant Steklov eigenfunction, normalized by
$$\|f\|_{L^2(\Sigma)}=1,\quad \int_\Sigma f=0,$$
and let $q=Pf-c_\epsilon$, where the constant $c_\epsilon$ is chosen so that
$$g:=Rq,\quad \sum_{b\in B}\mu_b g_b=0.$$
Since $f$ is smooth, $g$ converges to $f|_\Sigma$ in the discrete boundary norm, and therefore
$$\|g\|_B^2=1+o(1),\quad E_\rho^D(q)\leq \sigma_1(M)+o(1).$$
Now let $w\in L^2(X)$ be supported on the boundary vertices and satisfy $w_b=g_b$ for $b\in B$. Using $m_b\leq C\epsilon\mu_b$, the boundedness of the $Q$-weights, and
$\sum_{d(x_b,x_j)<\rho}m_j\le C\rho^n$, we obtain
$$E_\rho^D(w) \leq C\frac{\epsilon}{\rho^2}\|g\|_B^2 \leq C_\delta\|g\|_B^2.$$
Moreover, since $q$ is obtained from a fixed smooth function, $q_i-q_j=O(\rho)$ on every $\rho$-edge. Crude counting bounds yield, for the bilinear form associated with $E_\rho^D$,
$$ |E_\rho^D(q,w)|\leq \frac{C}{\rho^{n+2}} \sum_b\sum_{i\not\in B: d(x_i,x_b)<\rho}m_im_b|q_i-q_b||g_b|+\frac{C}{\rho^{n+2}} \sum_b\sum_{i\in B: d(x_i,x_b)<\rho}m_im_b(q_i-q_b)^2$$
$$ \leq \frac{C_f}{\rho^{n+2}} \sum_b\sum_{i\not\in B: d(x_i,x_b)<\rho}m_im_b\rho+\frac{C_f}{\rho^{n+2}} \sum_b\sum_{i\in B: d(x_i,x_b)<\rho}m_im_b\rho^2$$
$$ \leq \frac{C_f}{\rho^{n+2}} \sum_b m_b\rho^{n+1}+\frac{C_f}{\rho^{n+2}} \sum_b m_b\rho^{n+1}$$
$$\leq \frac{C_f}{\rho^{n+2}} \sum_b \mu_b\epsilon\rho^{n+1}+\frac{C_f}{\rho^{n+2}} \sum_b \mu_b\epsilon\rho^{n+1}$$
$$=O\!\left(\frac{\epsilon}{\rho}\right)=o(1).$$
For $0\leq t\leq 1$, let
$$v_t=tq+(1-t)w.$$
Then $Rv_t=g$, and therefore $Rv_t$ is orthogonal to the constants. Therefore $\operatorname{span}\{1,v_t\}$ is an admissible two-dimensional space for $\sigma_1^N(X,B)$, and
$$
\sigma_1^N(X,B) \leq \frac{E_\rho^D(v_t)}{\|g\|_B^2} \leq t^2\sigma_1(M)+(1-t)^2C_\delta+o(1). $$
Since 
$$t^2\sigma_1(M)+(1-t)^2C_\delta < \sigma_1(M)\iff t>\frac{C_\delta-\sigma_1(M)}{C_\delta+\sigma_1(M)},$$
it suffices to fix $t<1$ sufficiently close to $1$ to obtain that 
$$\limsup \sigma_1^N(X,B) <\sigma_1(M).$$
Thus spectral convergence fails whenever $\rho^2/\epsilon$ does not tend to zero.
\end{section}

\section{Discretization of the Neumann problem}
As mentioned before, the discretization of the Neumann problem is obtained almost directly from the work done for Steklov. Indeed, we can reuse the energy estimates without issue. The only new tools needed are simpler quadrature bounds related to $L^2(M)$ and $L^2(X)$.
\subsection{Upper bound on $\mathbf{\lambda_k(X)}$ in terms of $\mathbf{\lambda_k(M)}$}
For $f\in C^2(M)$, since 
$$(Pf)_i=\frac{1}{m_i}\int_{V_i}f,$$
then 
\begin{equation}\label{eqq}
\|f-P^*Pf\|^2_{L^2(M)}=\sum_i\int_{V_i}|f-(Pf)_i|^2.
\end{equation}
Applying the bound
$$\int_{V_i}|f-(Pf)_i|^2=\frac{1}{2m_i}\int_{V_i}\int_{V_i}|f(x)-f(y)|^2dxdy\leq C\epsilon^{-n}\int_{V_i}\int_{V_i}|f(x)-f(y)|^2$$
to (\ref{eqq}) yields
$$\|f-P^*Pf\|^2_{L^2(M)}\leq C\epsilon^{-n}\sum_i\int_{V_i}\int_{V_i}|f(x)-f(y)|^2 $$
$$\leq C\epsilon^{-n}\int_{M}\int_{B_{2\epsilon}(x)}|f(x)-f(y)|^2.$$
The first assertion of Lemma \ref{directcomparison} with $r=2\epsilon$ yields 
$$\|f-P^*Pf\|^2_{L^2(M)} \leq C\epsilon^2 \|\nabla f \|^2_{L^2(M)}+C\epsilon^4 \|f\|^2_{H^2(M)}.$$
 It follows immediately that 
$$|\|Pf\|_X-\|f\|_{L^2(M)}|=|\|P^*Pf\|_{L^2(M)}-\|f\|_{L^2(M)}|\leq \|f-P^*Pf\|_{L^2(M)},$$
\begin{equation}\label{neumannquad}
    \leq C\epsilon\|\nabla f\|_{L^2(M)}+C\epsilon^{2} \|f\|_{H^2(M)}.\end{equation}
    Let $$E_k(M):=\text{span}\{f_0,f_1,f_2,\dots,f_k\},$$
    be the first $k+1$ Neumann eigenfunction span. Fix $f\in E_k(M)$. It follows that 
    $$\|\nabla f\|^2_{L^2(M)}\leq \lambda_k(M)\| f\|^2_{L^2(M)}.$$
    Note also that by (\ref{hessianboundneum}) and Lemma \ref{bndc1}, 
    $$ \|f\|^2_{H^2(M)}\leq C(1+\lambda_k(M))^2\|f\|^2_{L^2(M)}, \quad \|\nabla f\|^2_{L^2(N_{\rho+2\epsilon})}\leq C(\rho+2\epsilon)(1+\lambda_k(M))^{2}\|f \|^2_{L^2(M)}.$$
    Since $\Lambda\left(\frac{d(x,\Sigma)}{\rho+2\epsilon}\right)\equiv 1$ on $M-N_{\rho+2\epsilon}$ and is bounded on $M$, we have that
    $$\int_M \Lambda\left(\frac{d(x,\Sigma)}{\rho+2\epsilon}\right)|\nabla f(x)|^2dx=:\|\nabla f\|^2_{L^2(\Lambda)}\leq \|\nabla f\|^2_{L^2(M)}+C\|\nabla f\|^2_{L^2(N_{\rho+2\epsilon})}$$
    and thus  
    \begin{equation}\label{LambdaL2bound}
    \|\nabla f\|^2_{L^2(\Lambda)}\leq  \|\nabla f\|^2_{L^2(M)}+C(\rho+2\epsilon)(1+\lambda_k(M))^{2}\|f\|^2_{L^2(M)}.
    \end{equation}
Substituting these bounds in (\ref{neumannquad}) yields 
\begin{equation}\label{ndrayleighquad}
\|Pf\|_X^2\geq (1-C\epsilon(1+\lambda_k(M))^{2})\|f\|^2_{L^2(M)}=:D_k(\rho,\epsilon)\|f\|^2_{L^2(M)}\end{equation}
    Thus, for $\epsilon,\rho$ smaller than a constant depending only on the family of manifolds and the eigenvalue $\lambda_k(M),$ $P$ is injective and $P(E_k(M))$ has dimension $k+1$ in $L^2(X)$. Thus $Pf$ is an admissible test function. We now need to estimate its discrete energy. By (\ref{energyctod}) combined with the usual Garding's inequality trick, we have that
$$E_\rho^D(Pf)\leq \left(1+C\rho+C\frac{\epsilon}{\rho}\right)(\lambda_k(M)+C(\rho+\epsilon)(1+\lambda_k(M))^{2})\|f\|^2_{L^2(M)}$$
$$+ C\left(\rho+\frac{\epsilon}{\rho}\right) (1+\lambda_k(M))^2\|f\|^2_{L^2(M)}$$
\begin{equation}\label{ndrayleighenergy}
=: N_k(\rho,\epsilon)\|f\|^2_{L^2(M)}.  \end{equation}
It follows from the Rayleigh quotient for $\lambda_k(X)$ combined with equations (\ref{ndrayleighquad}) and (\ref{ndrayleighenergy}) that 
$$\lambda_k(X)\leq \max_{0\not=f \in E_k(M)}R_\text{Neu}^D(Pf)\leq \frac{N_k(\rho,\epsilon)}{D_k(\rho,\epsilon)}$$
$$= \frac{\left( 1+C\rho+C\frac{\epsilon}{\rho}\right) (\lambda_k(M)+C(\rho+\epsilon)(1+\lambda_k(M))^{2})+C\left(\rho+\frac{\epsilon}{\rho}\right)(1+\lambda_k(M))^2}{(1-C\epsilon(1+\lambda_k(M))^{2})}.$$
Note that this bound tends to $\lambda_k(M)$ as $\epsilon, \rho$ and $\frac{\epsilon}{\rho}$ go to $0$. 
\subsection{Upper bound on $\mathbf{\lambda_k(M)}$ in terms of $\mathbf{\lambda_k(X)}$}
By lemma 5.4 of~\cite{BIK}, we have that for every $f\in L^2(M),$
\begin{equation}
    \|\Lambda_{\rho-4\epsilon}f-f\|^2_{L^2(M)}\leq \frac{C}{(\rho-4\epsilon)^n}\int_M \int_{B_{\rho-4\epsilon}(x)}|f(y)-f(x)|^2dydx.
\end{equation}
Note that our definition of $\Lambda_r$ is different from that of~\cite{BIK}. The proof however can be copied verbatim in our case, as it only relies on the fact that $\Lambda_r$ preserves constants as well as some common bounds on $\Theta$ terms and $k_r$ terms. Once again, by the same argument as in lemma 6.2 of~\cite{BIK},
$$|\|I_\rho u\|_{L^2(M)}-\|u\|_X|\leq C\rho \left(\frac{1}{\rho^{n+2}} \sum_i \sum_{j:d(x_j,x_i)<\rho}  m_im_j (u_j-u_i)^2\right)^{1/2}.$$
Since the $\mathcal{Q}$ terms are bounded, 
\begin{equation}\label{nquadraturedtos}
|\|I_\rho u\|_{L^2(M)}-\|u\|_X|\leq C\rho E_\rho^D(u)^{1/2}.
\end{equation}
Let 
$$E_k(X):= \text{span}\{u_0,u_1,\dots,u_k\},$$
be the first $k+1$ Neumann eigenfunction span. Fix $u\in E_k(X)$. It follows that 
$$E_\rho^D(u)\leq \lambda_k(X) \| u\|^2_X.$$
Then by (\ref{nquadraturedtos}),
\begin{equation}\label{ncrayleighquad}
    \|I_\rho u\|_{L^2(M)}\geq (1-C\rho\sqrt{\lambda_k(X)})\|u\|_X\end{equation}
Thus, by the bound on $\lambda_k(X)$ in terms of $\lambda_k(M)$, for $\epsilon,\rho$ smaller than a constant depending only on the family of manifolds and the eigenvalue $\lambda_k(M),$ $I_\rho$ is injective and $I_\rho(E_k(X))$ has dimension $k+1$ in $H^1(M)$. Thus $I_\rho u$ is an admissible test function. We now need to estimate its Dirichlet energy. From lemma \ref{energydtoc}, 
  \begin{equation}\label{ncrayleighenergy}
  \|\nabla  I_\rho u\|^2_{L^2(M)}\leq \left(1+C\left(\rho+\frac{\epsilon}{\rho}\right)\right)\lambda_k(X) \| u\|^2_X.\end{equation}
It follows from the Rayleigh quotient for $\lambda_k(M)$ combined with equations (\ref{ncrayleighquad}) and (\ref{ncrayleighenergy}) that 
$$\lambda_k(M)\leq \max_{0\not=u\in E_k(X)}R^D_\text{Neu}(I_\rho u)\leq \frac{\left(1+C\left(\rho+\frac{\epsilon}{\rho}\right)\right)\lambda_k(X)}{(1-C\rho\sqrt{\lambda_k(X)})^2}.$$

Note that this bound tends to $\lambda_k(X)$ as $\epsilon, \rho$ and $\frac{\epsilon}{\rho}$ go to $0$. Combining this bound with the eigenvalue bound of the previous subsection yields Theorem 3. 
\section{AI Usage statement}
During the preparation of this manuscript, the author used ChatGPT to assist with language, exposition, and the review of some mathematical arguments. It was also used to explore the literature surrounding the needed regularity bounds and the minimal geometric requirements of those bounds to hold uniformly.  All mathematical claims, all proofs, and references in the final manuscript were constructed and verified by the author, who takes full responsibility for its content. 

\bibliographystyle{plain}
\bibliography{ref}

@article {BIK,
    AUTHOR = {Burago, Dmitri and Ivanov, Sergei and Kurylev, Yaroslav},
     TITLE = {A graph discretization of the {L}aplace-{B}eltrami operator},
   JOURNAL = {J. Spectr. Theory},
  FJOURNAL = {Journal of Spectral Theory},
    VOLUME = {4},
      YEAR = {2014},
    NUMBER = {4},
     PAGES = {675--714},
      ISSN = {1664-039X,1664-0403},
   MRCLASS = {58J50 (05C50 53C21 58J60 65N25)},
  MRNUMBER = {3299811},
MRREVIEWER = {Emil\ Saucan},
       DOI = {10.4171/JST/83},
       URL = {https://doi-org.acces.bibl.ulaval.ca/10.4171/JST/83},
}

@article {JL,
    AUTHOR = {Lu, Jinpeng},
     TITLE = {Graph approximations to the {L}aplacian spectra},
   JOURNAL = {J. Topol. Anal.},
  FJOURNAL = {Journal of Topology and Analysis},
    VOLUME = {14},
      YEAR = {2022},
    NUMBER = {1},
     PAGES = {111--145},
      ISSN = {1793-5253,1793-7167},
   MRCLASS = {58J50 (05C50 53C21 53C23 58J60 65N25)},
  MRNUMBER = {4411102},
MRREVIEWER = {Shiping\ Liu},
       DOI = {10.1142/S1793525320500442},
       URL = {https://doi-org.acces.bibl.ulaval.ca/10.1142/S1793525320500442},
}

@article {CGR,
    AUTHOR = {Colbois, Bruno and Girouard, Alexandre and Raveendran, Binoy},
     TITLE = {The {S}teklov spectrum and coarse discretizations of manifolds
              with boundary},
   JOURNAL = {Pure Appl. Math. Q.},
  FJOURNAL = {Pure and Applied Mathematics Quarterly},
    VOLUME = {14},
      YEAR = {2018},
    NUMBER = {2},
     PAGES = {357--392},
      ISSN = {1558-8599,1558-8602},
   MRCLASS = {53C21 (58C40)},
  MRNUMBER = {4047402},
MRREVIEWER = {Changwei\ Xiong},
       DOI = {10.4310/pamq.2018.v14.n2.a3},
       URL = {https://doi-org.acces.bibl.ulaval.ca/10.4310/pamq.2018.v14.n2.a3},
}

@misc{Aubry2013,
      title={Approximation of the spectrum of a manifold by discretization}, 
      author={Erwann Aubry},
      year={2013},
      eprint={1301.3663},
      archivePrefix={arXiv},
      primaryClass={math.AP},
      url={https://arxiv.org/abs/1301.3663}, 
}

@article {Dodziuk,
    AUTHOR = {Dodziuk, J. and Patodi, V. K.},
     TITLE = {Riemannian structures and triangulations of manifolds},
   JOURNAL = {J. Indian Math. Soc. (N.S.)},
  FJOURNAL = {The Journal of the Indian Mathematical Society. New Series},
    VOLUME = {40},
      YEAR = {1976},
    NUMBER = {1-4},
     PAGES = {1--52},
      ISSN = {0019-5839,2455-6475},
   MRCLASS = {58G10 (57D20)},
  MRNUMBER = {488179},
MRREVIEWER = {M.\ L.\ Gromov},
}

@article {DT,
    AUTHOR = {Dias, Manuel and Tewodrose, David},
     TITLE = {Spectral properties of symmetrized {AMV} operators},
   JOURNAL = {J. Spectr. Theory},
  FJOURNAL = {Journal of Spectral Theory},
    VOLUME = {16},
      YEAR = {2026},
    NUMBER = {1},
     PAGES = {93--143},
      ISSN = {1664-039X,1664-0403},
   MRCLASS = {58J50 (35J05 35P05)},
  MRNUMBER = {5032301},
       DOI = {10.4171/jst/587},
       URL = {https://doi-org.acces.bibl.ulaval.ca/10.4171/jst/587},
}

@article{Guijarro,
  title={Restrictions on submanifolds via focal radius bounds},
  author={Luis Guijarro and Frederick Wilhelm},
  journal={Mathematical Research Letters},
  year={2020},
  url={https://api.semanticscholar.org/CorpusID:119132465}
}

@article{Anderson1990,
author = {Anderson, Michael T.},
journal = {Inventiones mathematicae},
number = {2},
pages = {429-446},
title = {Convergence and rigidity of manifolds under Ricci curvature bounds.},
url = {http://eudml.org/doc/143841},
volume = {102},
year = {1990},
}

@book{Stein,
 ISBN = {9780691080796},
 author = {ELIAS M. STEIN},
 publisher = {Princeton University Press},
 title = {Singular Integrals and Differentiability Properties of Functions (PMS-30)},
 year = {1970}
}

@book{Tubes,
 ISBN = {9783764369071},
 author = {ALFRED GRAY},
 publisher = {Birkhäuser Basel},
 title = {Tubes},
 year = {1990}
}

@book{Chavel, place={Cambridge},  title={Riemannian Geometry: A Modern Introduction}, publisher={Cambridge University Press}, author={Chavel, Isaac}, year={2006}}

@book{Griswald,
author = {Grisvard, Pierre},
title = {Elliptic Problems in Nonsmooth Domains},
publisher = {Society for Industrial and Applied Mathematics},
year = {2011},
doi = {10.1137/1.9781611972030}
}

@misc{Lyons,
      title={Spectral Convergence of the Discrete Steklov Problem}, 
      author={Andrew Lyons},
      year={2026},
      eprint={2610.40920},
      archivePrefix={arXiv},
      url={https://arxiv.org/abs/2610.40920}, 
}

@inproceedings{ColinDeVerdiere1986,
  author    = {Colin de Verdi{\`e}re, Yves},
  title     = {Spectre de vari{\'e}t{\'e}s riemanniennes et spectre de graphes},
  booktitle = {Proceedings of the International Congress of Mathematicians, Berkeley, 1986},
  pages     = {522--530},
  year      = {1987},
  publisher = {American Mathematical Society},
  note      = {Congress held in 1986; often cited as ICM (Berkeley), 1986}
}

@article{Mantuano2005,
  author  = {Mantuano, Tatiana},
  title   = {Discretization of compact {Riemannian} manifolds applied to the spectrum of {Laplacian}},
  journal = {Annals of Global Analysis and Geometry},
  volume  = {27},
  pages   = {33--46},
  year    = {2005},
  doi     = {10.1007/s10455-005-5215-0}
}

@article{HuaHuangWang2017,
  author  = {Hua, Bobo and Huang, Yan and Wang, Zuoqin},
  title   = {First eigenvalue estimates of {Dirichlet-to-Neumann} operators on graphs},
  journal = {Calculus of Variations and Partial Differential Equations},
  volume  = {56},
  number  = {6},
  pages   = {Article 178},
  year    = {2017},
  doi     = {10.1007/s00526-017-1260-3}
}

@article{Perrin2019,
  author  = {Perrin, H{\'e}l{\`e}ne},
  title   = {Lower bounds for the first eigenvalue of the {Steklov} problem on graphs},
  journal = {Calculus of Variations and Partial Differential Equations},
  volume  = {58},
  number  = {2},
  pages   = {Article 67},
  year    = {2019},
  doi     = {10.1007/s00526-019-1516-1}
}

@article{Perrin2021,
  author  = {Perrin, H{\'e}l{\`e}ne},
  title   = {Isoperimetric Upper Bound for the First Eigenvalue of Discrete {Steklov} Problems},
  journal = {Journal of Geometric Analysis},
  volume  = {31},
  number  = {8},
  pages   = {8144--8155},
  year    = {2021},
  doi     = {10.1007/s12220-020-00572-2}
}

@article{Tschanz2022,
  author  = {Tschanz, L{\'e}onard},
  title   = {Upper bounds for {Steklov} eigenvalues of subgraphs of polynomial growth {Cayley} graphs},
  journal = {Annals of Global Analysis and Geometry},
  volume  = {61},
  number  = {1},
  pages   = {37--55},
  year    = {2022},
  doi     = {10.1007/s10455-021-09799-w}
}

@article{Tschanz2023,
  author  = {Tschanz, L{\'e}onard},
  title   = {The {Steklov} Problem on Triangle-Tiling Graphs in the Hyperbolic Plane},
  journal = {Journal of Geometric Analysis},
  volume  = {33},
  pages   = {Article 161},
  year    = {2023},
  doi     = {10.1007/s12220-023-01208-x}
}

@article{VaughnBerryAntil2024,
  author  = {Vaughn, Ryan and Berry, Tyrus and Antil, Harbir},
  title   = {Diffusion maps for embedded manifolds with boundary with applications to {PDEs}},
  journal = {Applied and Computational Harmonic Analysis},
  volume  = {68},
  pages   = {Article 101593},
  year    = {2024},
  doi     = {10.1016/j.acha.2023.101593}
}

@article{PeoplesHarlim2026,
  author  = {Peoples, J. Wilson and Harlim, John},
  title   = {Spectral Convergence of Symmetrized Graph Laplacian on Manifolds with Boundary},
  journal = {Foundations of Data Science},
  volume  = {8},
  pages   = {119--167},
  year    = {2026},
  doi     = {10.3934/fods.2025005}
}

@article{AndreevTodorov2004,
  author  = {Andreev, Andrey B. and Todorov, Todor D.},
  title   = {Isoparametric finite-element approximation of a {Steklov} eigenvalue problem},
  journal = {IMA Journal of Numerical Analysis},
  volume  = {24},
  number  = {2},
  pages   = {309--322},
  year    = {2004},
  doi     = {10.1093/imanum/24.2.309}
}

@article{ArmentanoPadra2008,
  author  = {Armentano, Mar{\'i}a G. and Padra, Claudio},
  title   = {A posteriori error estimates for the {Steklov} eigenvalue problem},
  journal = {Applied Numerical Mathematics},
  volume  = {58},
  number  = {5},
  pages   = {593--601},
  year    = {2008},
  doi     = {10.1016/j.apnum.2007.01.011}
}

@article{Bogosel2016,
  author  = {Bogosel, Beniamin},
  title   = {The method of fundamental solutions applied to boundary eigenvalue problems},
  journal = {Journal of Computational and Applied Mathematics},
  volume  = {306},
  pages   = {265--285},
  year    = {2016},
  doi     = {10.1016/j.cam.2016.04.008}
}

@article{WangBenChenPolterovichSolomon2018,
  author  = {Wang, Yu and Ben-Chen, Mirela and Polterovich, Iosif and Solomon, Justin},
  title   = {{Steklov} Spectral Geometry for Extrinsic Shape Analysis},
  journal = {ACM Transactions on Graphics},
  volume  = {38},
  number  = {1},
  pages   = {Article 7},
  year    = {2018},
  doi     = {10.1145/3152156}
}
\end{document}